\documentclass{llncs}
\usepackage{epsfig}
\usepackage{amsmath}
\usepackage{amsfonts}
\usepackage{amssymb}

\newcommand{\R}{{I\!\!R}}
\newcommand{\N}{{I\!\!N}}

\def\R{{\rm I}\! {\rm R}}

\newtheorem{algorithm}[theorem]{Algorithm}

\newcommand{\OT}{\mathcal O}

\begin{document}

\pagestyle{headings}

\title{Splitting Methods For Solving Multi-Component Transport Model: A Multicomponent Mixture for Hydrogen Plasma}
\author{J\"urgen Geiser}
\institute{Ruhr University of Bochum, \\
The Institute of Theoretical Electrical Engineering, \\
Universit\"atsstrasse 150, D-44801 Bochum, Germany \\
\email{juergen.geiser@ruhr-uni-bochum.de}}
\maketitle

\begin{abstract}

In this paper, we present a splitting algorithm to 
solve multicomponent transport models.
These models are related to plasma simulations, in which
we consider the local thermodynamic equilibrium and 
weakly ionised plasma-mixture
models that are used for medical and technical applications in etching processes.
These multi-component transport-mixture models can be derived by 
approximately solving a linearised multi-component
Boltzmann equation with an approximation of the collision terms in the
mass, momentum and energy equations.
The model-equations are nonlinear partial differential equations
and they are known as Stefan-Maxwell equations. However, these partial differential
equations are delicate to solve and we propose to use noniterative and iterative splitting methods.
In the numerical experiments, we see the benefit of the iterative
splitting methods, while these methods can relax the nonlinear terms.

\end{abstract}

{\bf Keywords}:  multi-component model, Boltzmann equation, Chapman-Enskog expansion, Stefan-Maxwell equations, splitting methods, iterative splitting methods. \\

{\bf AMS subject classifications.} 35K25, 35K20, 74S10, 70G65.

\section{Introduction}

Understanding normal pressure, room temperature plasma applications
is important because of their use in medical and technical processes.
The increasing importance of plasma chemistry based on
multi-component plasma is a key factor for this trend, for 
low pressure plasma see \cite{sene06} and for atmospheric 
pressure regimes see \cite{tanaka2004}.
Both the influence of the mass transfer in the multi-component
mixture and the standard conservation laws have to be improved.
Although these improvements are well-known in fusion research---see,
for example, the modelling of high ionised plasmas \cite{igit2011}---, 
only a little 
work has been done for a weak-ionised plasma in atmospheric pressure regimes.

In this paper, we concentrate on an extension of the multicomponent transport model with 
respect to the reaction terms, see \cite{sene06}, while we
can approximate the collision integrals. The diffusive velocity
is simulated by the Stefan-Maxwell problem transport algorithm, see \cite{geiser2016}.
Based on the nonlinear diffusion term, we have to apply numerical schemes that
can solve the Stefan-Maxwell problem. We propose iterative schemes in combination 
with splitting approaches, which means that we decompose the transport- and 
reaction-parts, see \cite{kelley95}, \cite{mclach2002} and \cite{geiser2016_1}.
These combinations are efficient and the numerical error can be reduced by the iterative approaches.

The rest of this paper is structured as follows. In section \ref{modell} we present our mathematical model.
In section \ref{simple-model}, we present a simplification of the mathematical model to obtain a 
computable model.
The different solver parts are presented in section \ref{solution}.
The numerical algorithms and examples are given in Section \ref{experiments}.
Finally, in Section \ref{concl}, 
we summarise our results. 

\section{Mathematical Model}
\label{modell}

The starting point for plasma gas mixtures is given in the following reference frame, see also \cite{giovang1999}, \cite{giovang2010} and \cite{orlach2018}.
We concentrate on the heavy particle description, which is discussed in \cite{sene06}.

The distribution function of the heavy particles are given as $f_i({\bf x}, {\bf c}_i, t)$, while ${\bf x}$ is the three-dimensional
spatial coordinate, ${\bf c}_i$ is the velocity of the molecule and $t$ is the time.

The heavy-particle species distribution are given as: 
\begin{eqnarray}
&& {\cal D}_i(f_i) = {\cal S}_i(f) + {\cal C}_i(f) , \; i \in I 
\end{eqnarray}
where ${\cal S}_i(f)$ is the scattering source term and given in \cite{giovang1999}.
${\cal C}_i(f)$ is the reactive source term and given in \cite{giovang1999}.
The differential operator is given as
\begin{eqnarray}
\label{boltz_1}
&& {\cal D}_i(f_i) = \frac{\partial}{\partial t} f_i + \bold c_i \cdot \nabla_{\bold x} f_i + \frac{q_i}{m_i} (\bold E + \bold c_i \times \bold B) \cdot \nabla_{\bold c_i} f_i , \; i \in I .
\end{eqnarray}
Further $q_i$ is the charge of the $i$-th species, $m_i$ the mass of the $i$-th species and $\bold E$, $\bold B$ are the electric and
magnetic fields, we also assume ${\bf b}_i = \frac{q_i}{m_i} (\bold E + \bold c_i \times \bold B)$ is the external force, related to the electro-magnetic field.

In the next step, we apply the Chapman-Enskog expansion, while the zero-th
terms correspond to a Maxwellian distribution and we obtain the
Euler equations. The first-order perturbed distribution function, where a
linearized Boltzmann equation is applied, lead to a Navier-Stokes equation,
see \cite{giovang1999}.

We rewrite the generalized Boltzmann-equation into an Enskog-expansion, see \cite{giovang1999}:
\begin{eqnarray}
\label{boltz_2}
&& {\cal D}_i(f_i) = \frac{1}{\epsilon} {\cal S}_i(f) + \epsilon {\cal C}_i(f) , \; i \in I , 
\end{eqnarray}
while $\epsilon$ is a scaling factor, while $\frac{1}{\epsilon}$ mean, that fast collisions or nonreactive
collisions drive the heavy species to the Maxwell equilibrium. 

The species distribution functions are given as:
\begin{eqnarray}
&& f_i = f_i^0 (1 + \epsilon \phi_i + \OT(\epsilon^2)) , \; i \in I .
\end{eqnarray}

\subsection{Zeroth order approximation}

For the equation with powers $\frac{1}{\epsilon}$ in (\ref{boltz_2}), we have:
\begin{eqnarray}
{\cal S}_i(f^0) = 0 , \; i \in I , 
\end{eqnarray}
with $f^0 = (f_i^0)_{i \in I}$ and it follows the Maxwell distribution function.

For the equations with power $\epsilon^0$, we obtain the zero-th order
macroscopic equations, which are given as the Euler's equations:
\begin{eqnarray}
&& \frac{\partial \rho_i}{\partial t} + \nabla_{\bold x} \cdot ( \rho_i \bold v)= m_i \omega_i^0 , \; i \in I , \\
&& \frac{\partial (\rho {\bf v})}{\partial t} + \nabla_{\bold x} \cdot \left(\rho \bold u  \otimes \bold u + p I \right) = \sum_{i=1}^I \rho_i {\bf b}_i , \\
&& \frac{\partial (\frac{1}{2} \rho {\bf v} \cdot {\bf v} + {\cal E} )}{\partial t} + \nabla_{\bold x} \cdot \left( ( \frac{1}{2} \rho {\bf v} \cdot {\bf v} + {\cal E} + p) {\bf v} \right) =  \sum_{i=1}^I \rho_i {\bf v} \cdot {\bf b}_i ,
\end{eqnarray}
where $\rho = \sum_{i = 1}^I \rho_i$ is the mass density of all species, $p$ is the thermodynamic pressure. $\omega_i^0$ is the zero-th order production rate of species $i$ with:
\begin{eqnarray}
&& \omega_i^0 = \sum_{{\cal I} \in Q_i} \int {\cal C}_i (f^0) \; d {\bf c}_i ,
\end{eqnarray}
where $Q_i$ is the set of the quantum internal energy states of ${\cal I}$ of species $i$.
The internal energy is given as:
\begin{eqnarray}
{\cal E }= \sum_{i=1}^I \sum_{{\cal I} \in Q_i} \int ( \frac{1}{2} m_i ({\bf c}_i - {\bf v}) \cdot ({\bf c}_i - {\bf v}) + E_{i I})  f_i^0 d c_i , 
\end{eqnarray}
see \cite{}.

\subsection{First order approximation}

For the first order approximation, a linearized Boltzmann operator
around the Maxwellian distribution is used, see \cite{giovang1999}.

We have a linearized Boltzmann equation, which is given as:
\begin{eqnarray}
{\cal J}_i^{\cal S}(\phi) = \Phi_i , \; i \in I , \\
 \Phi_i = - {\cal D}_i(\log f_i^0) +  \frac{{\cal C}_i(f^0)}{f_i^0}   , \; i \in I ,
\end{eqnarray}
with $({\cal J}_i^{\cal S})_{i \in I}$ is the linearized Boltzmann operator, see \cite{giovang1999}.

For the equations with power $\epsilon^1$, we obtain the first order
macroscopic equations, which are given as the macroscopic equations in the Navier-Stokes regime:
\begin{eqnarray}
\label{navier_1}
&& \frac{\partial \rho_i}{\partial t} + \nabla_{\bold x} \cdot ( \rho_i \bold v) +  \nabla_{\bold x} (\rho_i {\cal V}_i) = m_i \omega_i^0 , \; i \in I , \\
&& \frac{\partial (\rho {\bf v})}{\partial t} + \nabla_{\bold x} \cdot \left(\rho \bold u  \otimes \bold u + p I \right) + \nabla_{\bold x} {\cal P} = \sum_{i=1}^I \rho_i {\bf b}_i , \\
\label{navier_3}
&& \frac{\partial (\frac{1}{2} \rho {\bf v} \cdot {\bf v} + {\cal E} )}{\partial t} + \nabla_{\bold x} \cdot \left( ( \frac{1}{2} \rho {\bf v} \cdot {\bf v} + {\cal E} + p) {\bf v} \right) + \nabla_{\bold x} ({\cal Q} + {\cal P} \cdot {\bf v}) = \\
&& =  \sum_{i=1}^I \rho_i \cdot {\bf b}_i ( {\bf v} +  {\cal V}_i) ,
\end{eqnarray}
where we have the following operators:

\begin{itemize}
\item The species diffusion velocities ${\cal V}_i$:
\begin{eqnarray}
&& \rho_i {\cal V}_i = m_i \sum_{{\cal I} \in Q_i} \int ({\bf c}_i - {\bf v}) f_i^0 \; \phi_i \; d {\bf c}_i , \; i \in I ,
\end{eqnarray}
\item The viscous tensor ${\cal P}$:
\begin{eqnarray}
{\cal P }= \sum_{i=1}^I \sum_{{\cal I} \in Q_i} \int m_i ({\bf c}_i - {\bf v}) \otimes ({\bf c}_i - {\bf v}) f_i^0 \; \phi_i \; d c_i , 
\end{eqnarray}
\item and the heat flux ${\cal Q}$:
\begin{eqnarray}
{\cal Q }= \sum_{i=1}^I \sum_{{\cal I} \in Q_i} \int ( \frac{1}{2} m_i ({\bf c}_i - {\bf v}) \cdot ({\bf c}_i - {\bf v}) + E_{i I}) ({\bf c}_i - {\bf v})  f_i^0 \; \phi_i \;  d c_i , 
\end{eqnarray}
see \cite{}.

\end{itemize}

\section{Simplified mathematical Model for three species}
\label{simple-model}

This section will present a simplified mathematical model, which concentrates on the
first equation of the Navier-Stokes type equations for the heavy species, see Section \ref{modell}.

We assume that we have ${\bf} v = 0$ in a so called isobaric case, see \cite{bothe2011}.

Then, the Navier-Stokes regime (\ref{navier_1})-(\ref{navier_3}) reduced to a convection-diffusion reaction
equation, which are also developped in the works of \cite{sene06} and \cite{gobb96}.

This model considers the mass-transport of a 
hydrogen plasma. Here, we deal with a hydrogen plasma
that is a mixture of $H, H_2, H_2^+$ particles, means atoms, molecules and ions.

We take into account the dissociation and ionisation reactions, which are
given as:
\begin{eqnarray}
H_2 + e \; \underleftrightarrow{\lambda_1} \; H_2^+ + 2 e , \\
H_2 + e \; \underleftrightarrow{\lambda_2} \; 2 H + e ,
\end{eqnarray}
where the electron temperature is given as $T_e = 17400 \; [K]$ 
and the gas temperature values remain constant $T_h = 600 \; [K]$. \\

Furthermore, we have $\lambda_1 = 1.58 \; 10^{-15} \; T_e^{0.5} \exp(\frac{-15.378}{T_e}) = 2.082 \; 10^{-13}$ \\
 and $\lambda_2 = 1.413 \; 10^{-15} \; T_e^{2} \exp(\frac{-4.48}{T_e}) = 4.276 \; 10^{-7}$.

The diffusion coefficients are given in the following formula:
\begin{eqnarray}
D_{ij} = \frac{3}{16} \frac{f_{ij} k_B^2 T_i T_j}{p\; m_{ij} \Omega_{i j}^{(1,1)}(T_{i j})} ,
\end{eqnarray}
where the parameters are: \\
$f_{ij}$ is a correction factor of order unity, $m_{ij} = \frac{m_i \; m_j}{m_i + m_j}$ is 
the reduced mass, $m_i$ is the mass of species $i$, $m_j$ is the mass of species $j$, $p$ is 
pressure, $t_i, T_j$ is the temperature of the corresponding species, and
$\Omega_{ij}^{(1,1)}$ is a collision integral \cite{hirsch1966}.

We assume the following binary diffusion parameters for our experiments:
\begin{eqnarray}
D_{H_2, H_2^+} = 0.34 \; [cm^2/sec] , \\
D_{H_2, H} = 0.21 \; [cm^2/sec] , \\
D_{H_2^+, H} = 0.21 \; [cm^2/sec] .
\end{eqnarray}

We have used the following Stefan-Maxwell model as a transport model for the
gaseous species.
The modelling equation is given as:
\begin{eqnarray}
\label{ord_0}
&& \partial_t \xi_i + \nabla \cdot N_i = S_i , \; 1 \le i \le 3 , \\
&& \sum_{j=1}^3 N_j = 0 , \\
&& \frac{\xi_2 N_1 - \xi_1 N_2}{D_{12}} + \frac{\xi_3 N_1 - \xi_1 N_3}{D_{13}} =  -  \nabla \xi_1 , \\
 && \frac{\xi_1 N_2 - \xi_2 N_1}{D_{12}} + \frac{\xi_3 N_2 - \xi_2 N_3}{D_{23}} =  -  \nabla \xi_2 ,
\end{eqnarray}
where $\xi_i$ are the mole fractions and $N_i$ is the molar flux of
species $i$, see \cite{bothe2010} and \cite{boudin2012}.
Furthermore, the kinetic term or reaction term $S_i$ is given as:
\begin{eqnarray}
\label{ord_1_1}
&& S_i =  \sum_{j=1}^3 \lambda_{i,j} \xi_j  ,
\end{eqnarray}
where $\lambda_{i,j}$ are the reaction-rates.
The domain is given as $\Omega \in \R^d, d \in \N^+$ with $\xi_i \in C^2$.

We decompose the diffusion and the reaction part, and apply the following 
splitting approach to our problem, we compute $n = 1, \ldots, N$, $t_0, t_1, \ldots, t_n$ time-steps:
The first step is given as (Diffusion step):
\begin{eqnarray}
\label{ord_0}
&& \partial_t \tilde{\xi}_i + \nabla \cdot N_i = 0 , \; 1 \le i \le 3 , \\
&& \sum_{j=1}^3 N_j = 0 , \\
&& \frac{\tilde{\xi}_2 N_1 - \tilde{\xi}_1 N_2}{D_{12}} + \frac{\tilde{\xi}_3 N_1 - \xi_1 N_3}{D_{13}} =  -  \nabla \xi_1 , \\
 && \frac{\tilde{\xi}_1 N_2 - \tilde{\xi}_2 N_1}{D_{12}} + \frac{\tilde{\xi}_3 N_2 - \tilde{\xi}_2 N_3}{D_{23}} =  -  \nabla \tilde{\xi}_2 , \mbox{for} \; t \in [t^n, t^{n+1}], \\
&& \tilde{\xi}_i(t^n) = \xi_i(t^n) , \; i = 1,2,3,
\end{eqnarray}
and the next step is given as (Reaction step): 
\begin{eqnarray}
\label{ord_0}
&& \partial_t \xi_i = S_i , \; 1 \le i \le 3 , \mbox{for} \; t \in [t^n, t^{n+1}], \\
&& \xi_i(t^n) = \tilde{\xi}_i(t^{n+1}) , \; i = 1,2,3 .
\end{eqnarray}

In the following section, we will discuss the different treatments of the subproblems.

\section{Solution of the Transport-Reaction Equation}
\label{solution}

The transport-reaction equation can be solved in the two parts
of the transport part, which is a Stefan-Maxwell equation, and the
reaction part, which is a pure ODE.

These two different approaches are discussed in the following schemes:
\begin{enumerate}
\item Stefan-Maxwell Problem (Diffusion-part):

We concentrate on the three component system and
solve this system as a linear optimal problem (General Linear Optimal Problem).
We deal with:
\begin{eqnarray}
\label{ord_0}
&& \partial_t \xi_i + \nabla \cdot N_i = 0 , \; 1 \le i \le 3 , \\
&& \sum_{j=1}^3 N_j = 0 , \\
&& \frac{\xi_2 N_1 - \xi_1 N_2}{D_{12}} + \frac{\xi_3 N_1 - \xi_1 N_3}{D_{13}} =  -  \nabla \xi_1 , \\
 && \frac{\xi_1 N_2 - \xi_2 N_1}{D_{12}} + \frac{\xi_3 N_2 - \xi_2 N_3}{D_{23}} =  -  \nabla \xi_2 ,
\end{eqnarray}
where the domain is given as $\Omega \in \R^d, d \in \N^+$ with $\xi_i \in C^2$.

We could reduce this to a simpler model problem as:
\begin{eqnarray}
\label{ord_0}
&& \partial_t \xi_i + \nabla \cdot N_i = 0 , \; 1 \le i \le 2 , \\
&& \frac{1}{D_{13}} N_1 + \alpha N_1 \xi_2 - \alpha N_2 \xi_1 =  -  \nabla \xi_1 , \\
 && \frac{1}{D_{23}} N_2 - \beta N_1 \xi_2 + \beta N_2 \xi_1 =  -  \nabla \xi_2 ,
\end{eqnarray}
where we have $\alpha = \left(\frac{1}{D_{12}} - \frac{1}{D_{13}}\right)$, 
$\beta = \left(\frac{1}{D_{12}} - \frac{1}{D_{23}}\right)$. 

\vspace{1cm}

The optimal problem is derived in the following manner.

Second, we rewrite the MOR model equation (\ref{ord_0}) to a
set of $s$ linearised states $U_0, U_1, \ldots, U_s$ by using the linear system:
\begin{eqnarray}
\label{ord_0}
&&  U'_{i+1}  = J_i(t) U_{i+1} + \hat{B}(t) v ,  
\end{eqnarray}
where $J_i$ is the Jacobian of $B(U,t)$ and is given in (\ref{ord_0}),
the control operator is $\hat{B}(t) = \tilde{B}(t) - J_i$, and the system
input is $v = U_i$.

Third, we can now apply the GLCS, using the following
notations: $u = U_{i+1}, v = U_i$, $A_1(t) = J_i(t)$,
$A_2(t) = \tilde{B}(t)$.

The GLCS is then
\begin{eqnarray}
\label{control_1}
 && \frac{d u}{d t} = A_1(t) u + A_2(t) v , \\
\label{control_2}
 &&  \tilde{u} = C(t) u + D(t) v , \\
 && u(0) = u_0 , 
\end{eqnarray}
where the time-dependent operators are $A(t) \in {\bf X}^n \times {\bf X}^n$,
$B(t) \in {\bf X}^n \times {\bf X}^m$, $C(t) \in {\bf X}^p \times {\bf X}^n$,
$D(t) \in {\bf X}^p \times {\bf X}^m$, 
$v: {\bf X} \rightarrow {\bf X}^m$ denotes the system input,
$\tilde{u}:  {\bf X} \rightarrow {\bf X}^p$ is the system output and
$u:  {\bf X} \rightarrow {\bf X}^n$ denotes the state vector.
Furthermore, ${\bf X}$ is an appropriate Banach space; for example, $U$, a space of 
continuous or piece-wise continuous functions.

The analytical solution of (\ref{control_1}) and (\ref{control_2}) is
\begin{eqnarray}
\label{control_ana_1}
  u(t) &&  =  \exp(\int_{0}^{t} A_1(s) ds u_0 + \int_{0}^{t} \exp(\int_{s}^{t}  A_1(\tilde{s}) d\tilde{s}  A_2(s) v(s) ds, \\
\label{control_ana_2}
  \tilde{u}(t) && =  C(t) \; \exp(\int_{0}^{t} A_1(s) ds u_0 \\
 && + C(t) \; \int_{0}^{t} \exp(\int_{s}^{t}  A_1(\tilde{s}) d\tilde{s}  A_2(s) v(s) ds + D(t) v(t) , \nonumber
\end{eqnarray}
where we apply the fast computation of the exponential integral matrices
via the Magnus expansion, see \cite{blan08}, \cite{blan05} and \cite{casas2006}, and which is discussed in the following.

\item Kinetic Problem (Reaction-part):

We concentrate on the three component system and
we deal with:
\begin{eqnarray}
\label{ord_0}
&& \partial_t \xi_i = S_i , \; 1 \le i \le 3 , 
\end{eqnarray}
where the domain is given as $\Omega \in \R^d, d \in \N^+$ with $\xi_i \in C^2$.

We apply the reaction-rates and have the following linear ODE system:
\begin{eqnarray}
\label{ord_0}
&& \partial_t {\bf \xi} = S {\bf \xi} , 
\end{eqnarray}
where ${\bf \xi} = (\xi_1, \xi_2, \xi_3)^t$ and
$S = \left( \begin{array}{c c c} \lambda_{1,1} & \lambda_{1,2} & \lambda_{1,3} \\ \lambda_{2,1} & \lambda_{2,2} & \lambda_{2,3} \\ \lambda_{3,1} & \lambda_{3,2} & \lambda_{3,3}  \end{array} \right)$.

We can apply the analytical solution, which is given as:
\begin{eqnarray}
\label{ord_0}
&& {\bf \xi}(t^{n+1}) = \exp( S \Delta t) {\bf \xi}(t^n) , 
\end{eqnarray}
and $\Delta t = t^{n+1} - t^n$.

\end{enumerate}

\section{Numerical Algorithms and Numerical Experiments}
\label{experiments}

In this section, we discuss the different numerical algorithms that are
based on splitting approaches and which are to solve the multicomponent transport-reaction
equations.

We deal with the following two experiments:
\begin{itemize}
\item Pure diffusion problem, here we only apply the Stefan-Maxwell equation.
\item Hydrogen Plasma, here we apply the Stefan-Maxwell equation with the reaction equation.
\end{itemize}

\subsection{Pure Diffusion Problem}

We concentrate on the three component system:
\begin{eqnarray}
\label{ord_0}
&& \partial_t \xi_i + \partial_x N_i = 0 , \; 1 \le i \le 3 , \\
&& \sum_{j=1}^3 N_j = 0 , \\
&& \frac{\xi_2 N_1 - \xi_1 N_2}{D_{12}} + \frac{\xi_3 N_1 - \xi_1 N_3}{D_{13}} =  -  \partial_x \xi_1 , \\
 && \frac{\xi_1 N_2 - \xi_2 N_1}{D_{12}} + \frac{\xi_3 N_2 - \xi_2 N_3}{D_{23}} =  -  \partial_x \xi_2 ,
\end{eqnarray}
where the domain is given as $\Omega \in \R^d, d \in \N^+$ with $\xi_i \in C^2$.

The parameters and the initial and boundary conditions are given as:
\begin{itemize}
\item $D_{12} = D_{13} = 0.833$ (means $\alpha = 0$) and $D_{23} = 0.168$ (uphill diffusion, semi-degenerated Duncan and Toor experiment)
\item $D_{12} = 0.0833, D_{13} = 0.680$ and $D_{23} = 0.168$ (asymptotic behavior, Duncan and Toor experiment)
\item $J = 140$ (spatial grid points)
\item  The time-step-restriction for the explicit method is given as: \\
 $\Delta t \le (\Delta x)^2 \max \{ \frac{1}{2 \{D_{12}, D_{13}, D_{23}\}} \}$
\item The spatial domain is $\Omega = [0, 1]$, the time-domain $[0, T] = [0, 1]$
\item The initial conditions are:
\begin{enumerate}
\item Uphill example
\begin{eqnarray}
\label{init}
&& \xi_1^{in}(x) = \left\{ \begin{array}{l l}
0.8 & \mbox{if} \; 0 \le x < 0.25 , \\
1.6 (0.75 - x) & \mbox{if} \; 0.25 \le x < 0.75 , \\
0.0 & \mbox{if} \; 0.75 \le x \le 1.0 , 
\end{array} \right. , \\
&& \xi_2^{in}(x) = 0.2 , \; \mbox{for all} \; x \in \Omega = [0,1] ,
\end{eqnarray}
\item Diffusion example (asymptotic behavior)
\begin{eqnarray}
\label{init}
&& \xi_1^{in}(x) = \left\{ \begin{array}{l l}
0.8 & \mbox{if} \; 0 \le x  \in 0.5 , \\
0.0 & \mbox{else} , 
\end{array} \right. , \\
&&  \xi_2^{in}(x) = 0.2 ,  \; \mbox{for all} \; x \in \Omega = [0,1] ,
\end{eqnarray}
\end{enumerate}

\item The boundary conditions are of no-flux type:
\begin{eqnarray}
\label{init}
&& N_1 = N_2 = N_3 = 0 , \mbox{on} \; \partial \Omega \times [0,1] ,
\end{eqnarray}
\end{itemize}

We could reduce this to a simpler model problem, as follows:
\begin{eqnarray}
\label{ord_0}
&& \partial_t \xi_i + \partial_x \cdot N_i = 0 , \; 1 \le i \le 2 , \\
&& \frac{1}{D_{13}} N_1 + \alpha N_1 \xi_2 - \alpha N_2 \xi_1 =  -  \partial_x \xi_1 , \\
 && \frac{1}{D_{23}} N_2 - \beta N_1 \xi_2 + \beta N_2 \xi_1 =  -  \partial_x \xi_2 ,
\end{eqnarray}
where $\alpha = \left(\frac{1}{D_{12}} - \frac{1}{D_{13}}\right)$, 
$\beta = \left(\frac{1}{D_{12}} - \frac{1}{D_{23}}\right)$.

We then rewrite into:
\begin{eqnarray}
\label{ord_0}
&& \partial_t \xi_1 + \partial_x \cdot N_1 = 0 , \\
&& \partial_t \xi_2 + \partial_x \cdot N_2 = 0 , \\
&& \left( \begin{array}{c c}
\frac{1}{D_{13}} + \alpha \xi_2  & - \alpha \xi_1  \\
 - \beta \xi_2 &  \frac{1}{D_{23}}  + \beta \xi_1 
\end{array} \right)
 \left( \begin{array}{l}
N_1  \\
N_2 
\end{array} \right) = 
 \left( \begin{array}{l}
 -  \partial_x \xi_1 \\
 -  \partial_x \xi_2 
\end{array} \right)
\end{eqnarray}
and we have
\begin{eqnarray}
\label{part_0}
&& \partial_t \xi_1 + \partial_x \cdot N_1 = 0 , \\
\label{part_1}
&& \partial_t \xi_2 + \partial_x \cdot N_2 = 0 , \\
\label{part_2}
&& \left( \begin{array}{l}
N_1  \\
N_2 
\end{array} \right) = \frac{D_{13} D_{23}}{1 + \alpha D_{13} \xi_2 + \beta D_{23} \xi_1}
 \left( \begin{array}{c c}
\frac{1}{D_{23}} + \beta \xi_1  & \alpha \xi_1  \\
 \beta \xi_2 &  \frac{1}{D_{13}}  + \alpha \xi_2 
\end{array} \right)
 \left( \begin{array}{l}
 -  \partial_x \xi_1 \\
 -  \partial_x \xi_2 
\end{array} \right)
\end{eqnarray}

The next step is to apply the semi-discretisation of the 
partial differential operator $\frac{\partial}{\partial x}$.

We apply the first differential operator in equation (\ref{part_0}) and (\ref{part_1})
as a forward upwind scheme, which is given as
\begin{eqnarray}
\frac{\partial}{\partial x} & = & D_+ =  \frac{1}{\Delta x}\cdot \left(\begin{array}{rrrrr}
 -1 & 0 & \ldots & ~ & 0 \\
  1 & -1 & 0 & \ldots & 0 \\
 \vdots & \ddots & \ddots & \ddots & \vdots \\
 0 & ~ & 1 & -1 & 0 \\
 0 & \ldots & 0 & 1 & -1
\end{array}\right)~\in~\R^{(J+1) \times (J+1)}
\end{eqnarray}
and we apply the second differential operator in equation (\ref{part_2})
as a backward upwind scheme, which is given as
\begin{eqnarray}
\frac{\partial}{\partial x} & = & D_- =  \frac{1}{\Delta x}\cdot \left(\begin{array}{rrrrr}
 -1 & 1 & 0 & \ldots & 0 \\
 0 & -1 & 1 & 0 & \ldots \\
 \vdots & \ddots & \ddots & \ddots & \ddots \\
 0 & \ldots & 0  & -1 & 1 \\
 0 & ~ & \ldots & 0 & -1
\end{array}\right)~\in~\R^{(J+1) \times (J+1)}
\end{eqnarray}
In the next part, we apply the iterative schemes to solve the
pure diffusion problem.
\subsubsection{Iterative Scheme in Time for the Pure Diffusion Problem}

In this section, we apply a global linearisation of the
Stefan-Maxwell equation.
Then, we consider the underlying semi-discretised equation 
with an iterative approach.

We solve the iterative scheme:
\begin{eqnarray}
\label{ord_0}
&& \xi_{1}^{n+1} = \xi_1^n - \Delta t \; D_+ N_{1}^n , \\
&& \xi_{2}^{n+1} = \xi_2^n - \Delta t \; D_+ N_{2}^n , \\
&& \left( \begin{array}{c c}
A & B \\
C & D
\end{array} \right)
 \left( \begin{array}{l}
N_1^{n+1}  \\
N_2^{n+1} 
\end{array} \right) =
 \left( \begin{array}{l}
 - D_- \xi_1^{n+1} \\
 - D_- \xi_2^{n+1} 
\end{array} \right) 
\end{eqnarray}
for $j = 0, \ldots, J$ , where $\xi_1^n = (\xi_{1,0}^n, \ldots, \xi_{1, J}^n)^T$,
$\xi_2^n = (\xi_{2,0}^n, \ldots, \xi_{2, J}^n)^T$ and $I_J \in \R^{J+1} \times \R^{J+1}$,
 $N_1^n = (N_{1,0}^n, \ldots, N_{1, J}^n)^T$,
$N_2^n = (N_{2,0}^n, \ldots, N_{2, J}^n)^T$ and $I_J \in \R^{J+1} \times \R^{J+1}$,
where $n=0,1,2, \ldots, N_{end}$ and $N_{end}$ are the number of time-steps, i.d. $N_{end} = T / \Delta t$.

The matrices are given as:
\begin{eqnarray}
&& A, B, C, D \in \R^{J+1} \times \R^{J+1}, \\
&& A_{j,j} = \frac{1}{D_{13}} + \alpha \xi_{2,j} , \; j = 0 \ldots, J ,\\
&& B_{j,j} = - \alpha \xi_{1,j} , \; j = 0 \ldots, J , \\
&& C_{j,j} = - \beta \xi_{2, j} , \; j = 0 \ldots, J , \\
&& D_{j,j} =  \frac{1}{D_{23}}  + \beta \xi_{1,j}  , \; j = 0 \ldots, J ,\\
&& A_{i,j} = B_{i,j} =  C_{i,j} = D_{i,j} = 0 , \; i,j = 0 \ldots, J, \; i \neq J , 
\end{eqnarray}
which means that the diagonal entries given as for the scale case in
equation (\ref{part_2}) and the outer-diagonal entries are zero. \\
The explicit form with time-discretisation is given as:

\begin{algorithm}

1.) Initialisation $n=0$:

\begin{eqnarray}
\label{ord_0}
&& \left( \begin{array}{l}
N_1^{0}  \\
N_2^{0} 
\end{array} \right) =
 \left( \begin{array}{c c}
\tilde{A} & \tilde{B} \\
\tilde{C} & \tilde{D}
\end{array} \right)
 \left( \begin{array}{l}
 - D_- \xi_1^{0} \\
 - D_- \xi_2^{0} 
\end{array} \right) 
\end{eqnarray}
where $\xi_1^{0} = (\xi_{1,0}^{0}, \ldots, \xi_{1, J}^{0})^T$, $\xi_2^0 = (\xi_{2,0}^0, \ldots, \xi_{2, J}^0)^T$ and $\xi_{1,j}^{0} = \xi_1^{in}(j \Delta x), \; \xi_{2,j}^{0} = \xi_2^{in}(j \Delta x)$, $j = 0, \ldots, J$ and, given as for the different initialisation, we have:
\begin{enumerate}
\item Uphill example
\begin{eqnarray}
\label{init}
&& \xi_1^{in}(x) = \left\{ \begin{array}{l l}
0.8 & \mbox{if} \; 0 \le x < 0.25 , \\
1.6 (0.75 - x) & \mbox{if} \; 0.25 \le x < 0.75 , \\
0.0 & \mbox{if} \; 0.75 \le x \le 1.0 , 
\end{array} \right. , \\
&& \xi_2^{in}(x) = 0.2 , \; \mbox{for all} \; x \in \Omega = [0,1] ,
\end{eqnarray}
\item Diffusion example (Asymptotic behavior)
\begin{eqnarray}
\label{init}
&& \xi_1^{in}(x) = \left\{ \begin{array}{l l}
0.8 & \mbox{if} \; 0 \le x  \in 0.5 , \\
0.0 & \mbox{else} , 
\end{array} \right. , \\
&&  \xi_2^{in}(x) = 0.2 ,  \; \mbox{for all} \; x \in \Omega = [0,1] ,
\end{eqnarray}
\end{enumerate}

The inverse matrices are given as:
\begin{eqnarray}
&& \tilde{A}, \tilde{B}, \tilde{C}, \tilde{D} \in \R^{J+1} \times \R^{J+1}, \\
&& \tilde{A}_{j,j} = \gamma_j (\frac{1}{D_{23}}  + \beta \xi_{1,j}^{0}) , \; j = 0 \ldots, J ,\\
&& B_{j,j} = \gamma_j \;  \alpha \xi_{1,j}^{0} , \; j = 0 \ldots, J , \\
&& C_{j,j} = \gamma_j \;  \beta \xi_{2, j}^{0} , \; j = 0 \ldots, J , \\
&& D_{j,j} = \gamma_j  (\frac{1}{D_{13}} + \alpha \xi_{2,j}^{0})  , \; j = 0 \ldots, J ,\\
&& \gamma_j = \frac{D_{13} D_{23}}{1 + \alpha D_{13} \xi_{2,j}^{0} + \beta D_{23} \xi_{1,j}^{0}} ,  \; j = 0 \ldots, J , \\
&& \tilde{A}_{i,j} = \tilde{B}_{i,j} =  \tilde{C}_{i,j} = \tilde{D}_{i,j} = 0 , \; i,j = 0 \ldots, J, \; i \neq J , 
\end{eqnarray}

The values of the first and the last grid points of $N$ are zero,
which means that $N_{1,0}^{0} = N_{1,J}^{0} = N_{2,0}^{0} = N_{2,J}^{0} = 0$ (boundary condition).

2.) Next time-steps (till $n = N_{end}$ ): \\

2.1) Computation of $\xi_1^{n+1}$ and $\xi_2^{n+1}$ 
\begin{eqnarray}
\label{ord_0}
&& \xi_{1}^{n+1} = \xi_1^n - \Delta t \; D_+ N_{1}^n , \\
&& \xi_{2}^{n+1} = \xi_2^n - \Delta t \; D_+ N_{2}^n ,
\end{eqnarray}

2.2) Computation of $N_1^{n+1}$ and $N_2^{n+1}$ 

\begin{eqnarray}
&& \left( \begin{array}{l}
N_1^{n+1}  \\
N_2^{n+1} 
\end{array} \right) =
 \left( \begin{array}{c c}
\tilde{A} & \tilde{B} \\
\tilde{C} & \tilde{D}
\end{array} \right)
 \left( \begin{array}{l}
 - D_- \xi_1^{n+1} \\
 - D_- \xi_2^{n+1} 
\end{array} \right) 
\end{eqnarray}
where $\xi_1^{n} = (\xi_{1,0}^{n}, \ldots, \xi_{1, J}^{n})^T$, $\xi_2^n = (\xi_{2,0}^n, \ldots, \xi_{2, J}^n)^T$ 
and the inverse matrices are given as:
\begin{eqnarray}
&& \tilde{A}, \tilde{B}, \tilde{C}, \tilde{D} \in \R^{J+1} \times \R^{J+1}, \\
&& \tilde{A}_{j,j} = \gamma_j (\frac{1}{D_{23}}  + \beta \xi_{1,j}^{n+1}) , \; j = 0 \ldots, J ,\\
&& B_{j,j} = \gamma_j \;  \alpha \xi_{1,j}^{n+1} , \; j = 0 \ldots, J , \\
&& C_{j,j} = \gamma_j \;  \beta \xi_{2, j}^{n+1} , \; j = 0 \ldots, J , \\
&& D_{j,j} = \gamma_j  (\frac{1}{D_{13}} + \alpha \xi_{2,j}^{n+1})  , \; j = 0 \ldots, J ,\\
&& \gamma_j = \frac{D_{13} D_{23}}{1 + \alpha D_{13} \xi_{2,j}^{n+1} + \beta D_{23} \xi_{1,j}^{n+1}} ,  \; j = 0 \ldots, J , \\
&& \tilde{A}_{i,j} = \tilde{B}_{i,j} =  \tilde{C}_{i,j} = \tilde{D}_{i,j} = 0 , \; i,j = 0 \ldots, J, \; i \neq J .
\end{eqnarray}

Furthermore, the values of the first and the last grid points of $N$ are zero,
which means that $N_{1,0}^{n} = N_{1,J}^{n} = N_{2,0}^{n} = N_{2,J}^{n} = 0$ (boundary condition).

3.) Do  $n = n+1$ and then goto 2.) 

\end{algorithm}

We have used the following examples: \\

\begin{figure}[ht]
\begin{center}  
\includegraphics[width=8.0cm,angle=-0]{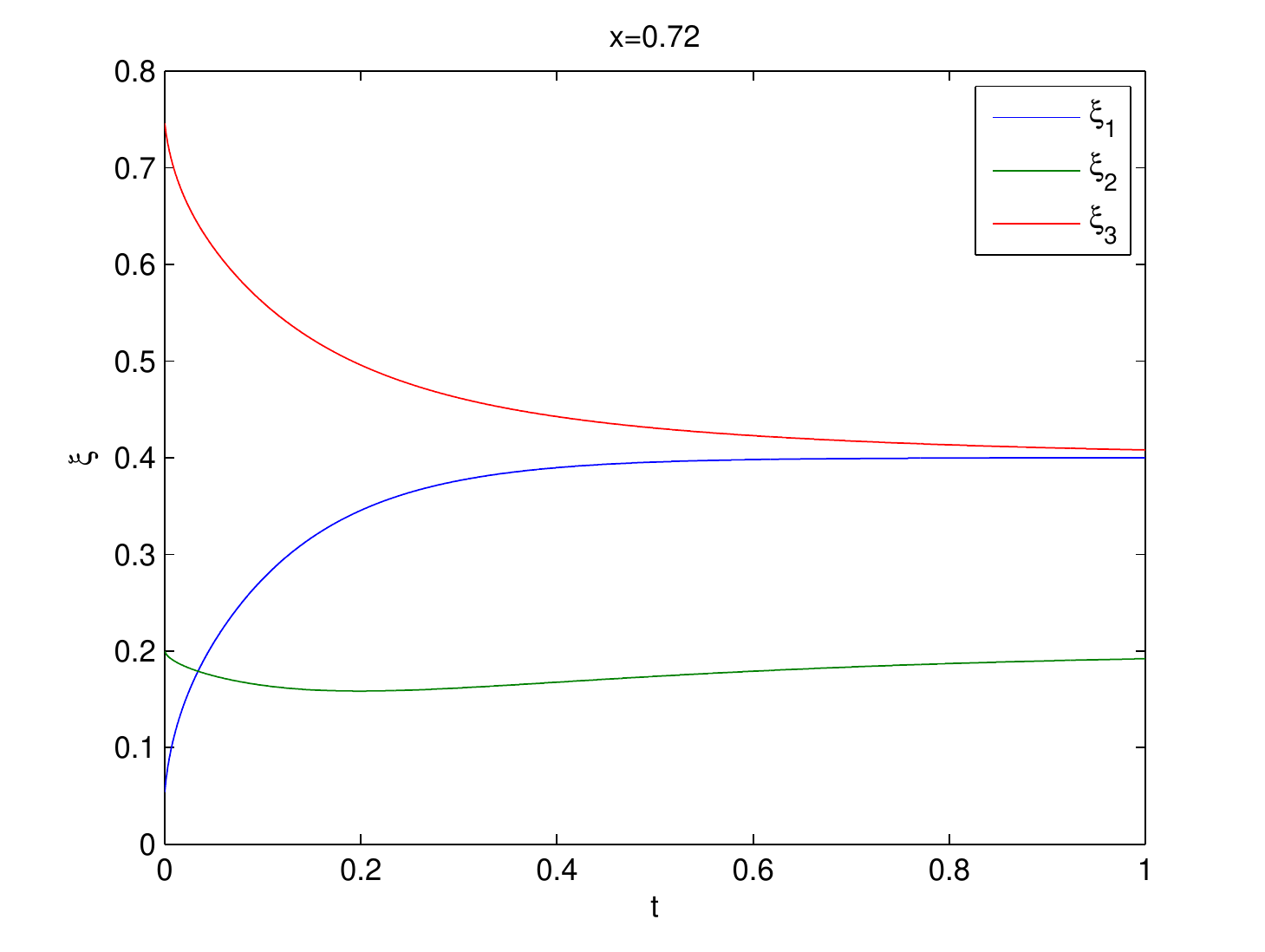}
\end{center}
\caption{\label{multi_1} The figures present the results of the
concentration $c_1$, $c_2$ and $c_3$.}
\end{figure}

We test the different schemes and obtain the results shown in 
Figure \ref{multi_1}.

\begin{figure}[ht]
\begin{center}  
\includegraphics[width=5.0cm,angle=-0]{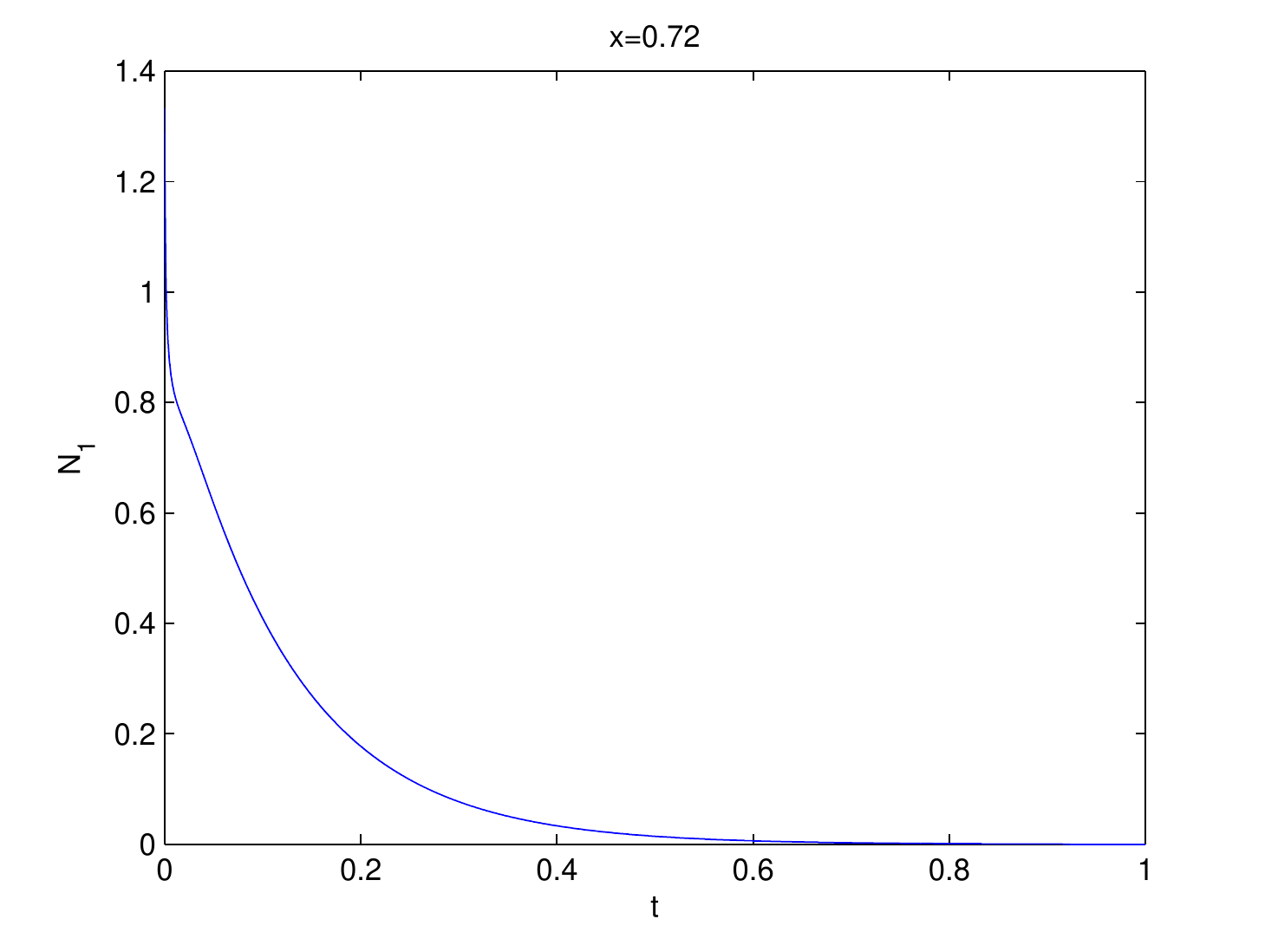}
\includegraphics[width=5.0cm,angle=-0]{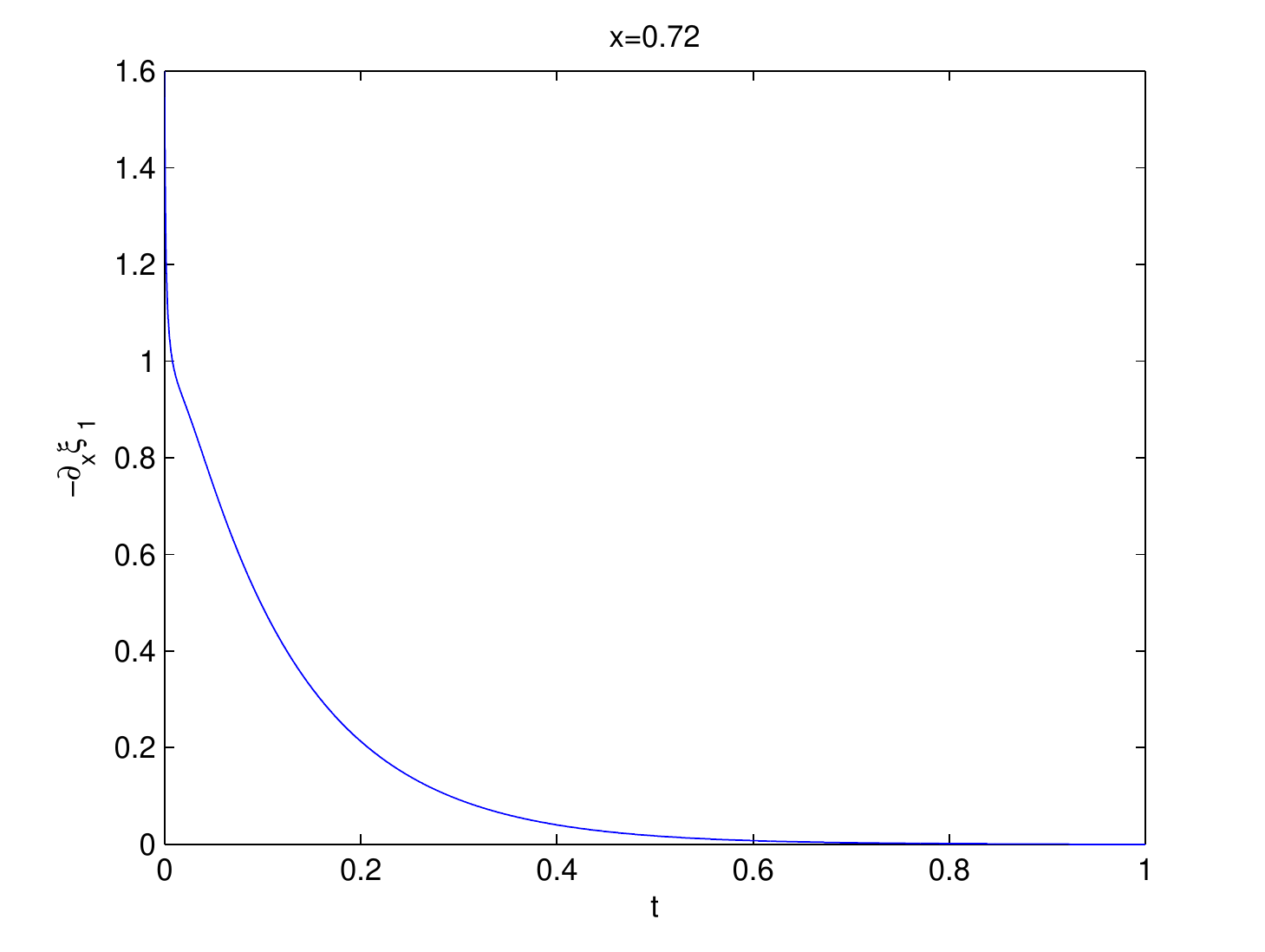} \\
\includegraphics[width=5.0cm,angle=-0]{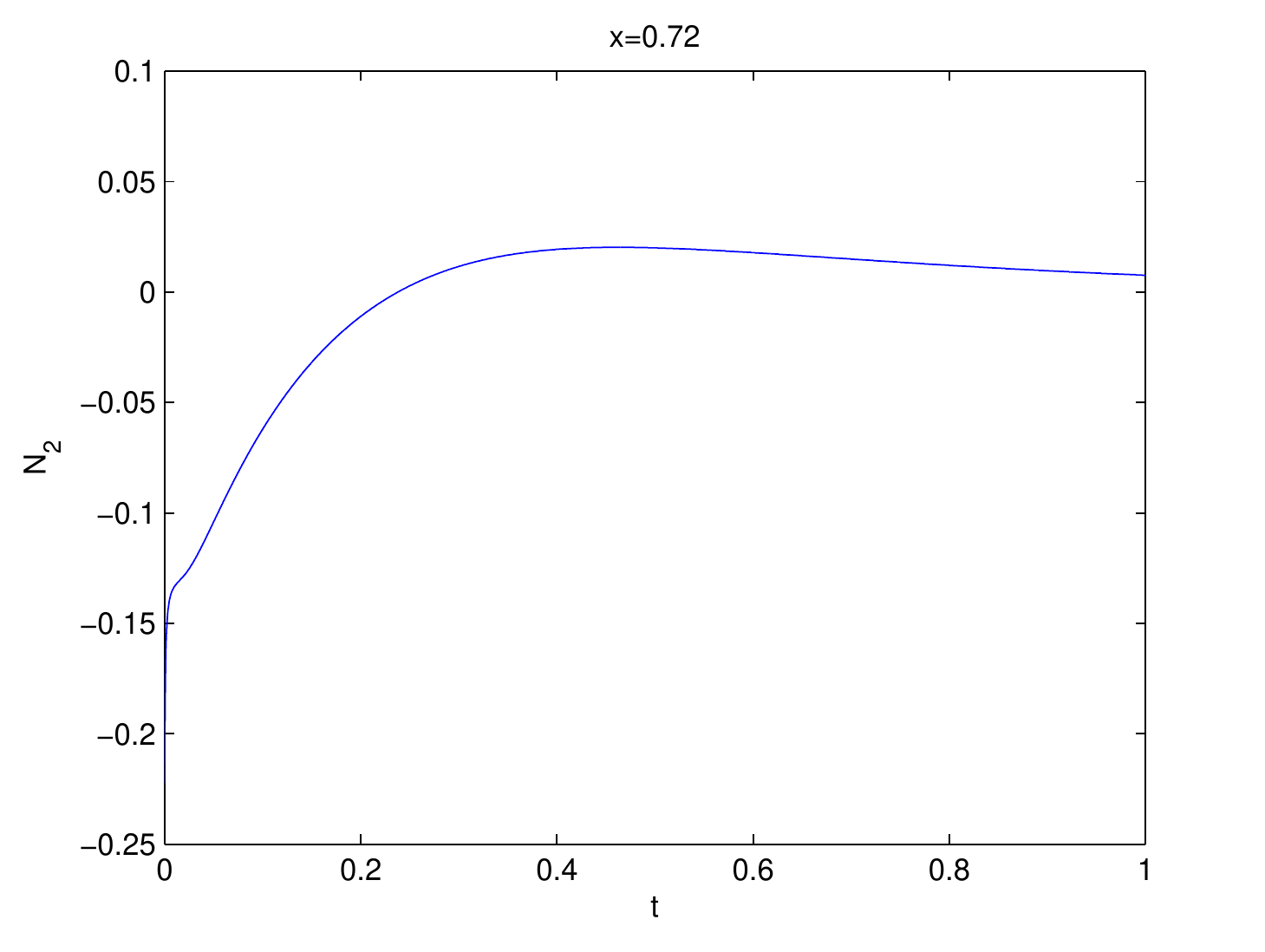}
\includegraphics[width=5.0cm,angle=-0]{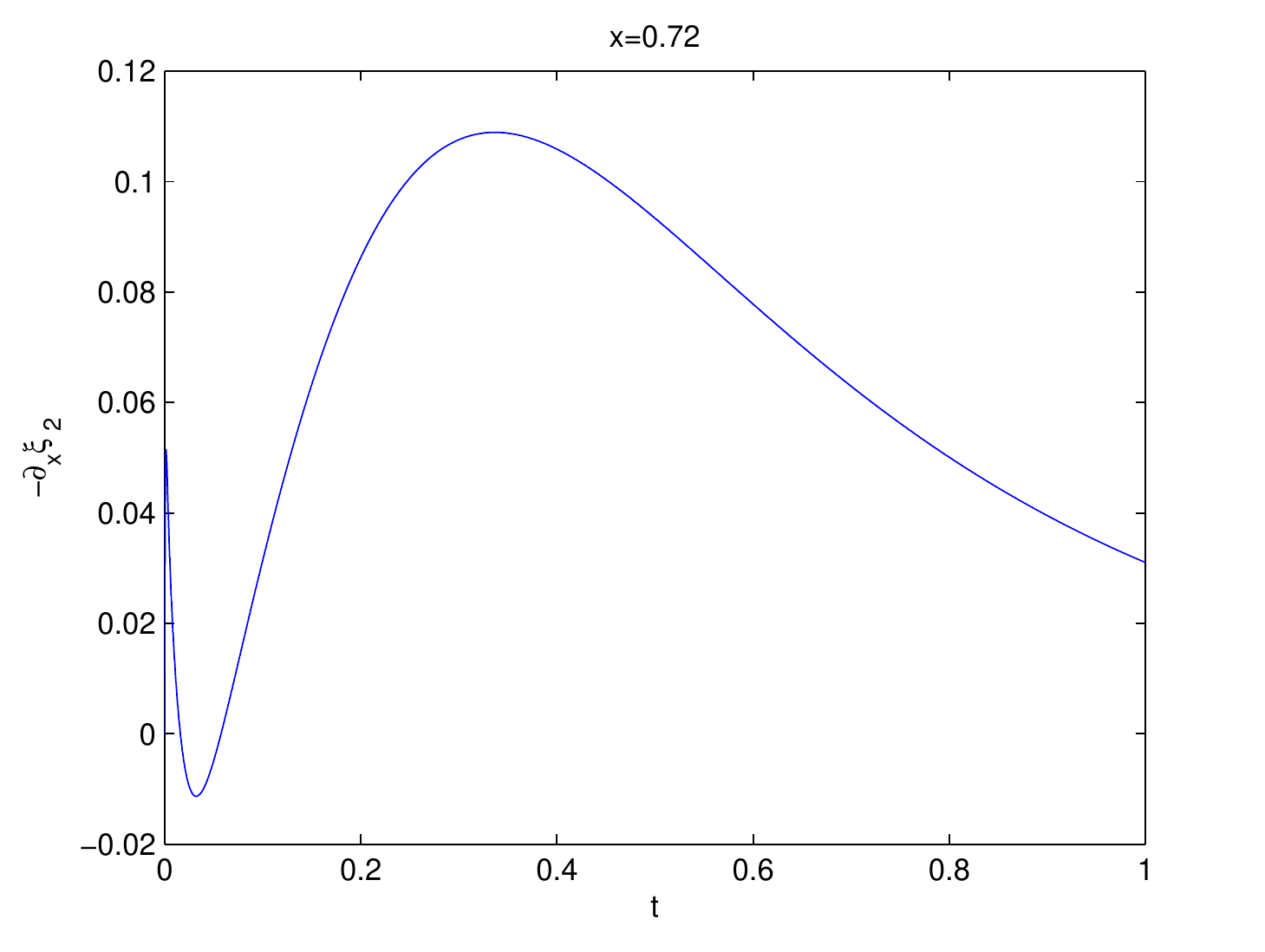}
\end{center}
\caption{\label{multi_2} The upper figures present the results of the
concentration $c_1$ and $- \partial_x \xi_1$. The lower figures present
the results of $c_2$ and $- \partial_x \xi_2$.}
\end{figure}

The concentration and their fluxes are given in Figure \ref{multi_2}.

The full plots in time and space of the concentrations and their fluxes are given in Figure \ref{multi_3}.
\begin{figure}[ht]
\begin{center}  
\includegraphics[width=5.0cm,angle=-0]{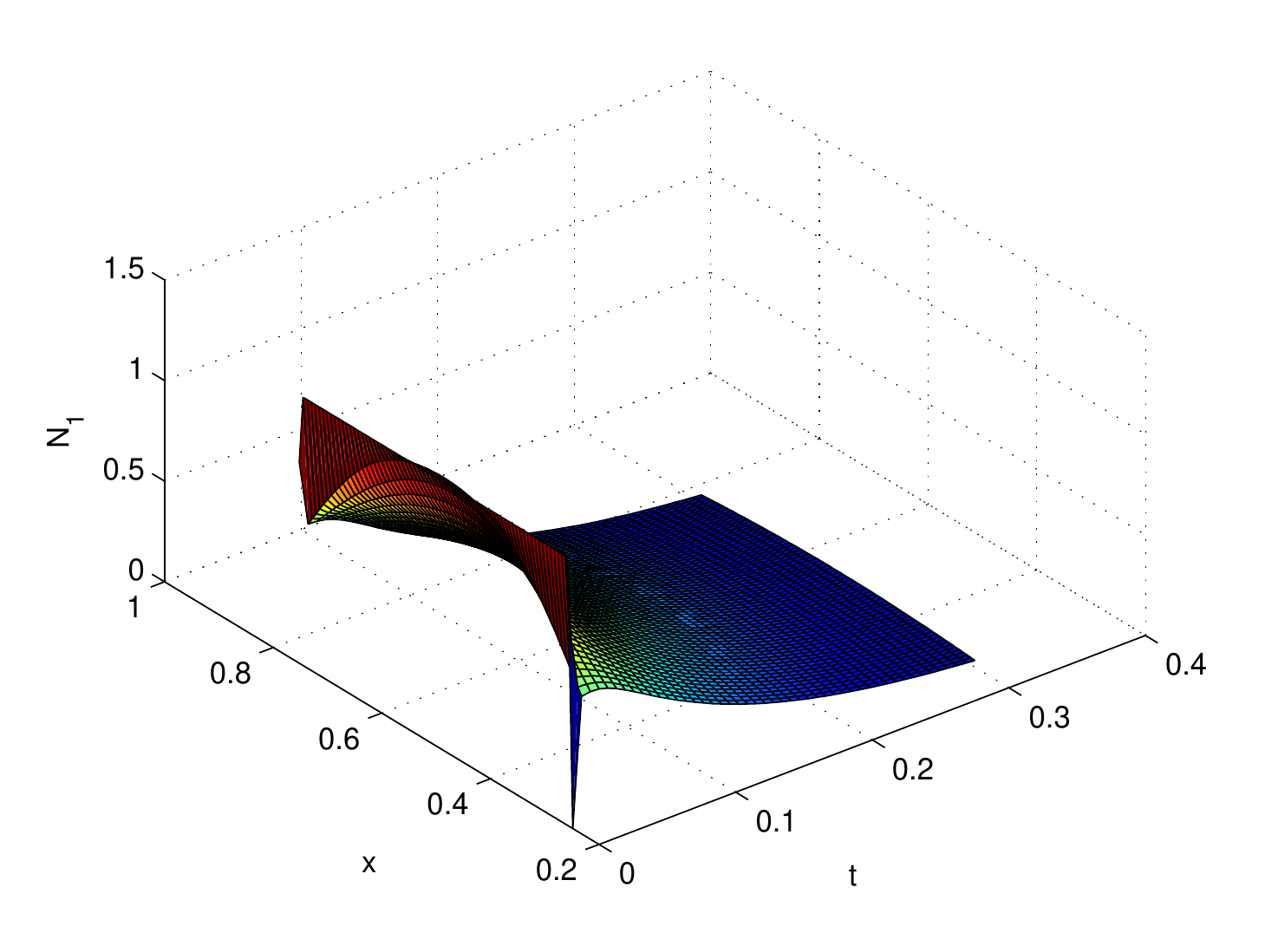}
\includegraphics[width=5.0cm,angle=-0]{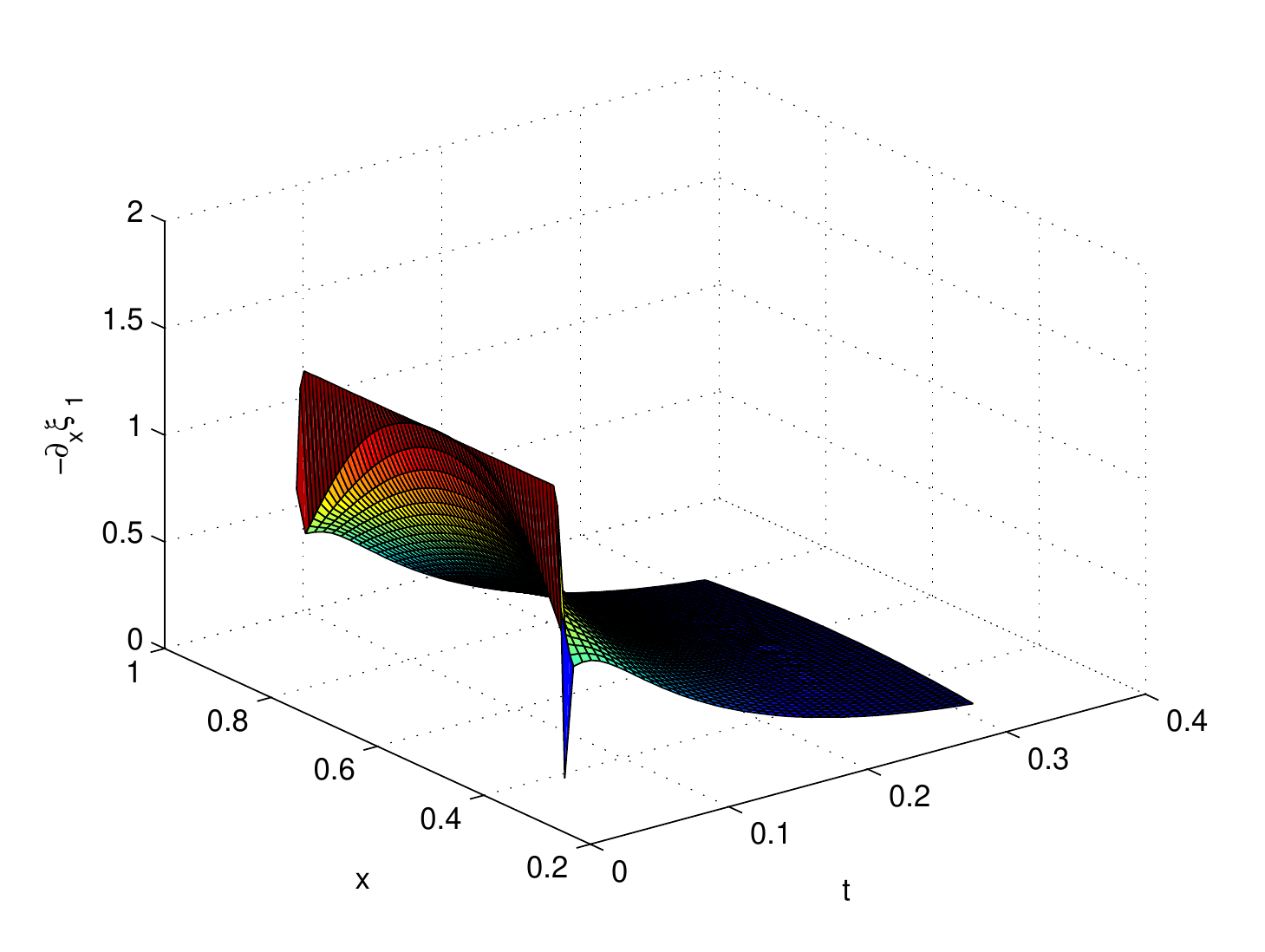} \\
\includegraphics[width=5.0cm,angle=-0]{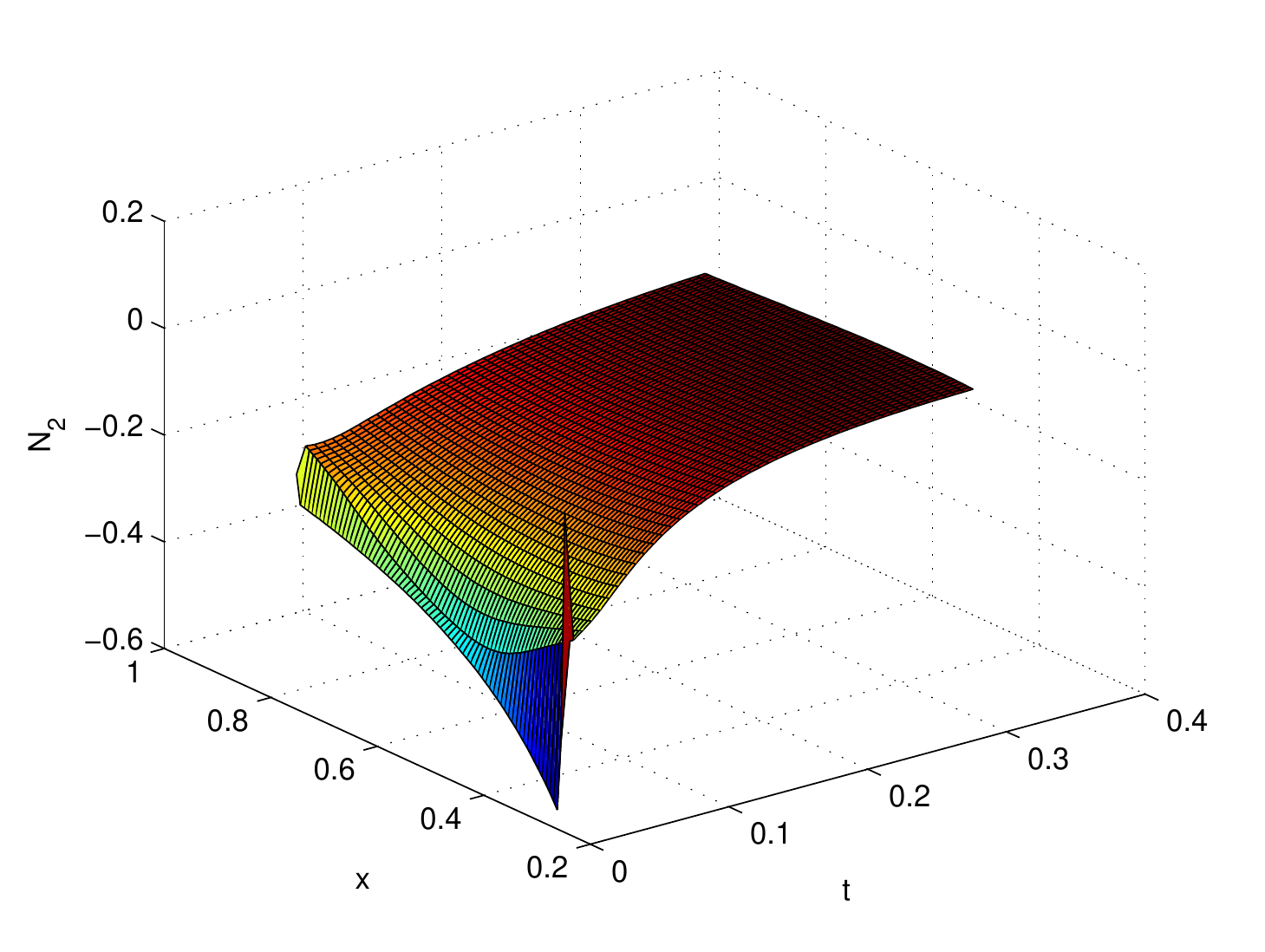}
\includegraphics[width=5.0cm,angle=-0]{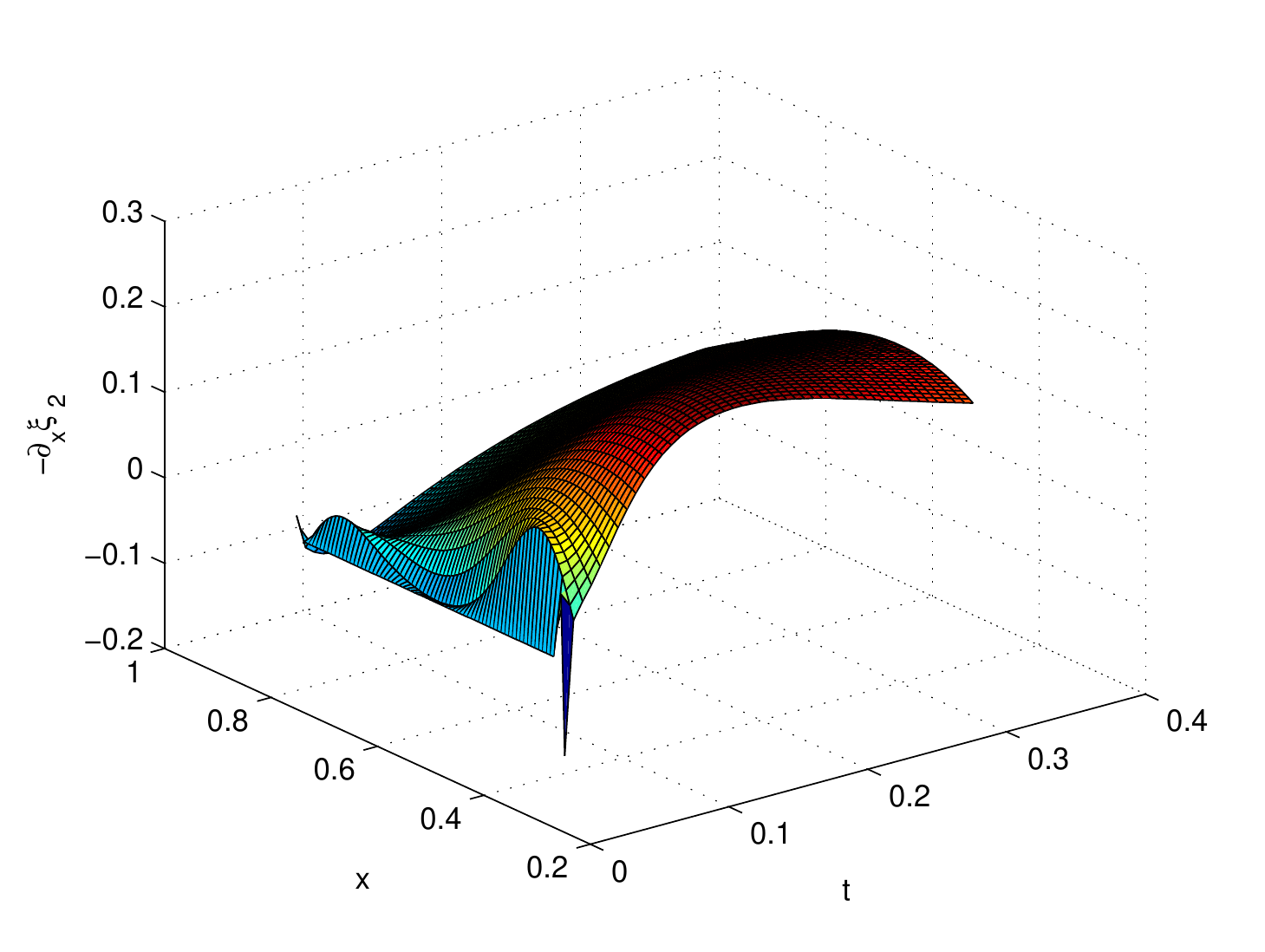}
\end{center}
\caption{\label{multi_3} The figures present the results of the 3d plots in time and space.  The upper figures present the results of the
concentration $c_1$ and $- \partial_x \xi_1$. The lower figures present
the results of $c_2$ and $- \partial_x \xi_2$.}
\end{figure}

The full plots in time and space of the concentrations and their fluxes are given in Figure \ref{multi_3}.
\begin{figure}[ht]
\begin{center}  
\includegraphics[width=5.0cm,angle=-0]{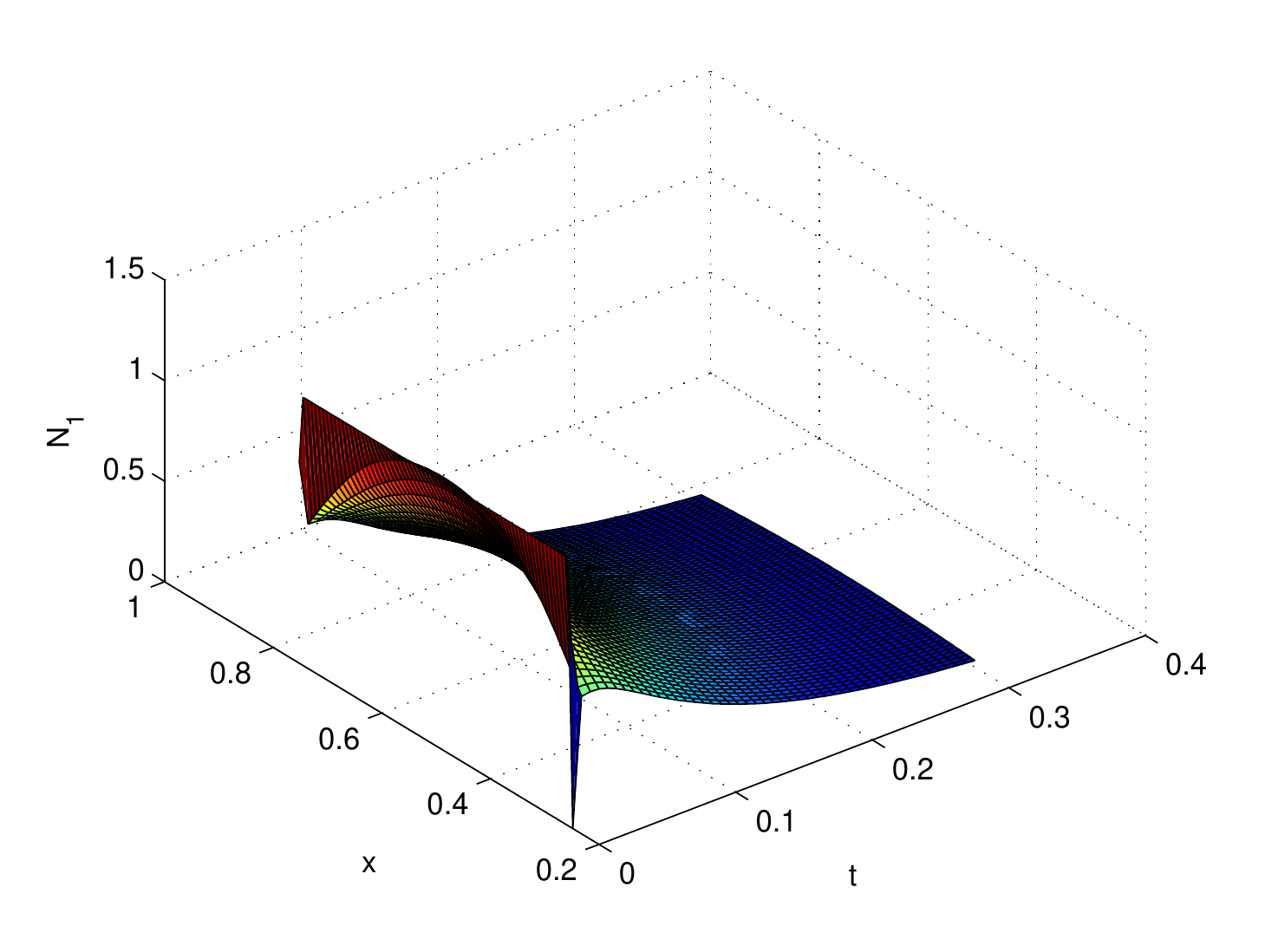}
\includegraphics[width=5.0cm,angle=-0]{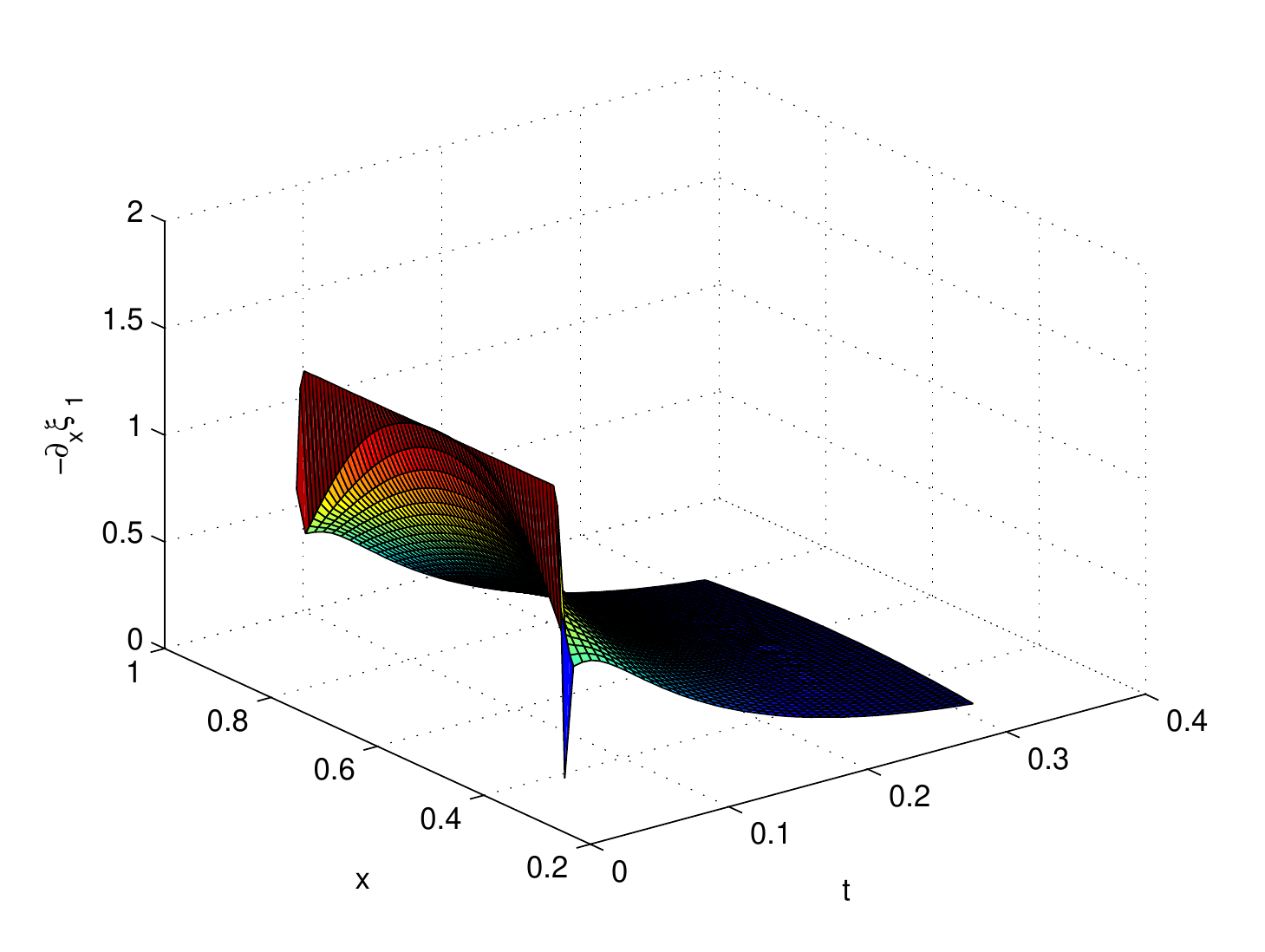} \\
\includegraphics[width=5.0cm,angle=-0]{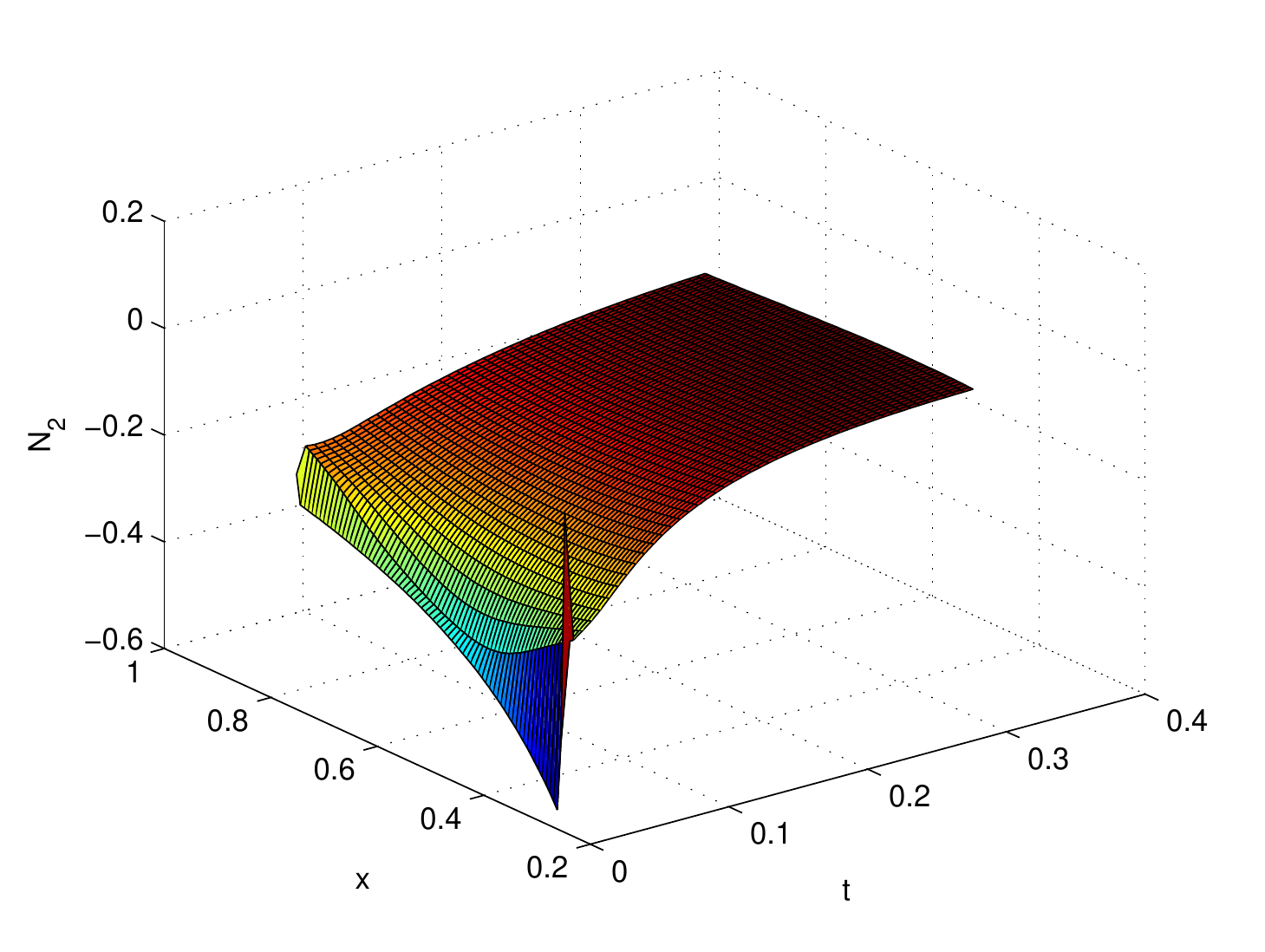}
\includegraphics[width=5.0cm,angle=-0]{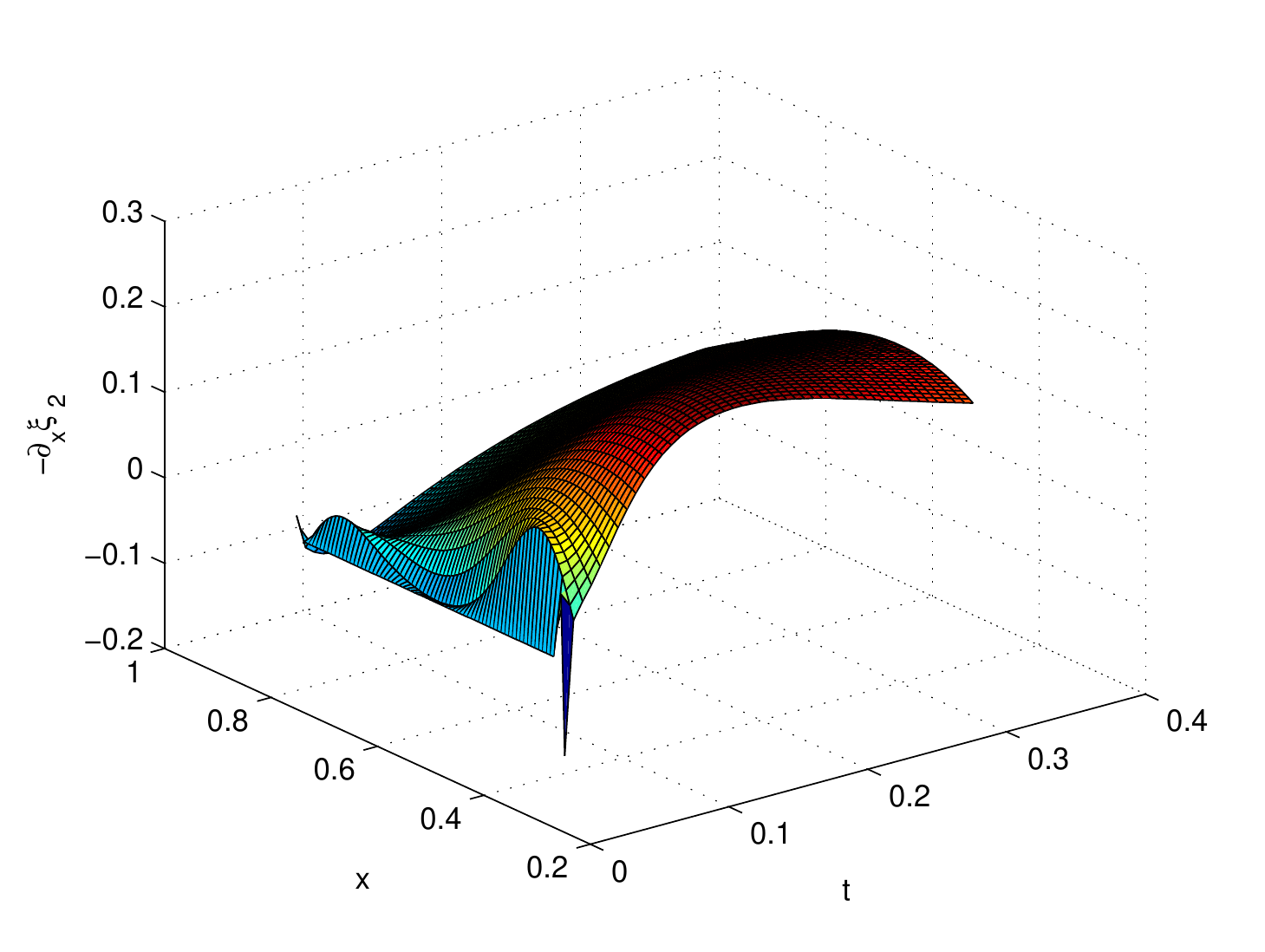}
\end{center}
\caption{\label{multi_3} The figures present the results of the 3d plots in time and space.  The upper figures present the results of the
concentration $c_1$ and $- \partial_x \xi_1$. The lower figures present
the results of $c_2$ and $- \partial_x \xi_2$.}
\end{figure}

\begin{figure}[ht]
\begin{center}  
\includegraphics[width=5.0cm,angle=-0]{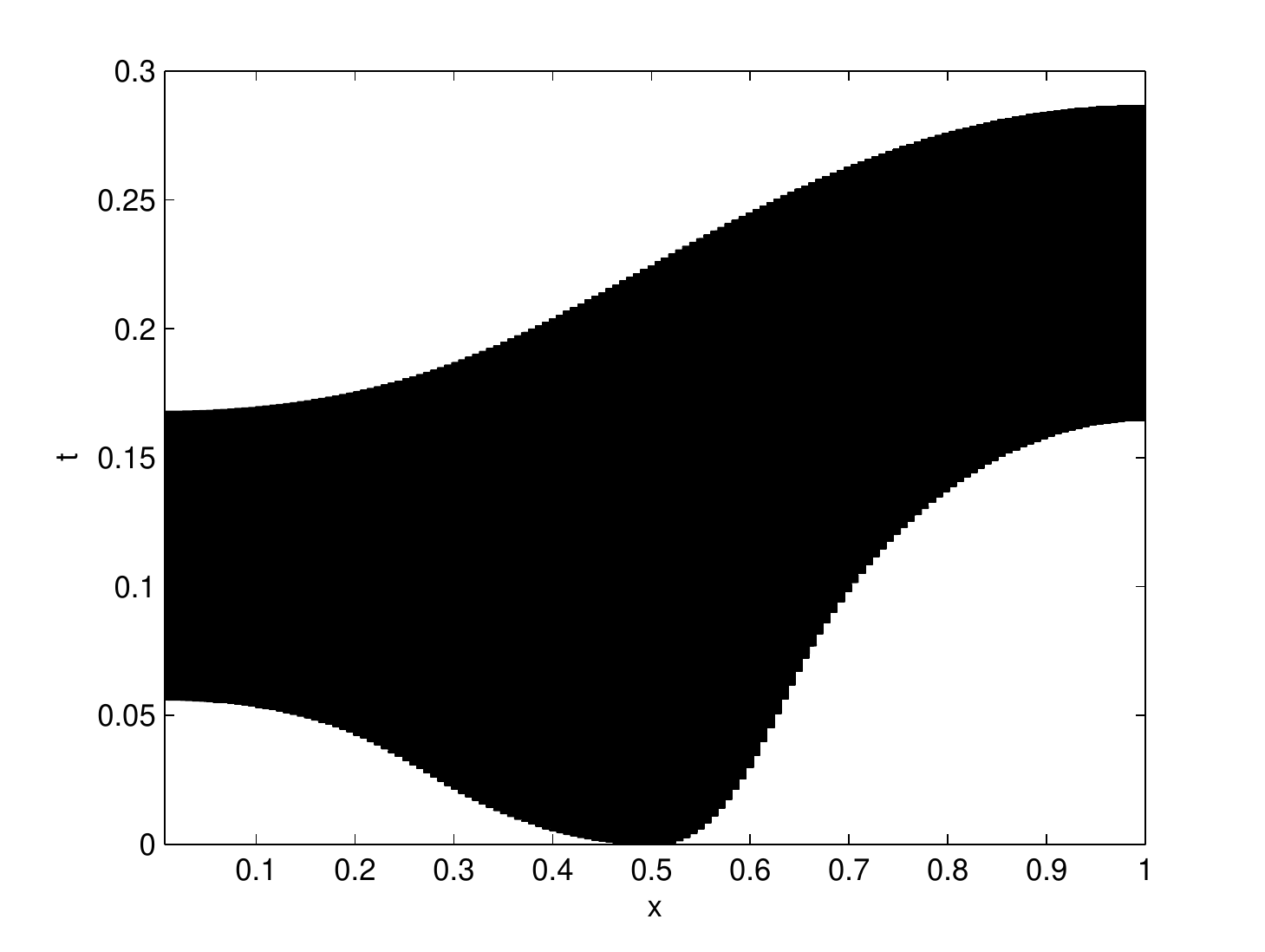}
\includegraphics[width=5.0cm,angle=-0]{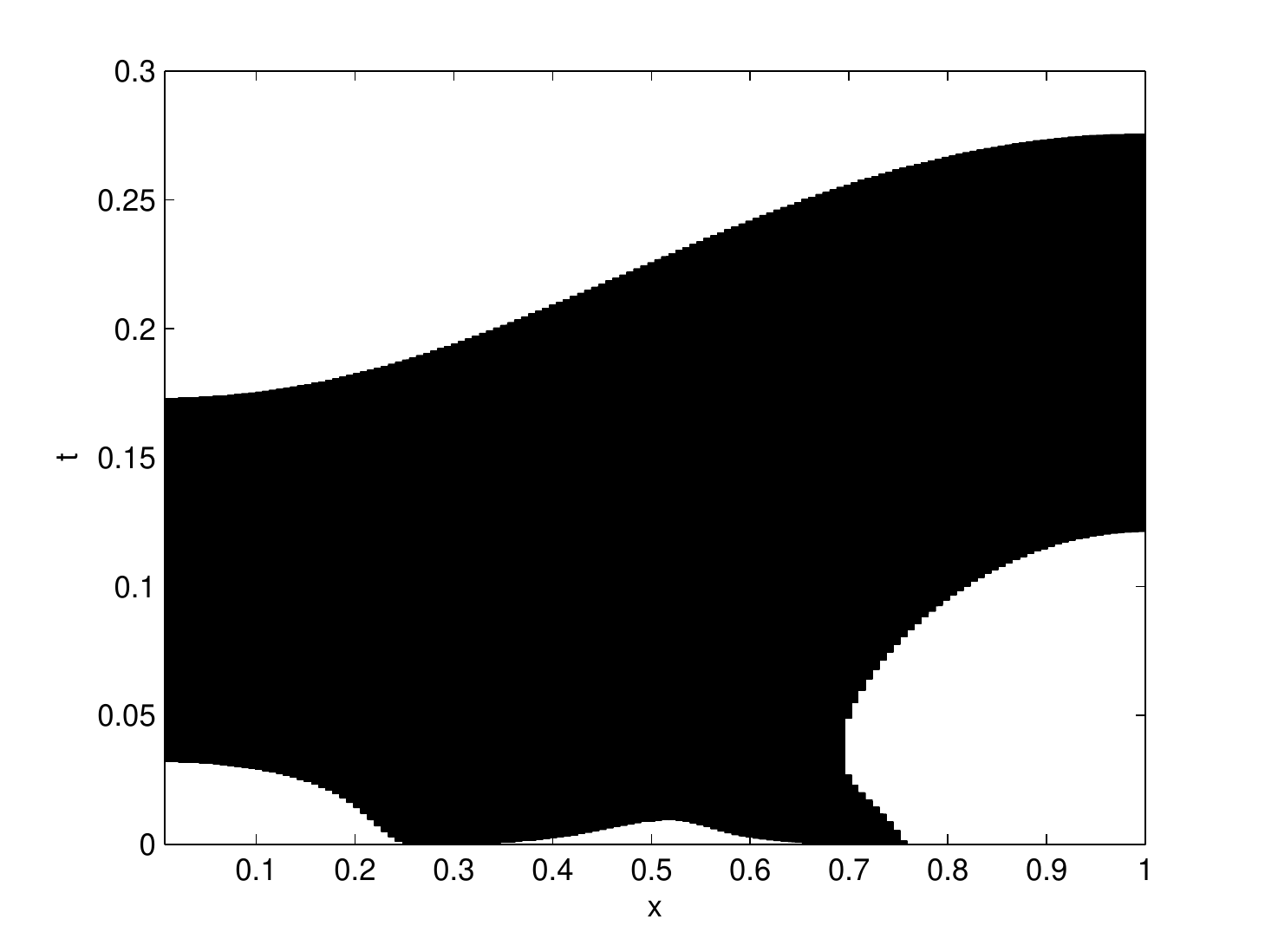} 
\end{center}
\caption{\label{multi_4} The figures present the asymptotic diffusion (left hand side) and uphill diffusion (right hand side) in the space-time region.}
\end{figure}

The space-time regions where $- N2 \partial_x \xi_2 \ge 0$ for the
uphill diffusion and asymptotic diffusion, given in Figure \ref{multi_4}.

\begin{remark}
The iterative scheme allows us to solve the pure diffusion problem effectively, see also \cite{geiser2016}. The improvement can be done with local linearisation in the pure diffusion problem, see also \cite{geiser2016} and in the next subsection.
\end{remark}

\subsection{Hydrogen Plasma: Diffusion-Reaction Problem}

In the following section we will discuss the different splitting approaches that are used to 
solve the diffusion-reaction problem.

We have explicit and implicit versions of the AB and ABA splitting approaches,
and also for the iterative splitting approach.

In the following, we have used the implicit version of the AB-splitting approach,
see Equation (\ref{ab-splitt}).

\begin{eqnarray}
\label{ab-splitt}
\ldots N^{n} \rightarrow \xi^n \rightarrow_{A} \tilde{\xi}^{n+1} \rightarrow_{B} \xi^{n+1} \rightarrow N^{n+1} \ldots
\end{eqnarray}

\begin{figure}[ht]
\begin{center}  
\includegraphics[width=8.0cm,angle=-0]{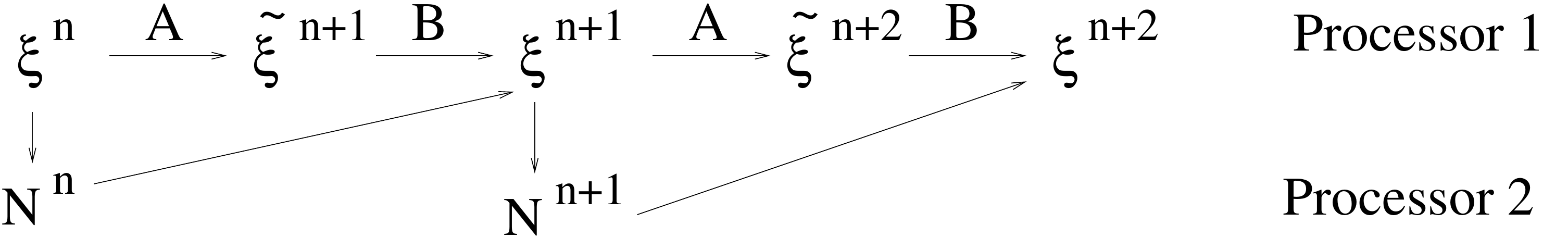} 
\end{center}
\caption{\label{ab-splitt_2} Explicit AB-splitting approach.}
\end{figure}

Furthermore, we could also apply a more explicit version of the AB-splitting
approach, which allows us to deal with a more parallel idea, see Figure \ref{ab-splitt_2}.

We concentrate on the three component system with reaction:
\begin{eqnarray}
\label{ord_0}
&& \partial_t \xi_i + \partial_x N_i = S_i , \; 1 \le i \le 3 , \\
&& \sum_{j=1}^3 N_j = 0 , \\
&& \frac{\xi_2 N_1 - \xi_1 N_2}{D_{12}} + \frac{\xi_3 N_1 - \xi_1 N_3}{D_{13}} =  -  \partial_x \xi_1 , \\
 && \frac{\xi_1 N_2 - \xi_2 N_1}{D_{12}} + \frac{\xi_3 N_2 - \xi_2 N_3}{D_{23}} =  -  \partial_x \xi_2 ,
\end{eqnarray}
where the domain is given as $\Omega \in \R^d, d \in \N^+$ with $\xi_i \in C^2$.

The parameters and the initial and boundary conditions are given as:
\begin{itemize}
\item component 1: $H_2$, component 2: $H_2^+$, component 3: $H$,
\item $D_{12} = 0.34$, $D_{13} = 0.21$ and $D_{23} = 0.21$
\item
\begin{itemize}
\item Example 1: 

$\lambda_{11} = - 4.276 \; 10^{-7}$, 
$\lambda_{21} = \lambda_{31} = - \frac{\lambda_{11}}{2}$ , \\
$\lambda_{22} = - 2.082 \; 10^{-13} $, $\lambda_{12} = \lambda_{23} = - \frac{\lambda_{22}}{2}$ , \\
$\lambda_{33} = - 4.276 \; 10^{-7}$, $\lambda_{31} = \lambda_{32} =  - \frac{\lambda_{33}}{2}$

\item Example 2:

 $\lambda_{11} = - 4.276 \; 10^{-2}$, 
$\lambda_{21} = \lambda_{31} = - \frac{\lambda_{11}}{2}$ , \\
$\lambda_{22} = - 2.082 \; 10^{-8} $, $\lambda_{12} = \lambda_{23} = - \frac{\lambda_{22}}{2}$ , \\
$\lambda_{33} = - 4.276 \; 10^{-8}$, $\lambda_{31} = \lambda_{32} =  - \frac{\lambda_{33}}{2}$

\item Example 3: 

 $\lambda_{11} = - 4.276 \; 10^{-1}$, 
$\lambda_{21} = \lambda_{31} = - \frac{\lambda_{11}}{2}$ , \\
$\lambda_{22} = - 2.082 \; 10^{-2} $, $\lambda_{12} = \lambda_{23} = - \frac{\lambda_{22}}{2}$ , \\
$\lambda_{33} = - 4.276 \; 10^{-2}$, $\lambda_{31} = \lambda_{32} =  - \frac{\lambda_{33}}{2}$

\end{itemize}
\item $J = 140$ (spatial grid points)
\item   The time-step-restriction for the explicit method is given as: \\
 $\Delta t \le (\Delta x)^2 \max \{ \frac{1}{2 \{D_{12}, D_{13}, D_{23}\}} \}$
\item The spatial domain is $\Omega = [0, 1]$ and the time-domain is $[0, T] = [0, 1]$
\item The initial conditions are:
\begin{enumerate}
\item Example uphill diff. dominant $H_2^+$:
\begin{eqnarray}
\label{init}
&& \xi_1^{in}(x) = \left\{ \begin{array}{l l}
0.8 & \mbox{if} \; 0 \le x < 0.25 , \\
1.6 (0.75 - x) & \mbox{if} \; 0.25 \le x < 0.75 , \\
0.0 & \mbox{if} \; 0.75 \le x \le 1.0 , 
\end{array} \right. ,  \\
&& \xi_1^{in}(x) = 0.2 , \; \mbox{for all} \; x \in \Omega = [0,1] ,
\end{eqnarray}
\item Example asymptotic diffusion, dominant $H_2^+$
\begin{eqnarray}
\label{init}
&& \xi_1^{in}(x) = \left\{ \begin{array}{l l}
0.8 & \mbox{if} \; 0 \le x  \in 0.5 , \\
0.0 & \mbox{else} , 
\end{array} \right. ,  \\
&&  \xi_2^{in}(x) = 0.2 ,  \; \mbox{for all} \; x \in \Omega = [0,1] ,
\end{eqnarray}
\end{enumerate}

\item The boundary conditions are of no-flux type:
\begin{eqnarray}
\label{init}
&& N_1 = N_2 = N_3 = 0 , \mbox{on} \; \partial \Omega \times [0,1] ,
\end{eqnarray}
\end{itemize}

We have used the following algorithm, which is given as AB-splitting:

\begin{algorithm}
{\bf The AB-splitting is given as:}

We start with $\xi_1(0), \xi_2(0)$ and $n = 1$:
\begin{itemize}
\item Step 1: Diffusion Step
\begin{eqnarray}
\label{part_0}
&& \partial_t \tilde{\xi_1} + \partial_x \cdot N_1 = 0 , \mbox{with} \; t \in [t^n, t^{n+1}] , \\
\label{part_1}
&& \partial_t \tilde{\xi_2} + \partial_x \cdot N_2 = 0 ,  \mbox{with} \; t \in [t^n, t^{n+1}] , \\
\label{part_2}
&& \left( \begin{array}{l}
N_1  \\
N_2 
\end{array} \right) = \frac{D_{13} D_{23}}{1 + \alpha D_{13}  \tilde{\xi}_2 + \beta D_{23}  \tilde{\xi}_1}
 \left( \begin{array}{c c}
\frac{1}{D_{23}} + \beta \tilde{\xi_1}  & \alpha \tilde{\xi_1}  \\
 \beta \tilde{\xi_2} &  \frac{1}{D_{13}}  + \alpha \tilde{\xi_2} 
\end{array} \right)
 \left( \begin{array}{l}
 -  \partial_x \tilde{\xi_1} \\
 -  \partial_x \tilde{\xi_2} 
\end{array} \right) \nonumber
\end{eqnarray}
where $\alpha = \left(\frac{1}{D_{12}} - \frac{1}{D_{13}}\right)$, 
$\beta = \left(\frac{1}{D_{12}} - \frac{1}{D_{23}}\right)$
and the initialisation is given as: $\tilde{\xi_1}(t^{n}) =\xi_1(t^{n}),\tilde{\xi_2}(t^{n}) =\xi_2(t^{n})$ (which means from the last second step).

We apply the explicit or implicit methods for the pure diffusion 
and obtain $\tilde{\xi_1}(t^{n+1}) , \tilde{\xi_2}(t^{n+1}), \tilde{\xi_3}(t^{n+1}) = 1 - \tilde{\xi_1}(t^{n+1}) - \tilde{\xi_2}(t^{n+1})$.

\item Step 2: Reaction Step
\begin{eqnarray}
\label{part_0}
 \xi_{1}(t^{n+1}) && = \tilde{\xi_1}(t^{n+1}) + \Delta t (\lambda_{11} - \lambda_{13})  \tilde{\xi}_{1}(t^{n+1}) \nonumber \\
    && + \Delta t (\lambda_{12} - \lambda_{13})  \tilde{\xi}_{2}(t^{n+1}) + \lambda_{13} , \\
 \xi_{2}(t^{n+1}) && = \tilde{\xi_2}(t^{n+1}) + \Delta t (\lambda_{21} - \Delta t \lambda_{23}) \tilde{\xi}_{1}(t^{n+1}) \nonumber \\
&&  + \Delta t (\lambda_{22} - \lambda_{23})  \tilde{\xi}_{2}(t^{n+1}) + \Delta t \lambda_{23} .
\end{eqnarray}

The solution-vectors are given as \\
 $\xi_1(t^{n+1}) = (\xi_{1,0}(t^{n+1}), \ldots, \xi_{1,J}(t^{n+1}))^t$, \\
 $\xi_2(t^{n+1}) = (\xi_{2,0}(t^{n+1}), \ldots, \xi_{2,J}(t^{n+1}))^t$, \\
 $\xi_3(t^{n+1}) = (\xi_{3,0}(t^{n+1}), \ldots, \xi_{3,J}(t^{n+1}))^t$,

\item Step 3: We go to Step 1 till $n = N$.

\end{itemize}

\end{algorithm}

\begin{figure}[ht]
\begin{center}  
\includegraphics[width=10.0cm,angle=-0]{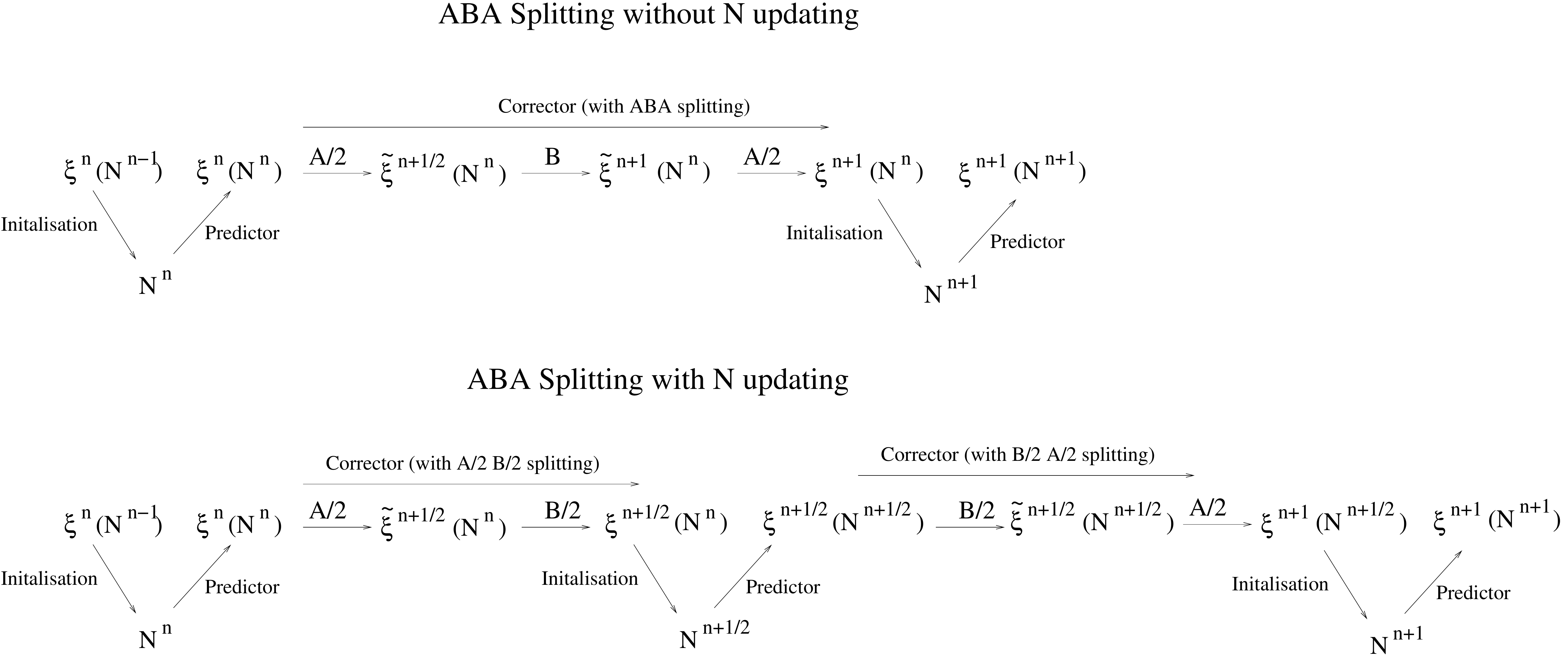} 
\end{center}
\caption{\label{aba-splitt} Explicit ABA-splitting approach with and without updating $N$.}
\end{figure}

We have used the following algorithm given as Strang-splitting in two versions,
see Algorithm \ref{algo_1_1} and Algorithm \ref{algo_1_2}.
We also explain the ideas of the splitting in Figure \ref{aba-splitt}.

\begin{algorithm}
\label{algo_1_1}

{\bf ABA-splitting (Strang-splitting) without updating $N$ is given as:}

We start with $\xi_1(0), \xi_2(0)$ and $n = 1$:
\begin{itemize}
\item Step 1: Predictor Step (updating N)
\begin{eqnarray}
\label{pred_1}
&& \left( \begin{array}{l}
N_1  \\
N_2 
\end{array} \right) = \frac{D_{13} D_{23}}{1 + \alpha D_{13}  \tilde{\xi}_2 + \beta D_{23}  \tilde{\xi}_1}
 \left( \begin{array}{c c}
\frac{1}{D_{23}} + \beta \tilde{\xi_1}  & \alpha \tilde{\xi_1}  \\
 \beta \tilde{\xi_2} &  \frac{1}{D_{13}}  + \alpha \tilde{\xi_2} 
\end{array} \right)
 \left( \begin{array}{l}
 -  \partial_x \tilde{\xi_1} \\
 -  \partial_x \tilde{\xi_2} 
\end{array} \right) \nonumber
\end{eqnarray}
where $\alpha = \left(\frac{1}{D_{12}} - \frac{1}{D_{13}}\right)$, 
$\beta = \left(\frac{1}{D_{12}} - \frac{1}{D_{23}}\right)$
and the initialisation is given as: $\tilde{\xi_1}(t^{n}) =\xi_1(t^{n}),\tilde{\xi_2}(t^{n}) =\xi_2(t^{n})$ 
(which means that this is the result of the last computation in step 3).

\item Step 2: Corrector Step (updating $\xi$)

\begin{itemize}
\item Step 2.1: Diffusion Step (with $\Delta t/2$)
\begin{eqnarray}
\label{part_0}
&& \partial_t \tilde{\xi_1} + \partial_x \cdot N_1 = 0 , \mbox{with} \; t \in [t^n, t^{n+1/2}] , \\
\label{part_1}
&& \partial_t \tilde{\xi_2} + \partial_x \cdot N_2 = 0 ,  \mbox{with} \; t \in [t^n, t^{n+1/2}] , 
\end{eqnarray}
where $\tilde{\xi_1}(t^{n}) =\xi_1(t^{n}),\tilde{\xi_2}(t^{n}) =\xi_2(t^{n})$ and $N_1, \; N_2$ is computed by the Step 1.

We apply the explicit or implicit methods for the pure diffusion 
and obtain $\tilde{\xi_1}(t^{n+1/2}) , \tilde{\xi_2}(t^{n+1/2}), and \tilde{\xi_3}(t^{n+1/2}) = 1 - \tilde{\xi_1}(t^{n+1/2}) - \tilde{\xi_2}(t^{n+1/2})$.

\item Step 2.2: Reaction Step (with $\Delta t$):
\begin{eqnarray}
\label{part_0}
 \hat{\xi}_{1}(t^{n+1}) && = \tilde{\xi_1}(t^{n+1/2}) + \Delta t (\lambda_{11} - \lambda_{13})  \tilde{\xi}_{1}(t^{n+1/2}) \nonumber \\
    && + \Delta t (\lambda_{12} - \lambda_{13})  \tilde{\xi}_{2}(t^{n+1}) + \Delta t \lambda_{13} , \\
 \hat{\xi}_{2}(t^{n+1}) && = \tilde{\xi_2}(t^{n+1/2}) + \Delta t (\lambda_{21} - \lambda_{23}) \tilde{\xi}_{1}(t^{n+1/2}) \nonumber \\
&&  + \Delta t (\lambda_{22} - \lambda_{23})  \tilde{\xi}_{2}(t^{n+1/2}) + \Delta t \lambda_{23} .
\end{eqnarray}

\item Step 2.3: Diffusion Step (with $\Delta t/2$)
\begin{eqnarray}
\label{part_0}
&& \partial_t \xi_1 + \partial_x \cdot N_1 = 0 , \mbox{with} \; t \in [t^{n+1/2}, t^{n+1}] , \\
\label{part_1}
&& \partial_t \xi_2 + \partial_x \cdot N_2 = 0 ,  \mbox{with} \; t \in [t^{n+1/2}, t^{n+1}] ,
\end{eqnarray}
where $\xi_1(t^{n+1/2}) = \hat{\xi}_1(t^{n+1}), \xi_2(t^{n+1/2}) = \hat{\xi}_2(t^{n+1})$ and $N_1, \; N_2$ is given in Step 1
(which means that $N_1(\tilde{\xi}_1(t^n)), \; N_2(\tilde{\xi}_2(t^n))$).

We apply the explicit or implicit methods for the pure diffusion 
and obtain $\xi_1(t^{n+1}) , \xi_2(t^{n+1}), \xi_2(t^{n+1}) = 1 - \xi_1(t^{n+1}) - \xi_2(t^{n+1})$.

The solution-vectors are given as \\
 $\xi_1(t^{n+1}) = (\xi_{1,0}(t^{n+1}), \ldots, \xi_{1,J}(t^{n+1}))^t$, \\
 $\xi_2(t^{n+1}) = (\xi_{2,0}(t^{n+1}), \ldots, \xi_{2,J}(t^{n+1}))^t$, \\
 $\xi_3(t^{n+1}) = (\xi_{3,0}(t^{n+1}), \ldots, \xi_{3,J}(t^{n+1}))^t$,

\end{itemize}
\item Step 3: We do $n=n+1$ and go to Step 1 till $n = N$.

\end{itemize}

\end{algorithm}

\begin{algorithm}
\label{algo_1_2}

{\bf ABA-splitting (Strang-splitting) with updating $N$ is given as:}

We start with $\xi_1(0), \xi_2(0)$ and $n = 1$:
\begin{itemize}
\item Step 1: Predictor Step (updating N)
\begin{eqnarray}
\label{pred_1}
&& \left( \begin{array}{l}
N_1  \\
N_2 
\end{array} \right) = \frac{D_{13} D_{23}}{1 + \alpha D_{13}  \tilde{\xi}_2 + \beta D_{23}  \tilde{\xi}_1}
 \left( \begin{array}{c c}
\frac{1}{D_{23}} + \beta \tilde{\xi_1}  & \alpha \tilde{\xi_1}  \\
 \beta \tilde{\xi_2} &  \frac{1}{D_{13}}  + \alpha \tilde{\xi_2} 
\end{array} \right)
 \left( \begin{array}{l}
 -  \partial_x \tilde{\xi_1} \\
 -  \partial_x \tilde{\xi_2} 
\end{array} \right) \nonumber
\end{eqnarray}
where $\alpha = \left(\frac{1}{D_{12}} - \frac{1}{D_{13}}\right)$, 
$\beta = \left(\frac{1}{D_{12}} - \frac{1}{D_{23}}\right)$
and the initialisation is given as: $\tilde{\xi_1}(t^{n}) =\xi_1(t^{n}),\tilde{\xi_2}(t^{n}) =\xi_2(t^{n})$ 
(which means that this is the result of the last computation in step 3).

\item Step 2: Corrector Step (updating $\xi$)

\begin{itemize}
\item Step 2.1: Diffusion Step (with $\Delta t/2$)
\begin{eqnarray}
\label{part_0}
&& \partial_t \tilde{\xi_1} + \partial_x \cdot N_1 = 0 , \mbox{with} \; t \in [t^n, t^{n+1/2}] , \\
\label{part_1}
&& \partial_t \tilde{\xi_2} + \partial_x \cdot N_2 = 0 ,  \mbox{with} \; t \in [t^n, t^{n+1/2}] , 
\end{eqnarray}
where $\tilde{\xi_1}(t^{n}) =\xi_1(t^{n}),\tilde{\xi_2}(t^{n}) =\xi_2(t^{n})$ and $N_1, \; N_2$ is computed by the Step 1.

We apply the explicit or implicit methods for the pure diffusion 
and obtain $\tilde{\xi_1}(t^{n+1/2}) , \tilde{\xi_2}(t^{n+1/2}), \tilde{\xi_3}(t^{n+1/2}) = 1 - \tilde{\xi_1}(t^{n+1/2}) - \tilde{\xi_2}(t^{n+1/2})$.

\item Step 2.2: Reaction Step (with $\Delta t/2$)
\begin{eqnarray}
\label{part_0}
 \xi_{1}(t^{n+1/2}) && = \tilde{\xi_1}(t^{n+1/2}) + \Delta t/2 (\lambda_{11} - \lambda_{13})  \tilde{\xi}_{1}(t^{n+1/2}) \nonumber \\
    && + \Delta t/2 (\lambda_{12} - \lambda_{13})  \tilde{\xi}_{2}(t^{n+1/2}) + \Delta t/2 \lambda_{13} , \\
 \xi_{2}(t^{n+1/2}) && = \tilde{\xi_2}(t^{n+1/2}) + \Delta t/2 (\lambda_{21} - \lambda_{23}) \tilde{\xi}_{1}(t^{n+1/2}) \nonumber \\
&&  + \Delta t/2 (\lambda_{22} - \lambda_{23})  \tilde{\xi}_{2}(t^{n+1/2}) + \Delta t/2 \lambda_{23} .
\end{eqnarray}

\end{itemize}

\item Step 3: Predictor Step (updating N)
\begin{eqnarray}
\label{pred_1}
&& \left( \begin{array}{l}
N_1  \\
N_2 
\end{array} \right) = \frac{D_{13} D_{23}}{1 + \alpha D_{13} \xi_2^{n+1/2} + \beta D_{23} \xi_1^{n+1/2}}
 \left( \begin{array}{c c}
\frac{1}{D_{23}} + \beta \xi_1^{n+1/2}  & \alpha \xi_1^{n+1/2}  \\
 \beta \xi_2^{n+1/2} &  \frac{1}{D_{13}}  + \alpha \xi_2^{n+1/2}
\end{array} \right)
 \left( \begin{array}{l}
 -  \partial_x \xi_1^{n+1/2} \\
 -  \partial_x \xi_2^{n+1/2} 
\end{array} \right) \nonumber
\end{eqnarray}
where $\alpha = \left(\frac{1}{D_{12}} - \frac{1}{D_{13}}\right)$, 
$\beta = \left(\frac{1}{D_{12}} - \frac{1}{D_{23}}\right)$
and the initialisation is given as: $\xi_1^{n+1/2} =\xi_1(t^{n+1/2}), \xi_2^{n+1/2} =\xi_2(t^{n+1/2})$ (which means that this is the result of the last computation in step 2.2).

\item Step 4: Corrector Step (updating $\xi$)

\begin{itemize}

\item Step 4.1: Reaction Step (with $\Delta t/2$)
\begin{eqnarray}
\label{part_0}
 \tilde{\xi}_{1}(t^{n+1}) && = \xi_1(t^{n+1/2}) + \Delta t/2 (\lambda_{11} - \lambda_{13})  \xi_{1}(t^{n+1/2}) \nonumber \\
    && + \Delta t/2 (\lambda_{12} - \lambda_{13})  \xi_{2}(t^{n+1/2}) + \Delta t/2 \lambda_{13} , \\
 \tilde{\xi}_{2}(t^{n+1}) && = \xi_2(t^{n+1/2}) + \Delta t/2 (\lambda_{21} - \lambda_{23}) \xi_{1}(t^{n+1/2}) \nonumber \\
&&  + \Delta t/2 (\lambda_{22} - \lambda_{23})  \xi_{2}(t^{n+1/2}) + \Delta t/2 \lambda_{23} .
\end{eqnarray}

\item Step 4.2: Diffusion Step (with $\Delta t/2$)
\begin{eqnarray}
\label{part_0}
&& \partial_t \xi_1 + \partial_x \cdot N_1 = 0 , \mbox{with} \; t \in [t^{n+1/2}, t^{n+1}] , \\
\label{part_1}
&& \partial_t \xi_2 + \partial_x \cdot N_2 = 0 ,  \mbox{with} \; t \in [t^{n+1/2}, t^{n+1}] ,
\end{eqnarray}
where $\xi_1(t^{n+1/2}) = \tilde{\xi}_1(t^{n+1}), \xi_2(t^{n+1/2}) = \tilde{\xi}_2(t^{n+1})$ and $N_1, \; N_2$ is given in the
updated Step 3
(which means that $N_1(\xi_1(t^{n+1/2})), \; N_2(\xi_2(t^{n+1/2}))$).

We apply the explicit or implicit methods for the pure diffusion 
and obtain $\xi_1(t^{n+1}) , \xi_2(t^{n+1}), \xi_2(t^{n+1}) = 1 - \xi_1(t^{n+1}) - \xi_2(t^{n+1})$.

The solution-vectors are given as \\
 $\xi_1(t^{n+1}) = (\xi_{1,0}(t^{n+1}), \ldots, \xi_{1,J}(t^{n+1}))^t$, \\
 $\xi_2(t^{n+1}) = (\xi_{2,0}(t^{n+1}), \ldots, \xi_{2,J}(t^{n+1}))^t$, \\
 $\xi_3(t^{n+1}) = (\xi_{3,0}(t^{n+1}), \ldots, \xi_{3,J}(t^{n+1}))^t$,

\end{itemize}
\item Step 5: We do $n=n+1$ and go to Step 1 till $n = N$.

\end{itemize}

\end{algorithm}

\newpage

We have used the following algorithm, given as an iterative splitting approach,
while we solve the diffusion part and perturb over the reaction part:

\begin{algorithm}
{\bf The iterative splitting for reaction (Picard's fixpoint scheme) is given as:}

We start with $\xi_1(0), \xi_2(0)$ and $n = 1$:

\begin{itemize}

\item Step 0: Initialisation for $i = 0$ with  
$\xi_{1,0}(t^{n+1}) = \xi_{1}(t^n) , \xi_{2,0}(t^{n+1}) = \xi_{2}(t^n)$ \\
and $N_{1,0}(t^{n+1}) = N_{1}(t^n) , N_{2,0}(t^{n+1}) = N_{2}(t^n)$

\item Step 1: Iterative step $i$: Diffusion and Reaction Step (with $\Delta t$)

\begin{itemize}
\item Step 1.1. Computation of $\xi_{1, i}^{n+1}$ and $\xi_{2, i}^{n+1}$
\begin{eqnarray}
\label{part_0}
&& \xi_{1, i}^{n+1} = \xi_{1, i}^{n} - \Delta t D_+ N_{1,i-1}^{n+1} + \nonumber \\
&& + \Delta t \;  \left( \lambda_{11}  \xi_{1, i-1}^{n+1} + \lambda_{12}  \xi_{2, i-1}^{n+1} + \lambda_{13} (1 -  \xi_{1, i-1}^{n+1} - \xi_{2, i-1}^{n+1} ) \right) ,  \\
\label{part_1}
&& \xi_{2, i}^{n+1} = \xi_{2, i}^{n} - \Delta t D_+ N_{2,i-1}^{n+1} + \nonumber \\
&& + \Delta t \; \left( \lambda_{21}  \xi_{1, i-1}^{n+1} + \lambda_{22}  \xi_{2, i-1}^{n+1} + \lambda_{23} (1 -  \xi_{1, i-1}^{n+1} - \xi_{2, i-1}^{n+1} ) \right) ,
\end{eqnarray}

\item Step 1.2: $i = i + 1$ and we go to Step 1 till $i = I$ (else goto Step 2) 

\end{itemize}

\item Step 2. Computation of $N_{1, i}^{n+1}$ and $N_{2, i}^{n+1}$
\begin{eqnarray}
\label{part_2}
&& \left( \begin{array}{l}
N_{1,i}^{n+1}  \\
N_{2,i}^{n+1} 
\end{array} \right) = \\
&&  = \frac{D_{13} D_{23}}{1 + \alpha D_{13} \xi_{2,i}^{n+1} + \beta D_{23} \xi_{1,i}^{n+1}}
 \left( \begin{array}{c c}
\frac{1}{D_{23}} + \beta \xi_{1,i}^{n+1}  & \alpha \xi_{1,i}^{n+1}  \\
 \beta \xi_{2,i}^{n+1} &  \frac{1}{D_{13}}  + \alpha \xi_{2,i}^{n+1} 
\end{array} \right)
 \left( \begin{array}{l}
 -  \partial_x \xi_{1,i}^{n+1} \\
 -  \partial_x \xi_{2,i}^{n+1} 
\end{array} \right) \nonumber
\end{eqnarray}

where $\alpha = \left(\frac{1}{D_{12}} - \frac{1}{D_{13}}\right)$, 
$\beta = \left(\frac{1}{D_{12}} - \frac{1}{D_{23}}\right)$
and the initialisation is given as: $\xi_{1,i}(t^{n}) = \xi_1(t^{n}), \xi_{2,i}(t^{n}) =\xi_2(t^{n})$ (which is the means from the last computation).

\item Step 4: $n = n+1$ and we go to Step 0 till $n = N$.

\end{itemize}

We apply the explicit or implicit methods for the diffusion-reaction equation 
and obtain $\xi_{1,i}(t^{n+1}) , \xi_{2,i}(t^{n+1}), \xi_{3,i}(t^{n+1}) = 1 - \xi_{1,i}(t^{n+1}) - \xi_{2,i}(t^{n+1})$.

\end{algorithm}

\begin{algorithm}
{\bf The iterative splitting for diffusion and reaction (Inner and outer Picard's fixpoint scheme) is given as:}

We start with $\xi_1(0), \xi_2(0)$ and $n = 1$:

\begin{itemize}

\item Step 0: Initialisation for $i,j = 0$ with  
$\xi_{1,0}(t^{n+1}) = \xi_{1}(t^n) , \xi_{2,0}(t^{n+1}) = \xi_{2}(t^n)$ \\
and $N_{1,0}(t^{n+1}) = N_{1}(t^n) , N_{2,0}(t^{n+1}) = N_{2}(t^n)$.
We have $i = j = 1$ (initialisation of the loops).

\item Step 1: Outer Loop (Iterative step $j$): Diffusion Step, Computation of $N_{1, j}^{n+1}$ and $N_{2, j}^{n+1}$ (with $\Delta t$)

\begin{itemize}
\item Step 1.1.: Inner Loop (Iterative step $i$): Reaction Step, Computation of $\xi_{1, i}^{n+1}$ and $\xi_{2, i}^{n+1}$
\begin{eqnarray}
\label{part_0}
&& \xi_{1, i}^{n+1} = \xi_{1, i}^{n} - \Delta t D_+ N_{1,j-1}^{n+1} + \nonumber \\
&& + \Delta t \;  \left( \lambda_{11}  \xi_{1, i-1}^{n+1} + \lambda_{12}  \xi_{2, i-1}^{n+1} + \lambda_{13} (1 -  \xi_{1, i-1}^{n+1} - \xi_{2, i-1}^{n+1} ) \right) ,  \\
\label{part_1}
&& \xi_{2, i}^{n+1} = \xi_{2, i}^{n} - \Delta t D_+ N_{2,j-1}^{n+1} + \nonumber \\
&& + \Delta t \; \left( \lambda_{21}  \xi_{1, i-1}^{n+1} + \lambda_{22}  \xi_{2, i-1}^{n+1} + \lambda_{23} (1 -  \xi_{1, i-1}^{n+1} - \xi_{2, i-1}^{n+1} ) \right) ,
\end{eqnarray}

\item Step 1.2.: $i = i + 1$ and we go to Step 1.1. till $i = I * j $ (else goto Step 2) 

\end{itemize}

\item Step 2. Computation of $N_{1, j}^{n+1}$ and $N_{2, j}^{n+1}$
\begin{eqnarray}
\label{part_2}
&& \left( \begin{array}{l}
N_{1,j}^{n+1}  \\
N_{2,j}^{n+1} 
\end{array} \right) = \\
&&  = \frac{D_{13} D_{23}}{1 + \alpha D_{13} \xi_{2,i}^{n+1} + \beta D_{23} \xi_{1,I}^{n+1}}
 \left( \begin{array}{c c}
\frac{1}{D_{23}} + \beta \xi_{1,I}^{n+1}  & \alpha \xi_{1,i}^{n+1}  \\
 \beta \xi_{2,i}^{n+1} &  \frac{1}{D_{13}}  + \alpha \xi_{2,i}^{n+1} 
\end{array} \right)
 \left( \begin{array}{l}
 -  \partial_x \xi_{1,i}^{n+1} \\
 -  \partial_x \xi_{2,i}^{n+1} 
\end{array} \right) \nonumber
\end{eqnarray}

where $\alpha = \left(\frac{1}{D_{12}} - \frac{1}{D_{13}}\right)$, 
$\beta = \left(\frac{1}{D_{12}} - \frac{1}{D_{23}}\right)$
and the initialisation is given as: $\xi_{1,i}(t^{n}) = \xi_1(t^{n}), \xi_{2,i}(t^{n}) =\xi_2(t^{n})$ (means from the last computation).

\item Step 3: $j = j+1$ and we go to Step 1 till $j = J$.

\item Step 4: $n = n+1$ and we go to Step 0 till $n = N$.

\end{itemize}

We apply the explicit or implicit methods for the diffusion-reaction equation 
and obtain $\xi_{1,i}(t^{n+1}) , \xi_{2,i}(t^{n+1}), \xi_{3,i}(t^{n+1}) = 1 - \xi_{1,i}(t^{n+1}) - \xi_{2,i}(t^{n+1})$.

\end{algorithm}

For a run, we assume that $I = J = 2$, which means that we have two iterative loops in the
inner and two in the outer.
For the convergence threshold, we define the variance between a
reference solution and the numerical solutions, given as:
 Time-averaged mean-square value over time (scan over the time-space):
\begin{eqnarray}
 \sigma^2_{\xi_1, \Delta t} = \frac{1}{T} \sum_{i=1}^N \Delta t \; ( \xi_{1, \Delta t, Scheme}(i \; \Delta t) - \xi_{1, \Delta t, ref}(i \; \Delta t) )^2 .
\end{eqnarray}
\begin{eqnarray}
 \sigma^2_{\xi_2, \Delta t} = \frac{1}{T} \sum_{i=1}^N \Delta t \; ( \xi_{2, \Delta t, Scheme}(i \; \Delta t) - \xi_{2, \Delta t, ref}(i \; \Delta t) )^2 ,
\end{eqnarray}
where the time-space is given as $i=1, \ldots, N$, $\Delta t \; N = T = 1$.

Furthermore, the vectorial time-averaged means square value is:
\begin{eqnarray}
 \sigma^2_{\xi, \Delta t} && = \frac{1}{T} \sum_{i=1}^N \Delta t \; \left( ( \xi_{1, \Delta t, Scheme}(i \; \Delta t) - \xi_{1, \Delta t, ref}(i \; \Delta t) )^2  \right. \\
&& \left. +  ( \xi_{2, \Delta t, Scheme}(i \; \Delta t) - \xi_{2, \Delta t, ref}(i \; \Delta t) )^2 \right) , \nonumber
\end{eqnarray}
where the time-space is given as $i=1, \ldots, N$, $\Delta t \; N = T = 1$.

\begin{itemize}
\item Example 1:

$\lambda_{11} = - 4.276 \; 10^{-7}$, 
$\lambda_{21} = \lambda_{31} = - \frac{\lambda_{11}}{2}$ , \\
$\lambda_{22} = - 2.082 \; 10^{-13} $, $\lambda_{12} = \lambda_{23} = - \frac{\lambda_{22}}{2}$ , \\
$\lambda_{33} = - 4.276 \; 10^{-7}$, $\lambda_{31} = \lambda_{32} =  - \frac{\lambda_{33}}{2}$

The numerical solutions of the three hydrogen plasma in experiment 1 with the
asymptotic diffusion \ref{hydrogen_1_1}.
\begin{figure}[ht]
\begin{center}  
\includegraphics[width=5.0cm,angle=-0]{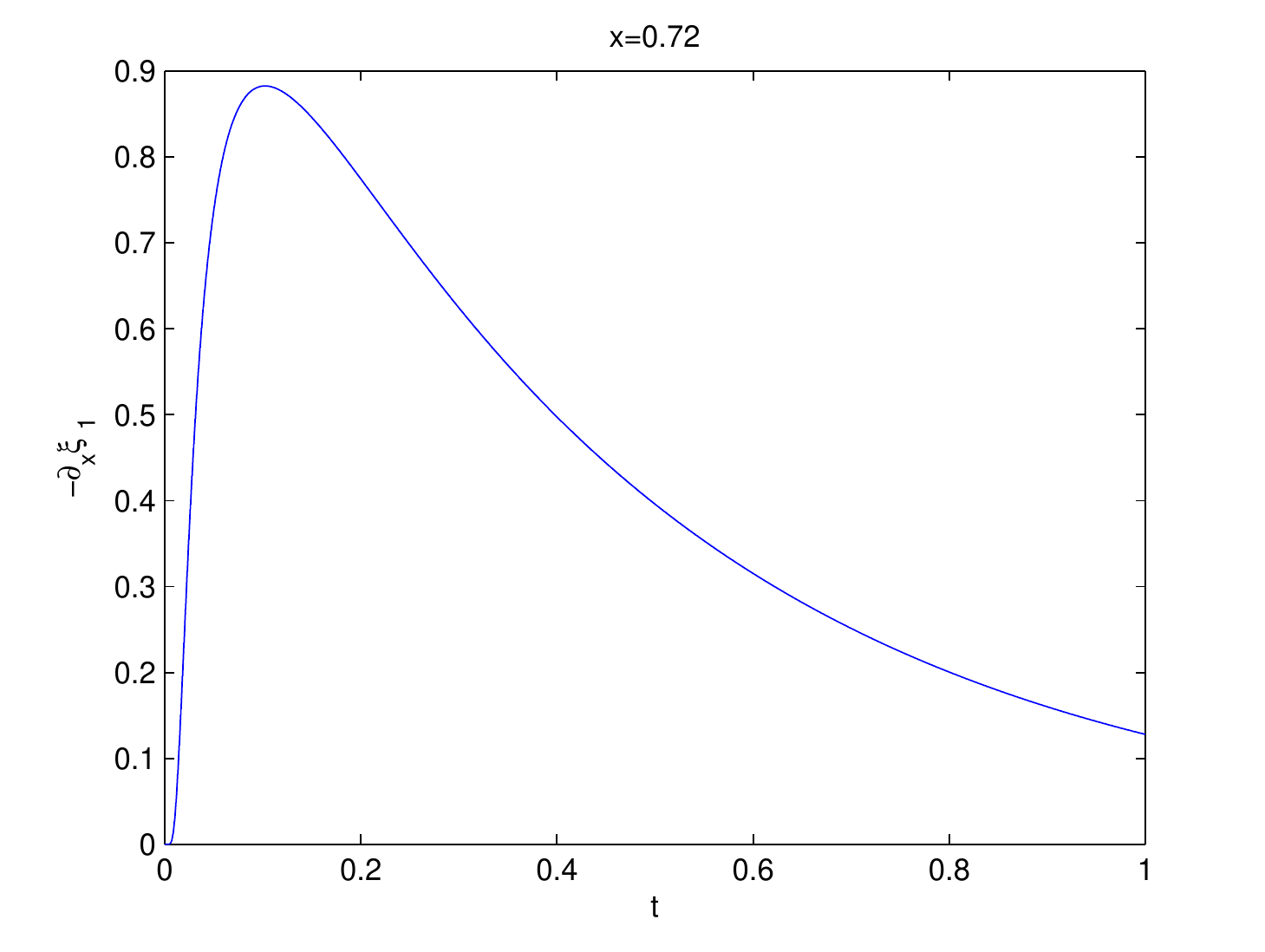} 
\includegraphics[width=5.0cm,angle=-0]{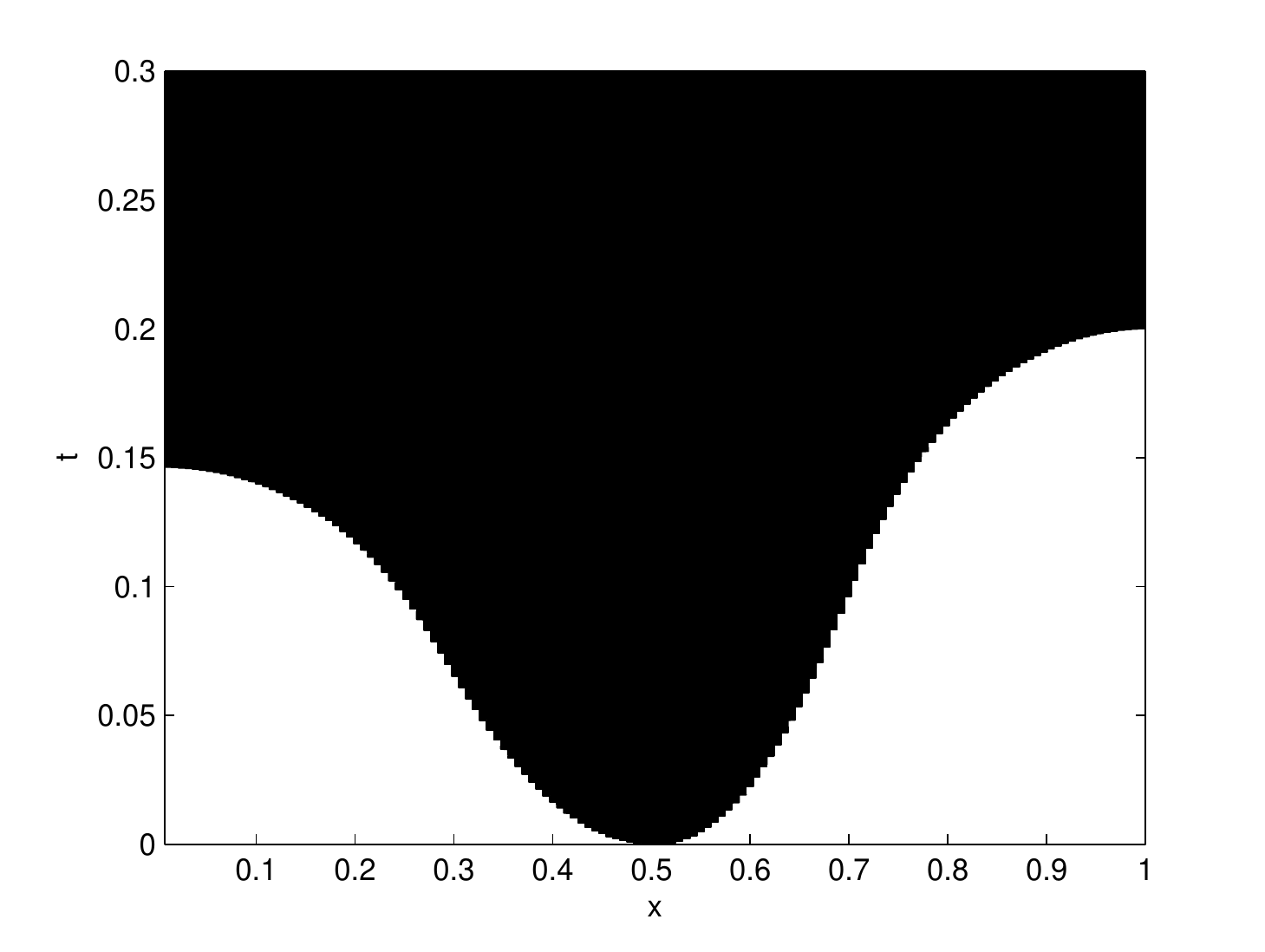} \\
\includegraphics[width=5.0cm,angle=-0]{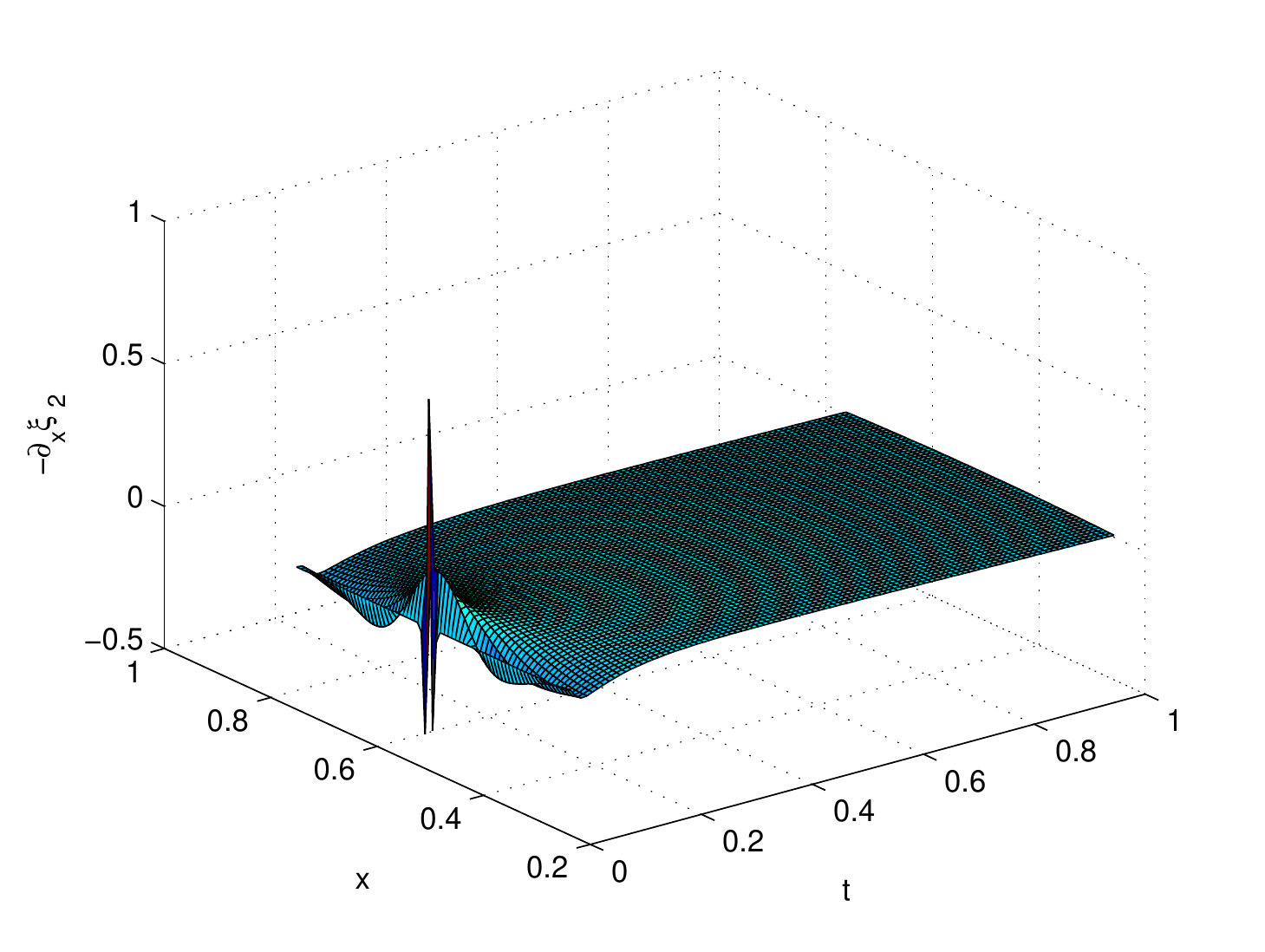} 
\includegraphics[width=5.0cm,angle=-0]{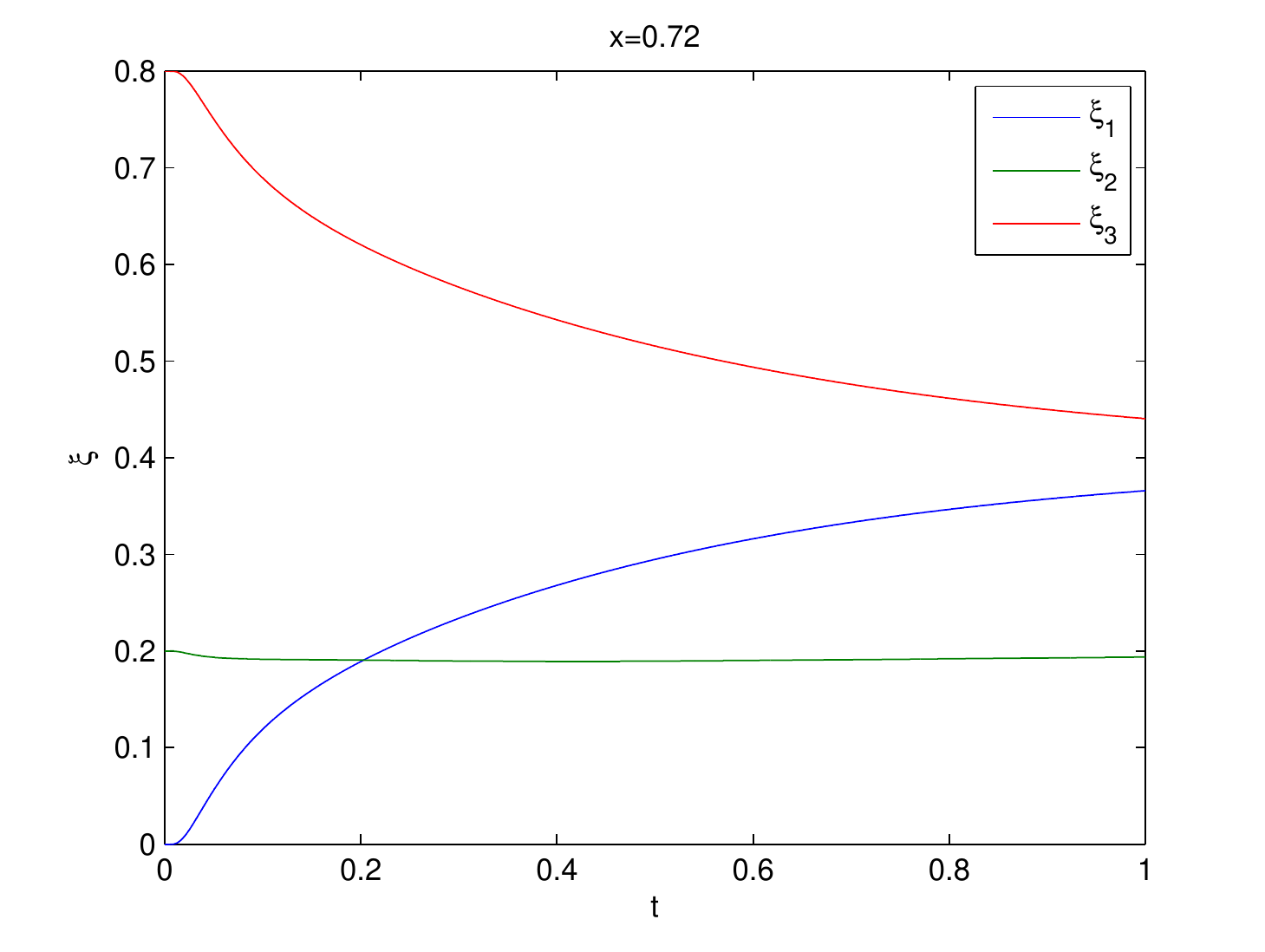} 
\end{center}
\caption{\label{hydrogen_1_1} The upper left figure presents the concentration of $x_1$ at spatial point $0.72$, the upper right is the result in the space time region, the lower left figure presents the the 3D plot of the second component and the lower right figure presents all of the components at spatial-point $0.72$.}
\end{figure}

The numerical solutions of the three hydrogen plasma in experiment 1 with the
uphill diffusion \ref{hydrogen_1_2}.
\begin{figure}[ht]
\begin{center}  
\includegraphics[width=5.0cm,angle=-0]{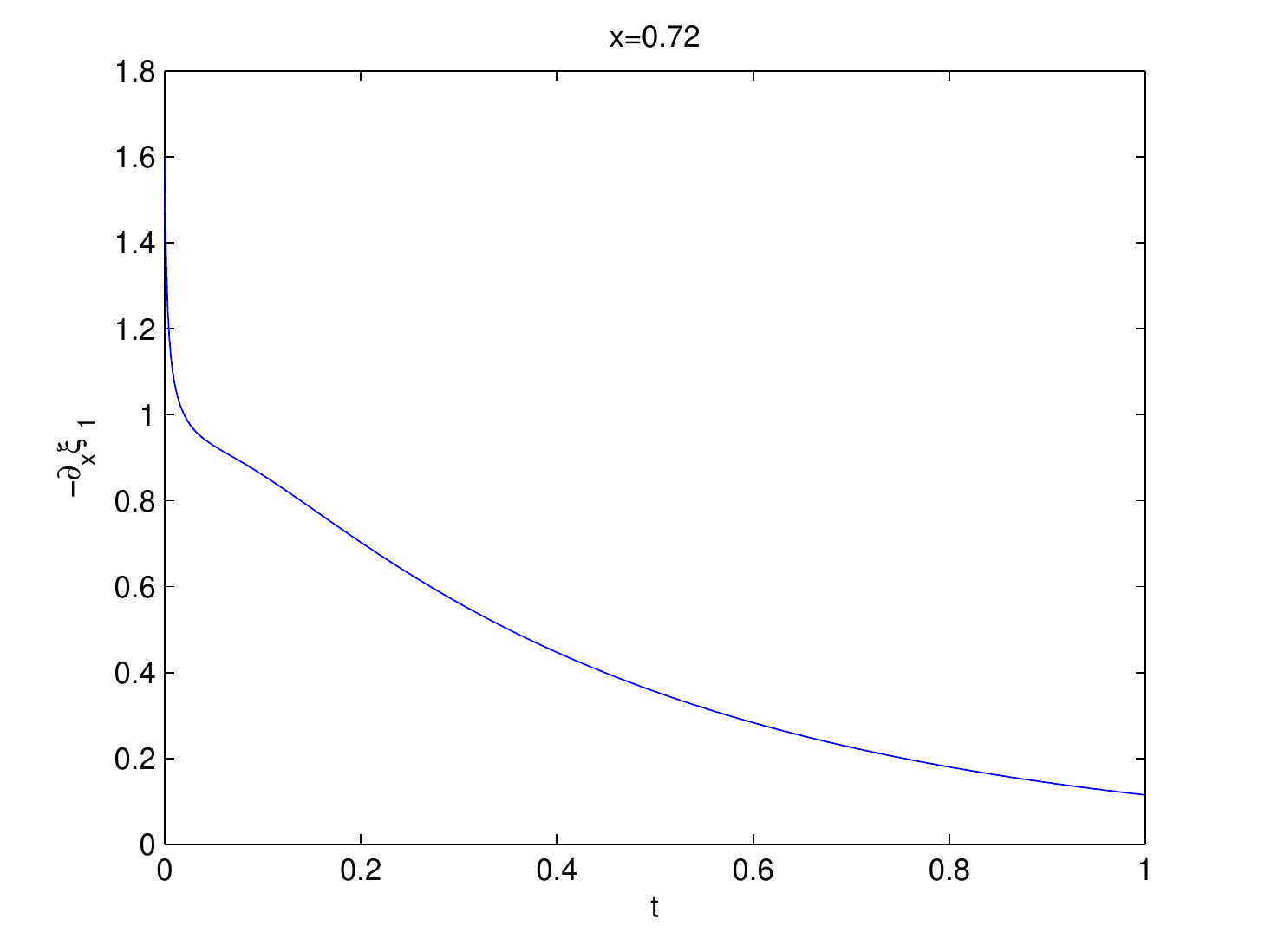} 
\includegraphics[width=5.0cm,angle=-0]{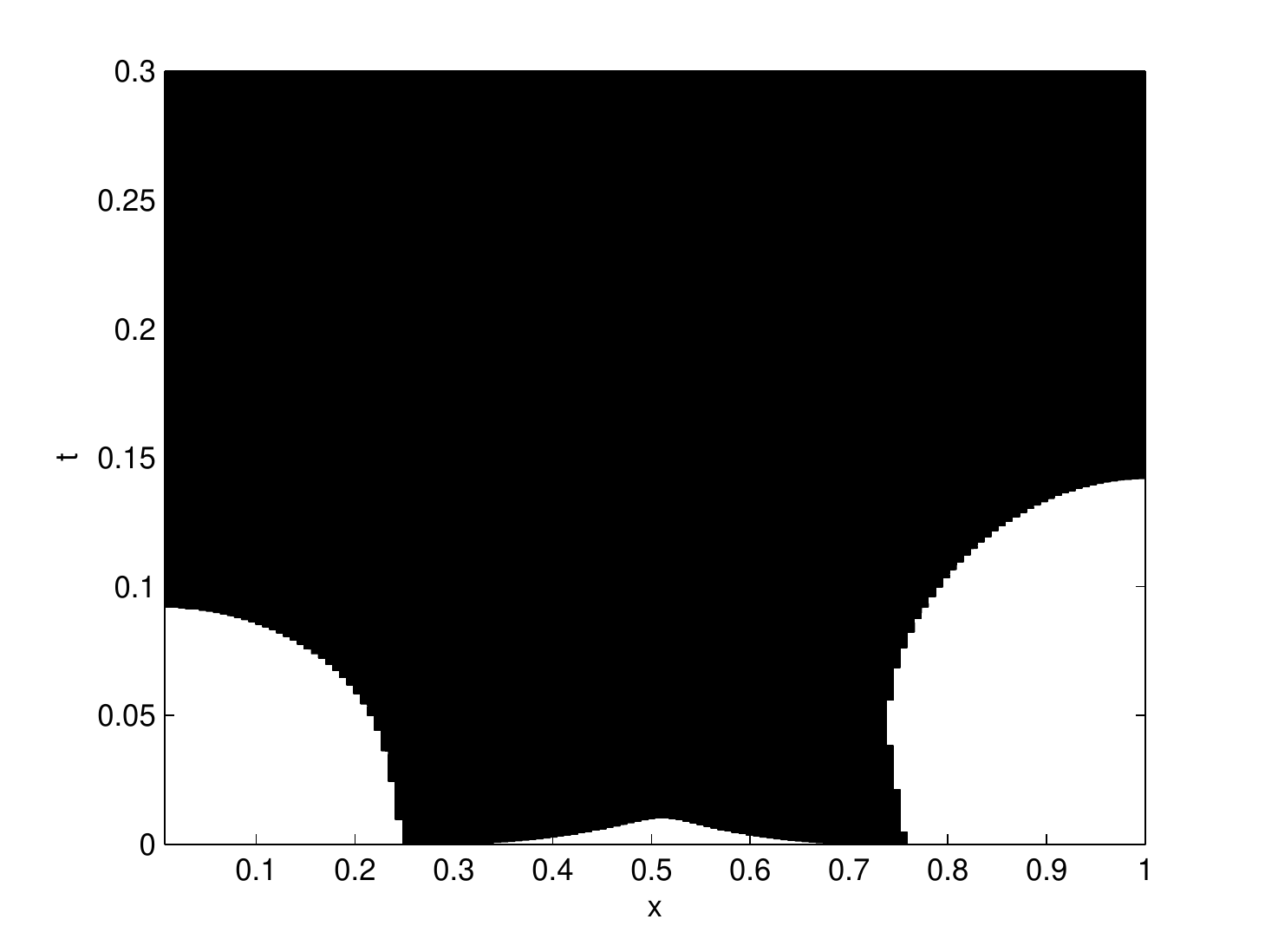} \\
\includegraphics[width=5.0cm,angle=-0]{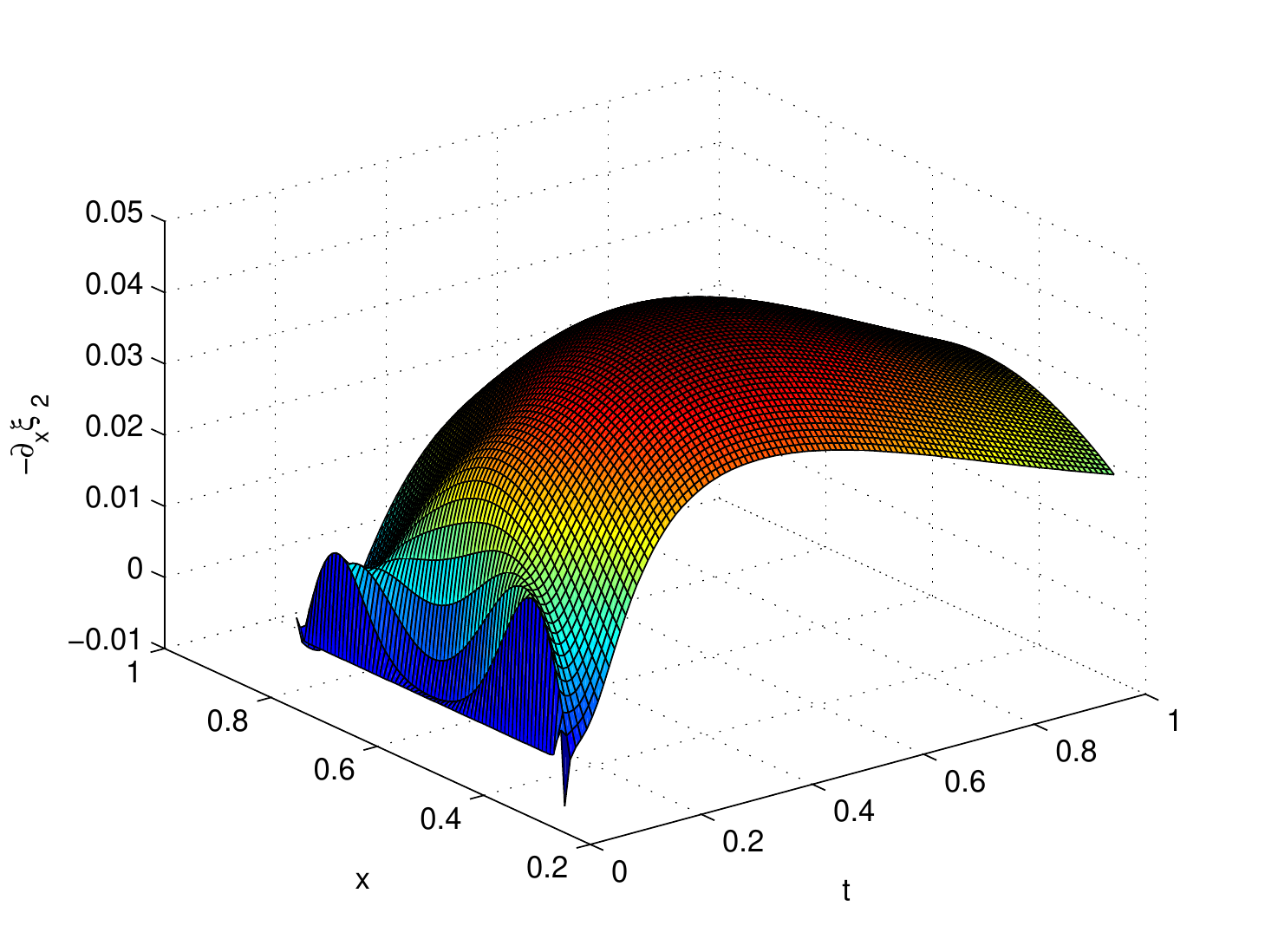} 
\includegraphics[width=5.0cm,angle=-0]{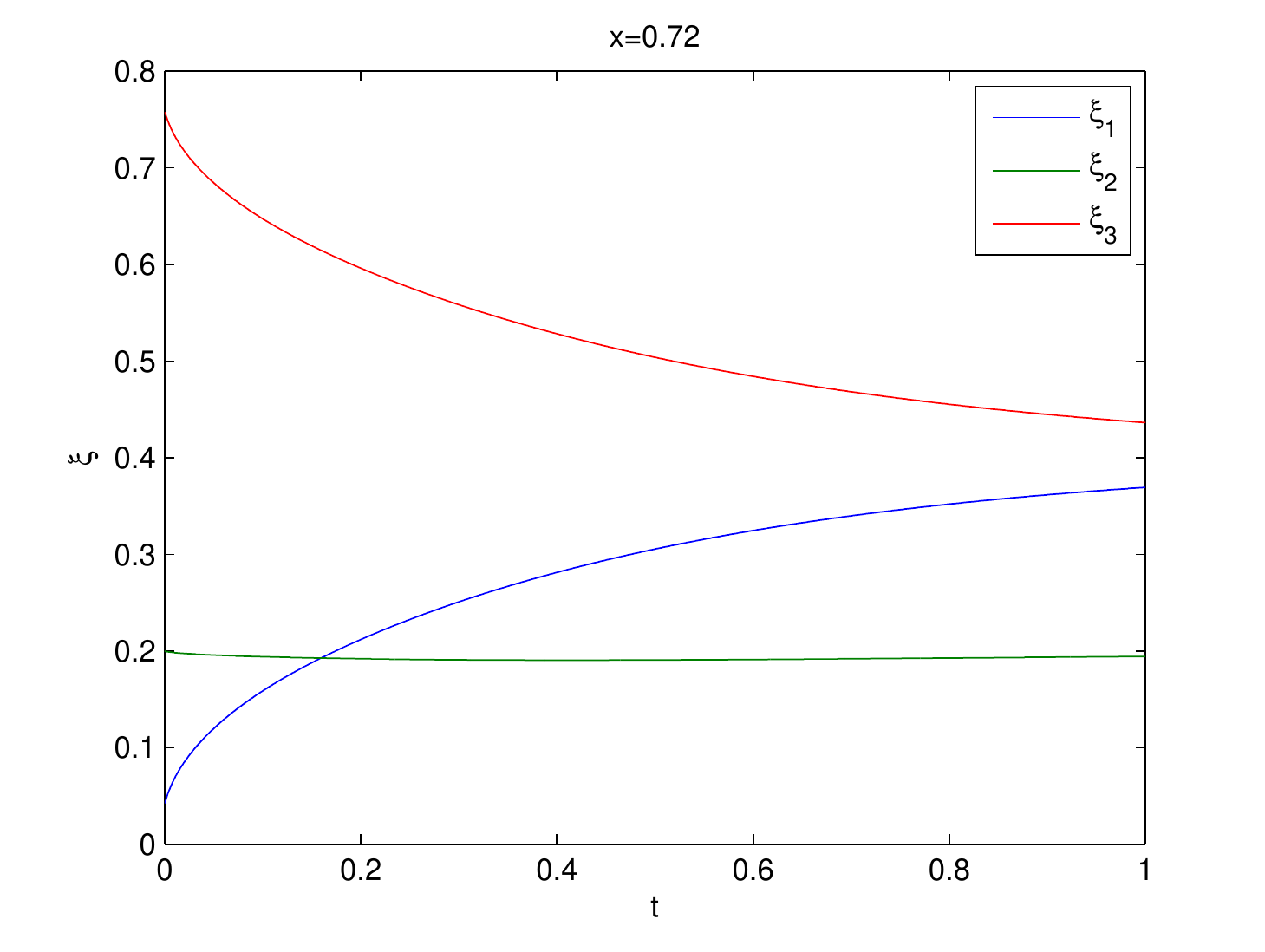} 
\end{center}
\caption{\label{hydrogen_1_2} The upper left figure presents the concentration of $x_1$ at spatial point $0.72$, the upper right is the result in the space time region, the lower left figure presents the the 3D plot of the second component and the lower right figure presents all of the components at spatial-point $0.72$.}
\end{figure}

\item Example 2:

 $\lambda_{11} = - 4.276 \; 10^{-2}$, 
$\lambda_{21} = \lambda_{31} = - \frac{\lambda_{11}}{2}$ , \\
$\lambda_{22} = - 2.082 \; 10^{-8} $, $\lambda_{12} = \lambda_{23} = - \frac{\lambda_{22}}{2}$ , \\
$\lambda_{33} = - 4.276 \; 10^{-8}$, $\lambda_{31} = \lambda_{32} =  - \frac{\lambda_{33}}{2}$

The numerical solutions of the three hydrogen plasma in experiment 2 with the
asymptotic diffusion \ref{hydrogen_2_1}.
\begin{figure}[ht]
\begin{center}  
\includegraphics[width=5.0cm,angle=-0]{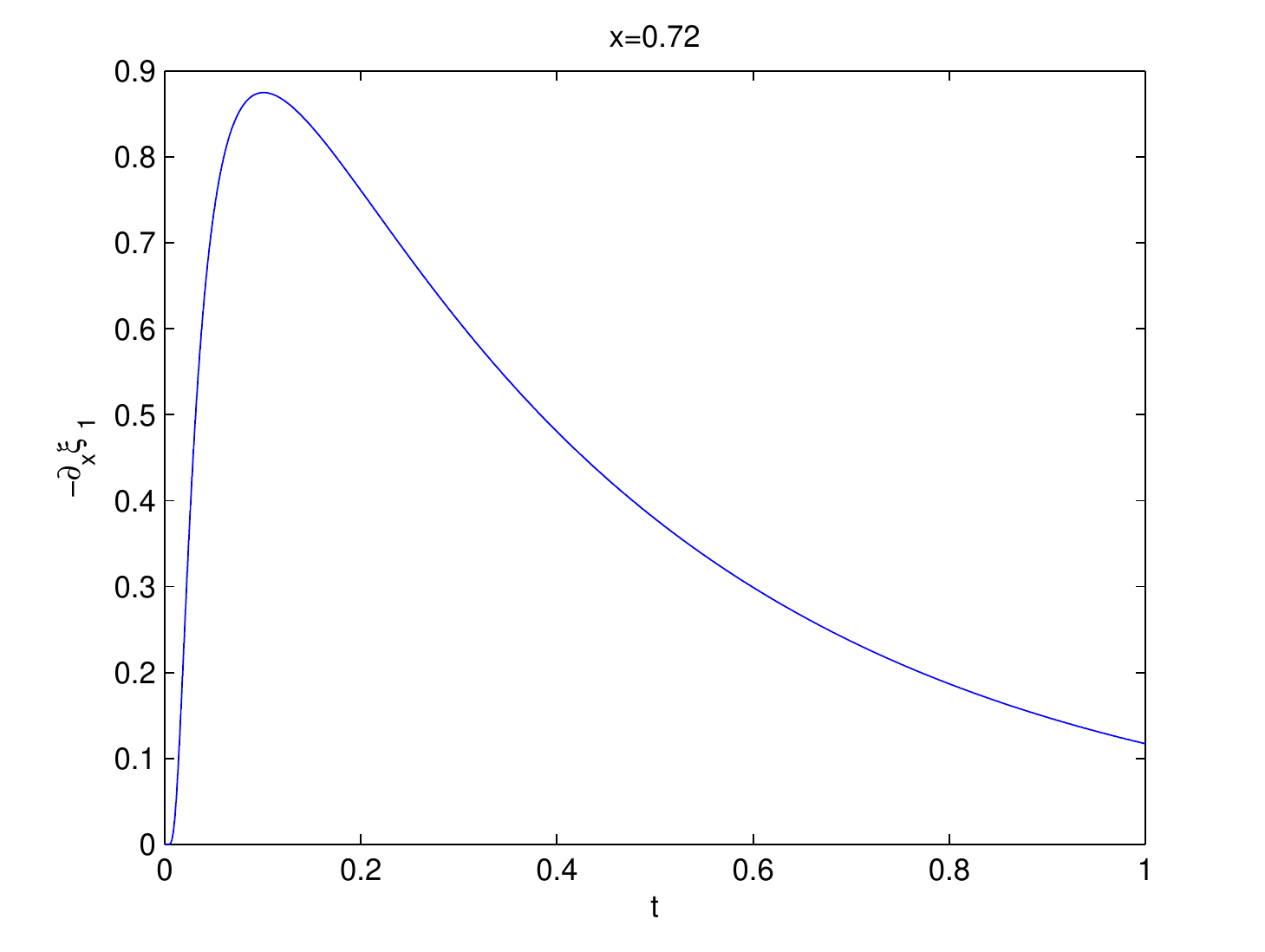} 
\includegraphics[width=5.0cm,angle=-0]{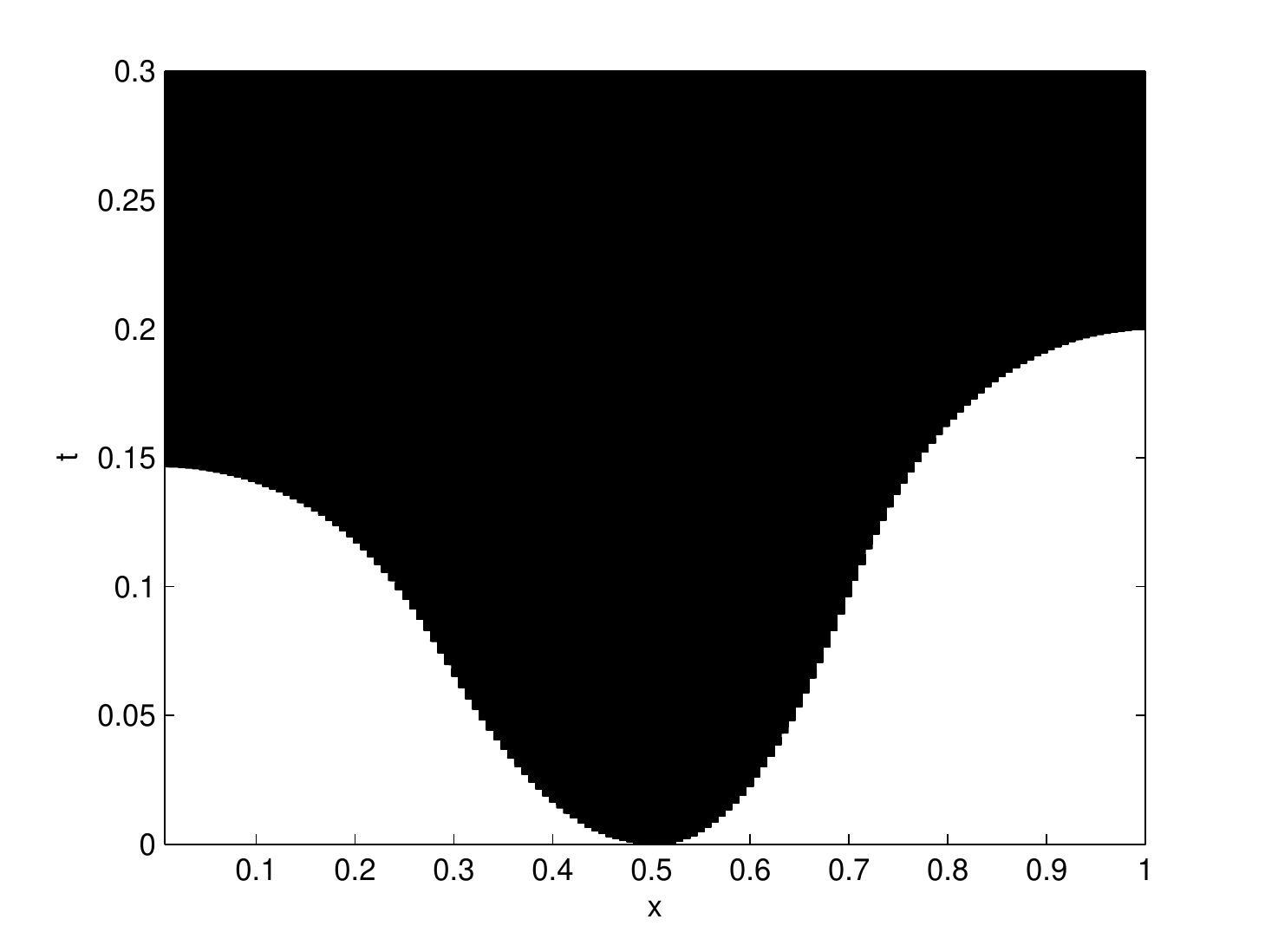} \\
\includegraphics[width=5.0cm,angle=-0]{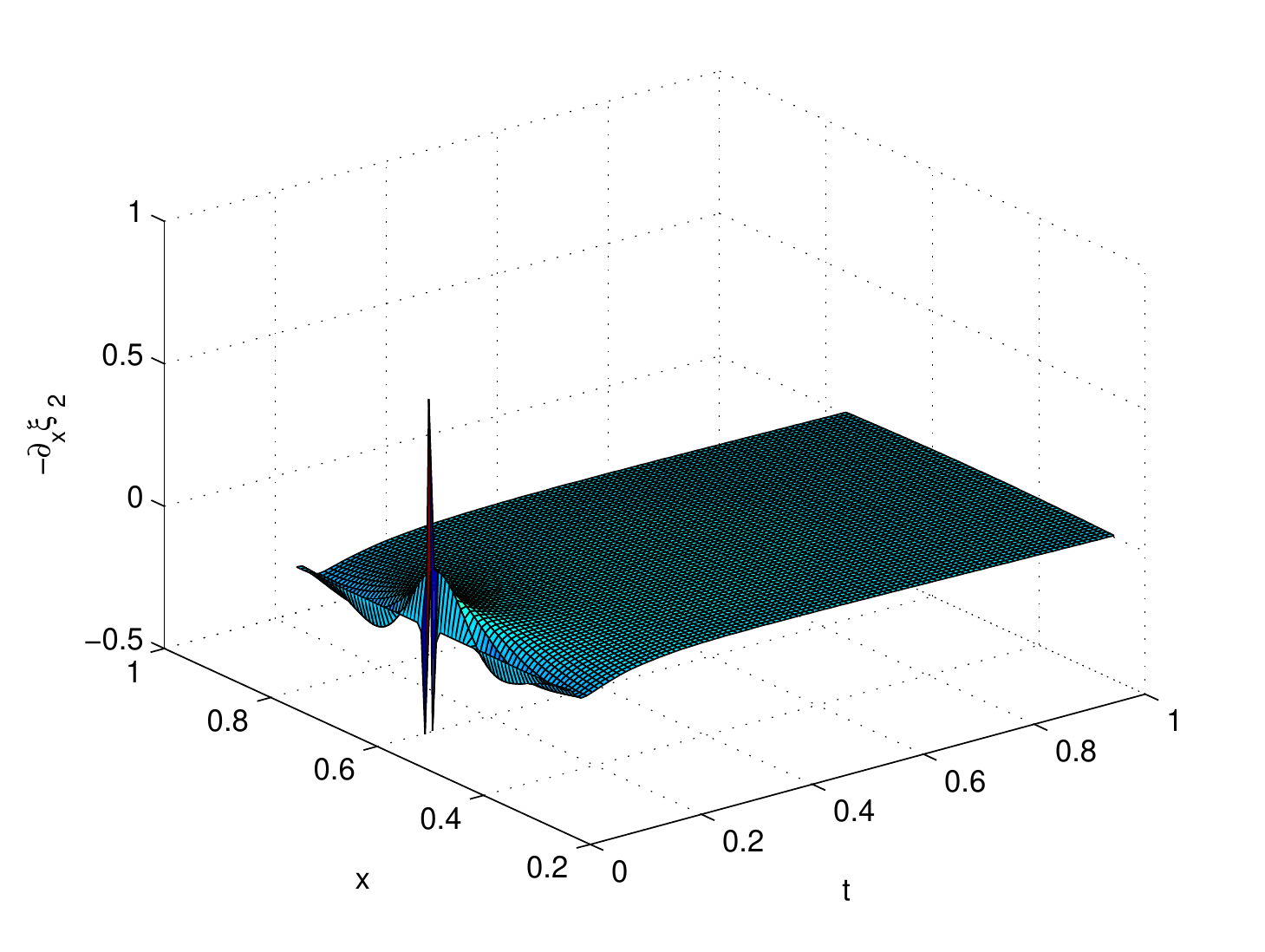} 
\includegraphics[width=5.0cm,angle=-0]{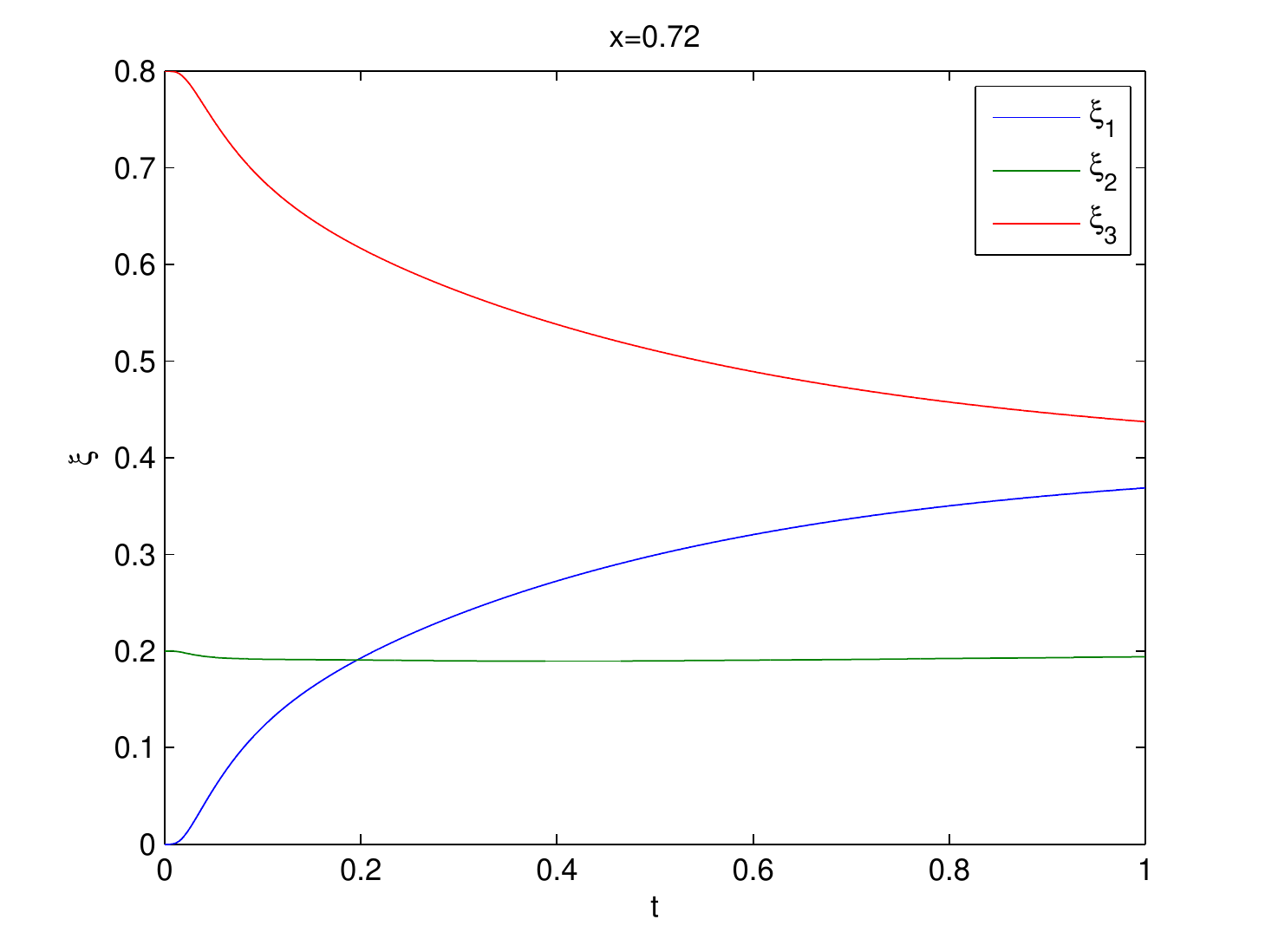} 
\end{center}
\caption{\label{hydrogen_2_1} The upper left figure presents the concentration of $x_1$ at spatial point $0.72$, the upper right is the result in the space time region, the lower left figure presents the the 3D plot of the second component and the lower right figure presents all of the components at spatial-point $0.72$.}
\end{figure}

The numerical solutions of the three hydrogen plasma in experiment 2 with the
uphill diffusion \ref{hydrogen_2_2}.
\begin{figure}[ht]
\begin{center}  
\includegraphics[width=5.0cm,angle=-0]{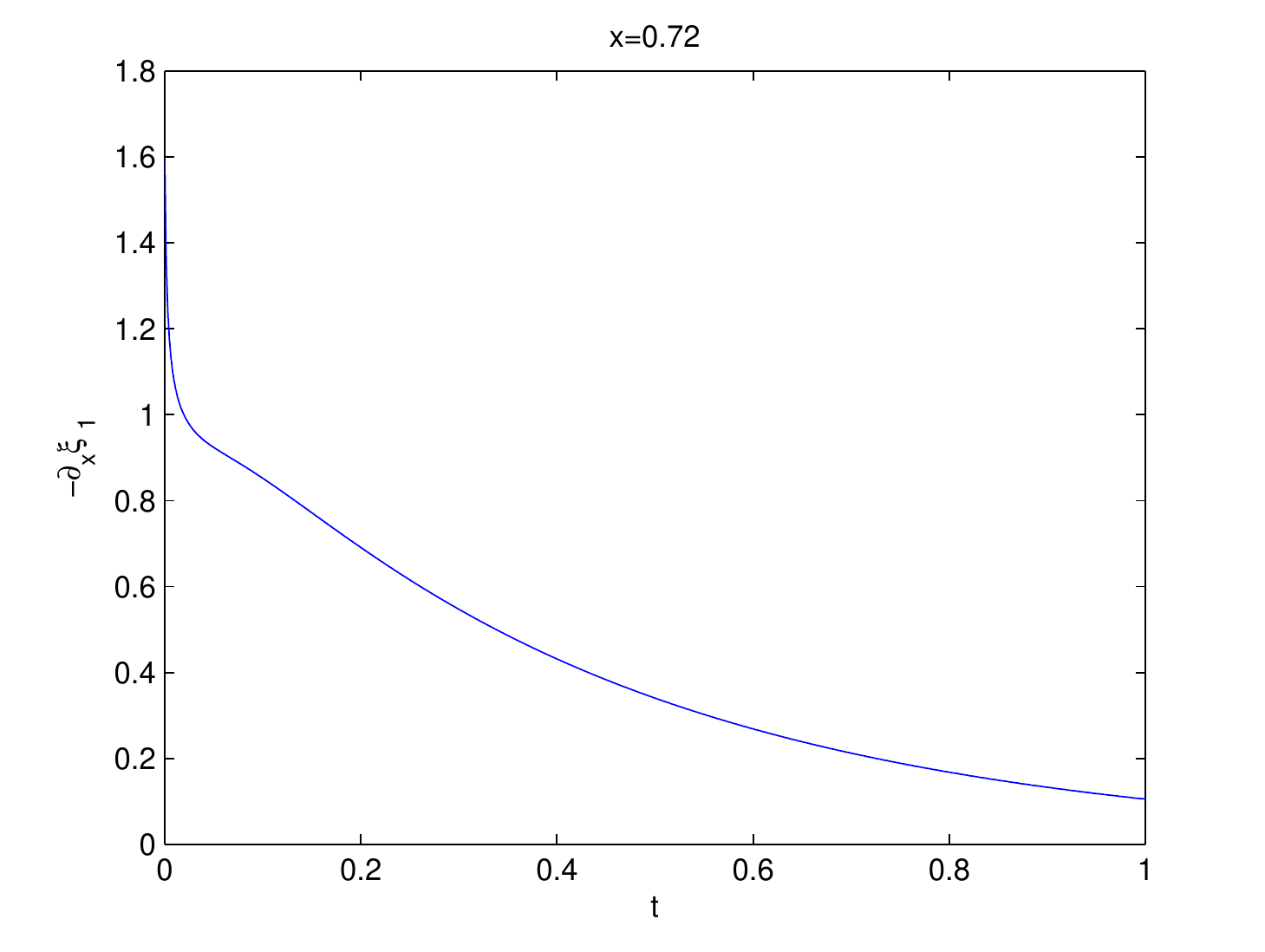} 
\includegraphics[width=5.0cm,angle=-0]{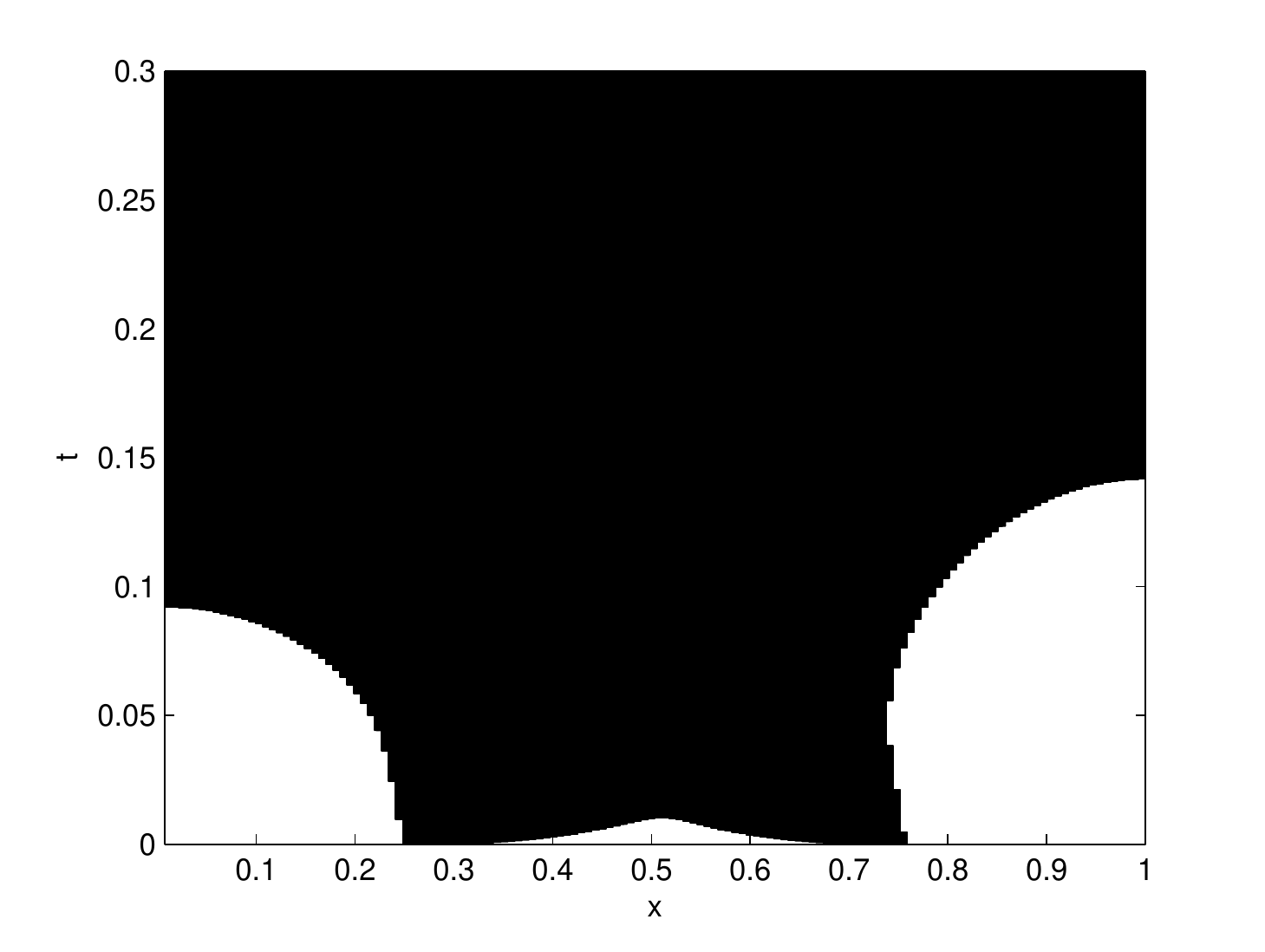} \\
\includegraphics[width=5.0cm,angle=-0]{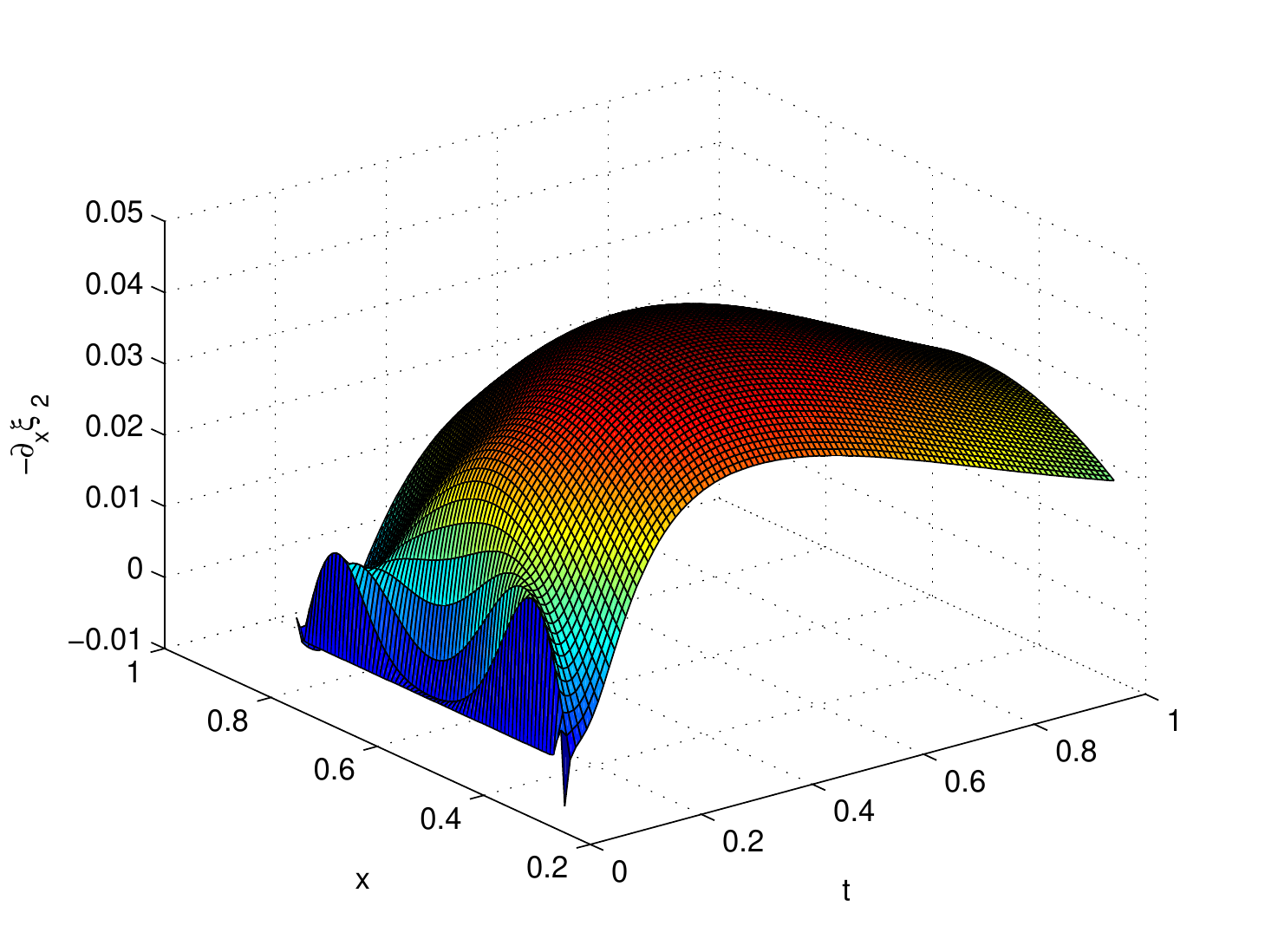} 
\includegraphics[width=5.0cm,angle=-0]{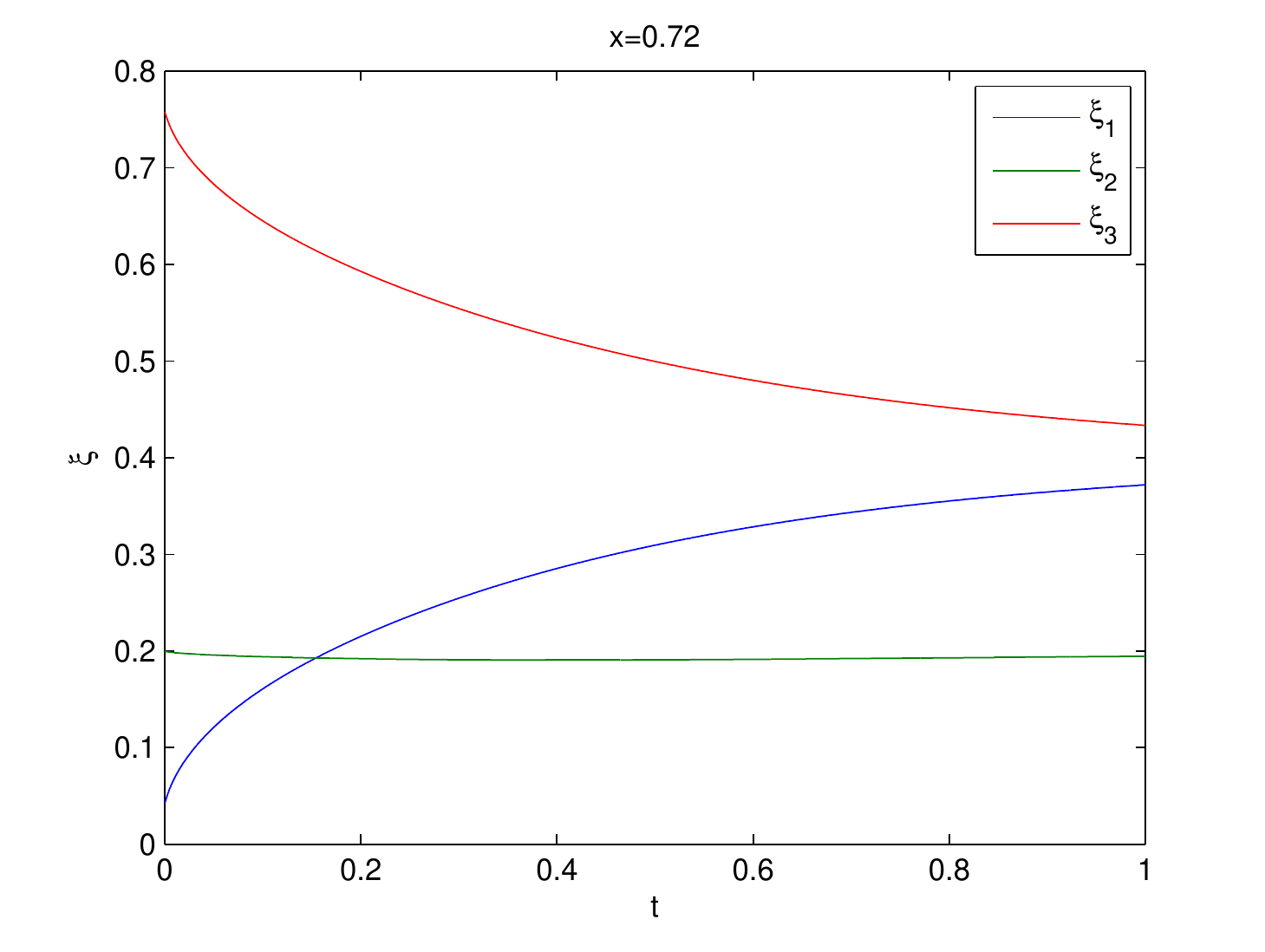} 
\end{center}
\caption{\label{hydrogen_2_2} The upper left figure presents the concentration of $x_1$ at spatial point $0.72$, the upper right is the result in the space time region, the lower left figure presents the the 3D plot of the second component and the lower right figure presents all of the components at spatial-point $0.72$.}
\end{figure}

\item Example 3:

 $\lambda_{11} = - 4.276 \; 10^{-1}$, 
$\lambda_{21} = \lambda_{31} = - \frac{\lambda_{11}}{2}$ , \\
$\lambda_{22} = - 2.082 \; 10^{-2} $, $\lambda_{12} = \lambda_{23} = - \frac{\lambda_{22}}{2}$ , \\
$\lambda_{33} = - 4.276 \; 10^{-2}$, $\lambda_{31} = \lambda_{32} =  - \frac{\lambda_{33}}{2}$

The numerical solutions of the three hydrogen plasma in experiment 3 with the
asymptotic diffusion \ref{hydrogen_3_1}.
\begin{figure}[ht]
\begin{center}  
\includegraphics[width=5.0cm,angle=-0]{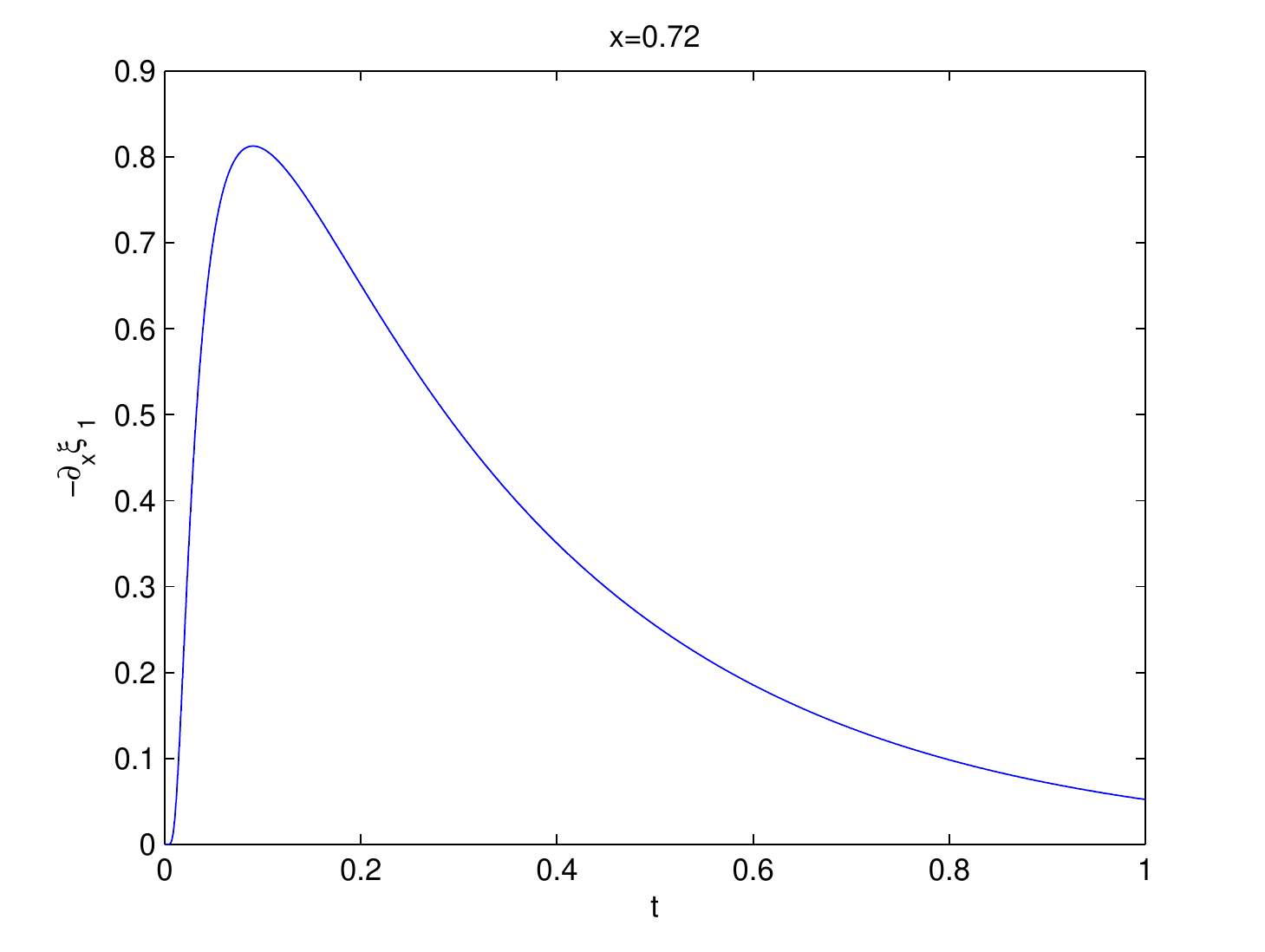} 
\includegraphics[width=5.0cm,angle=-0]{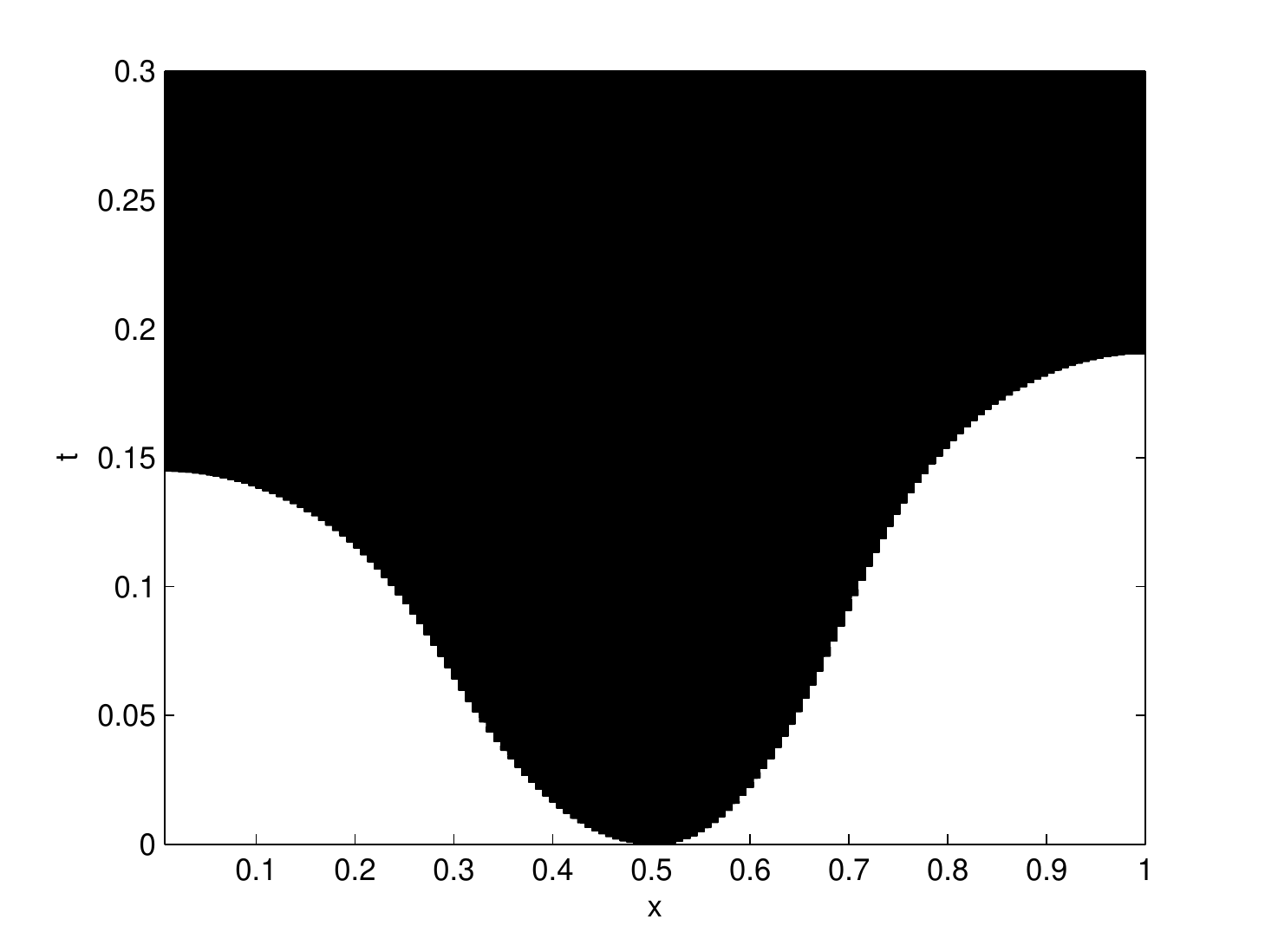} \\
\includegraphics[width=5.0cm,angle=-0]{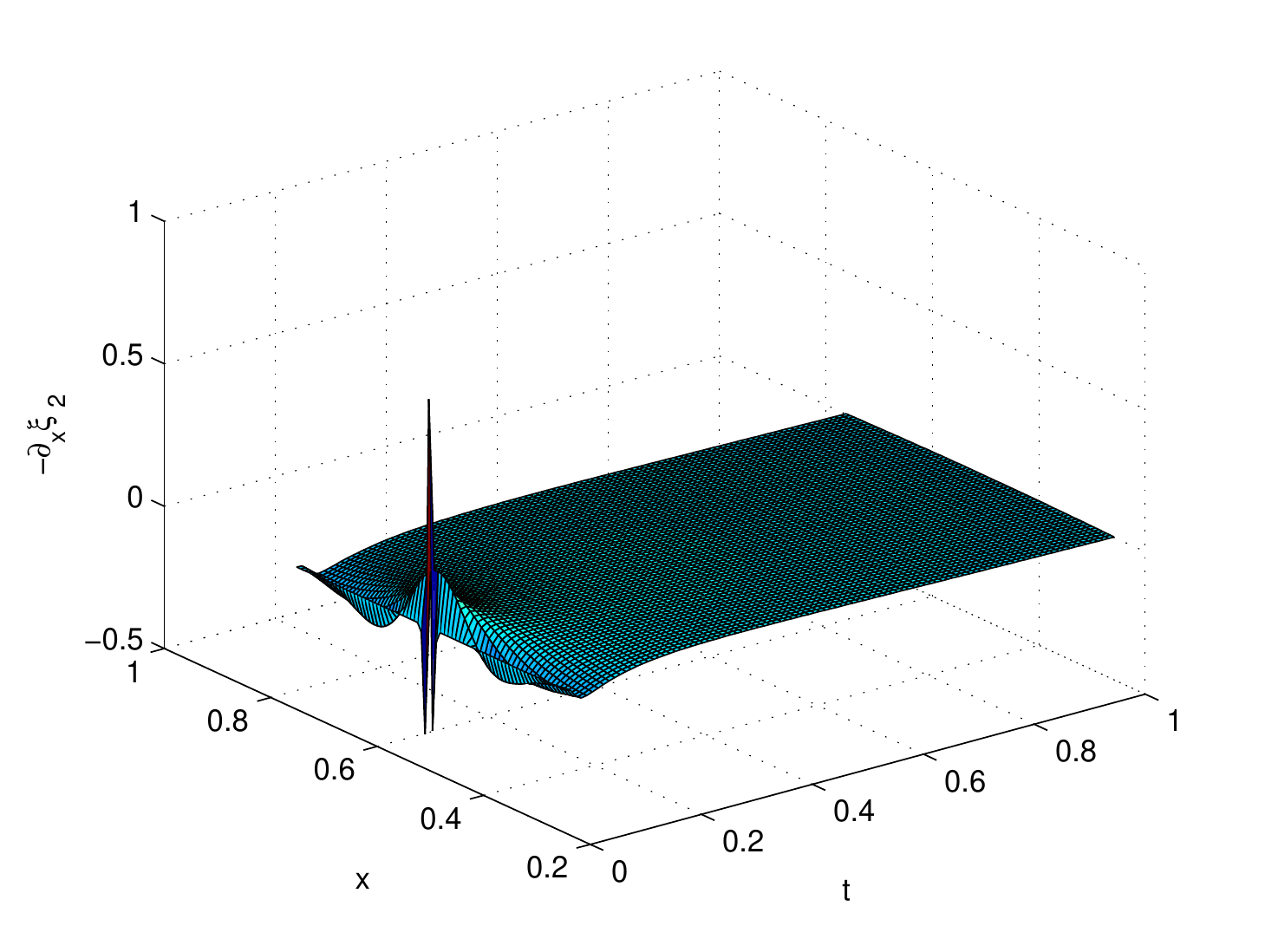} 
\includegraphics[width=5.0cm,angle=-0]{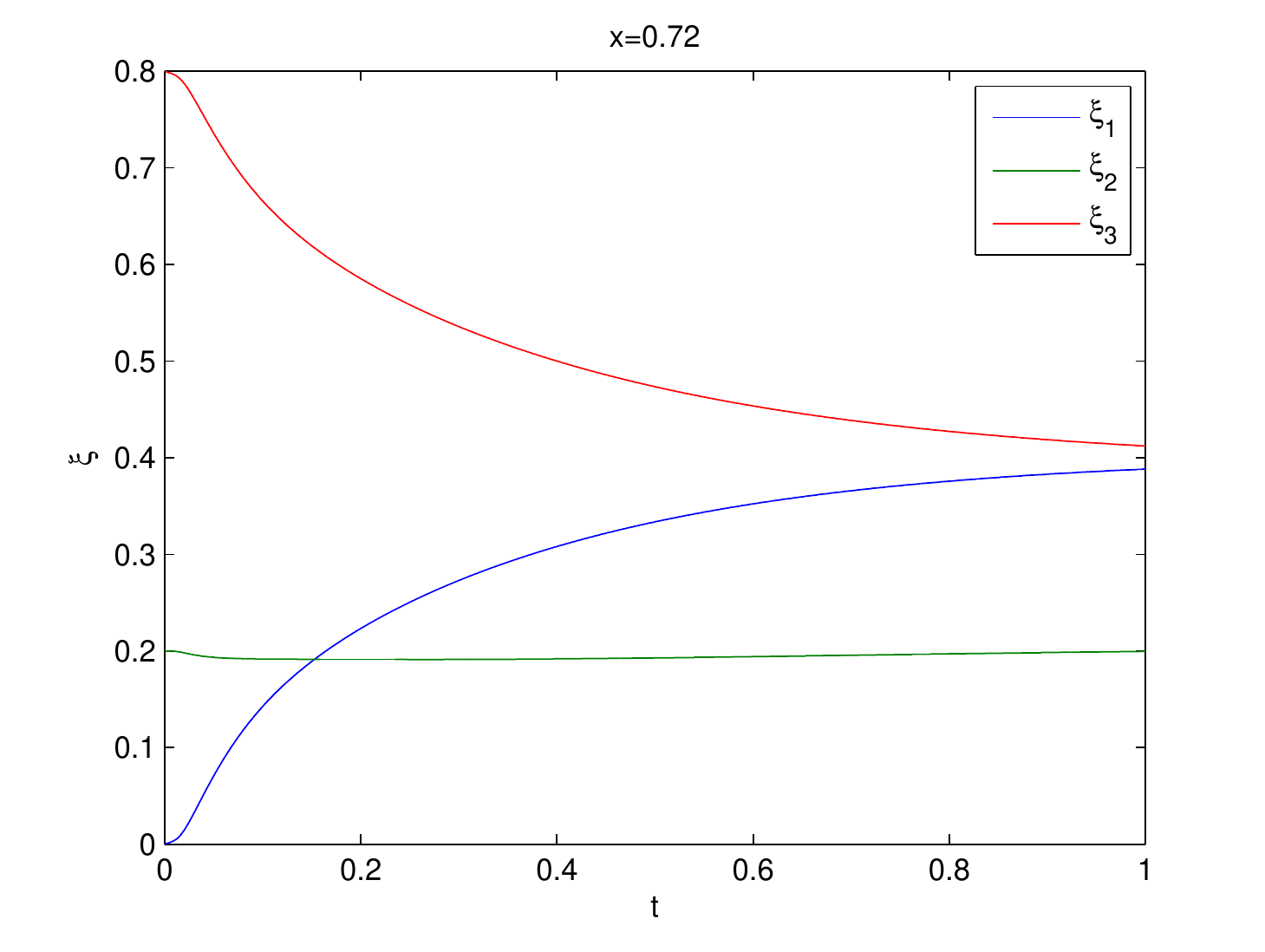} 
\end{center}
\caption{\label{hydrogen_3_1} The upper left figure presents the concentration of $x_1$ at spatial point $0.72$, the upper right is the result in the space time region, the lower left figure presents the the 3D plot of the second component and the lower right figure presents all of the components at spatial-point $0.72$.}
\end{figure}

The numerical solutions of the three hydrogen plasma in experiment 3 with the
uphill diffusion \ref{hydrogen_3_2}.
\begin{figure}[ht]
\begin{center}  
\includegraphics[width=5.0cm,angle=-0]{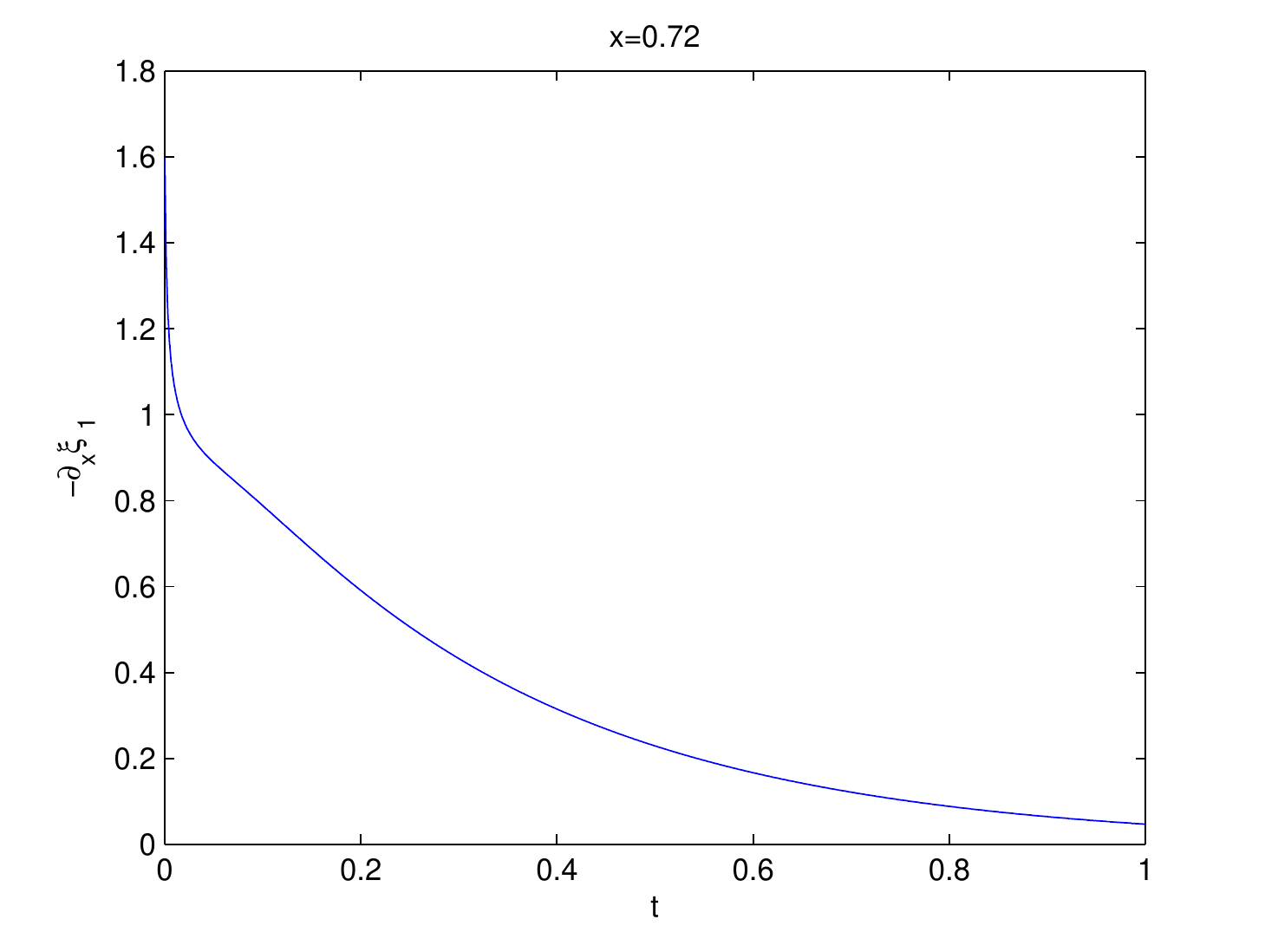} 
\includegraphics[width=5.0cm,angle=-0]{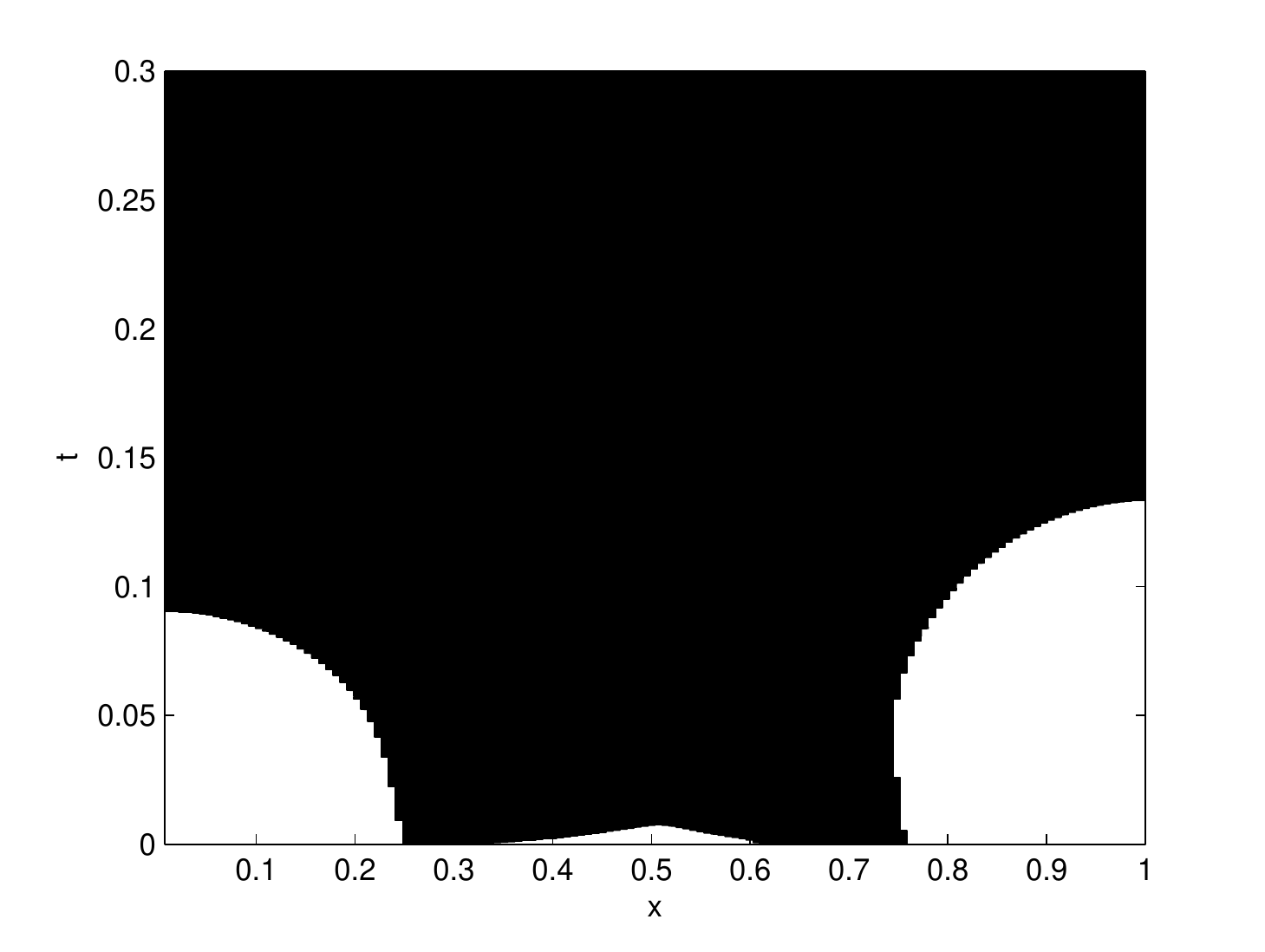} \\
\includegraphics[width=5.0cm,angle=-0]{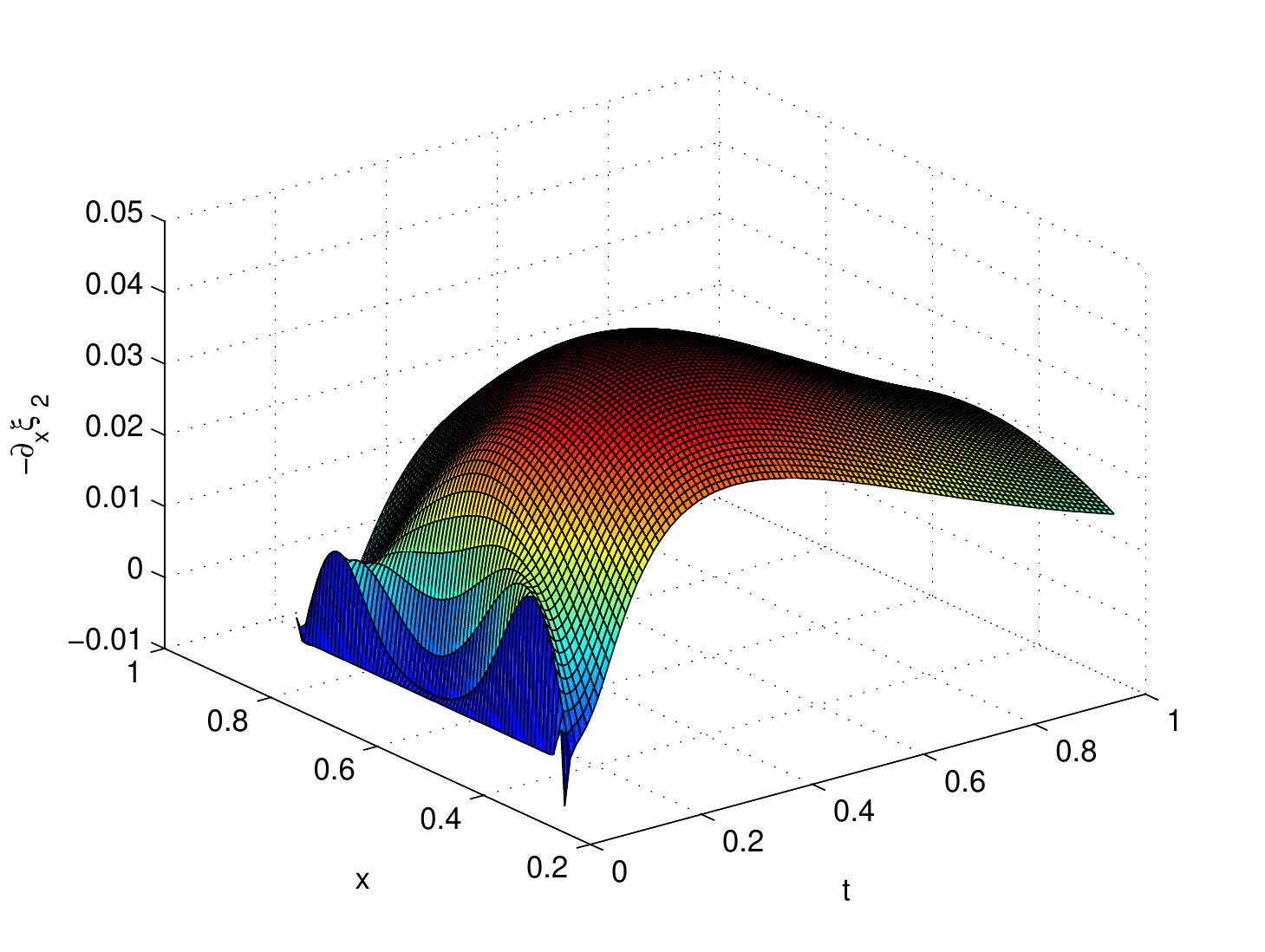} 
\includegraphics[width=5.0cm,angle=-0]{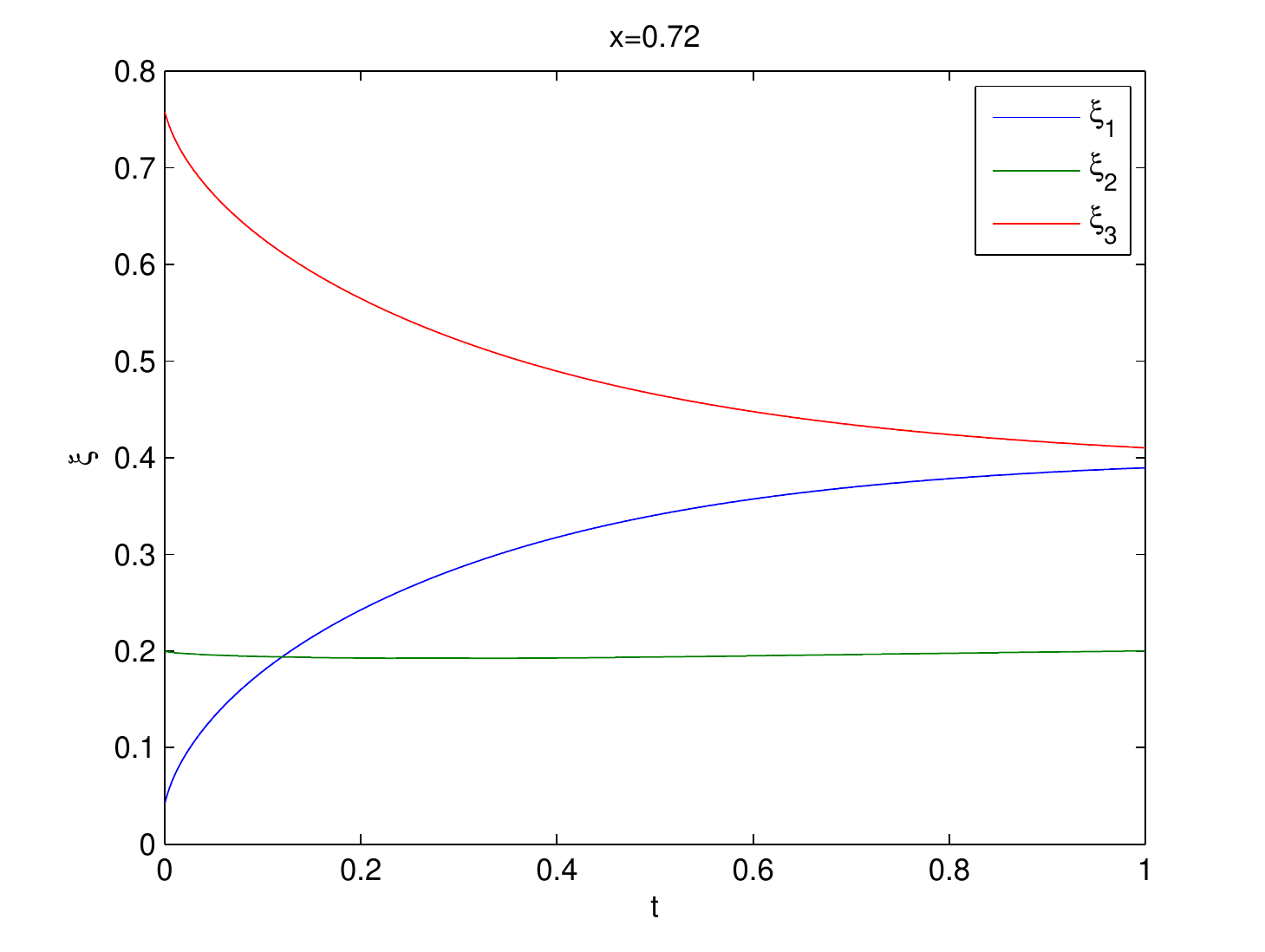} 
\end{center}
\caption{\label{hydrogen_3_2} The upper left figure presents the concentration of $x_1$ at spatial point $0.72$, the upper right is the result in the space time region, the lower left figure presents the the 3D plot of the second component and the lower right figure presents all of the components at spatial-point $0.72$.}
\end{figure}

\end{itemize}

In the following, we compare the different splitting methods based
on the first example with the uphill diffusion.

We deal with a CFL-grid means and we compare the results to the
optimal time- and spatial-grid size.
Based on this comparison, we are able to find the convergence-tableau
for the explicit methods.

We apply the following errors:

\begin{itemize}

\item Scalar for each $\xi_1, \xi_2$:
\begin{itemize}
\item Comparison in Time:
\begin{eqnarray}
 err_{i, j}(x, \Delta t) = \sum_{n=1}^N \Delta t \; | \xi_{i, ref}(x, t_n) - \xi_{i, j}(x, t_n) | ,
\end{eqnarray}
where we have the component index $i = 1, 2$ and the method index 
$j = \{ AB, ABA, iter \}$ and $x$ is given as an important spatial point, such as $x = 0.72$. Furthermore, $u_{ref}$ is a reference solution, such as with very small $\Delta t_{ref}$, and $u_j$ is the numerical solution of the method $j$ with $\Delta t = \{ \Delta t_{coarse}, \Delta t_{coarse}/2,  \Delta t_{coarse}/4, \Delta t_{coarse}/8 \}$ and the finest time-step is $\Delta t_{coarse}/16$. In space, we compare to the coarsest grid, which means that we interpolate the
finer space solutions to the coarsest grid.
\item Comparison in time and space:
\begin{eqnarray}
&& err_{i, j}(\Delta t) = \\
&& = \sum_{k=1}^{J_{coarse}} \Delta x_{coarse} \; \sum_{n=1}^N \Delta t \; | \xi_{i, ref}(x_k(\Delta t_{coarse}/16), t_n) - \xi_{i, j}(x_k(\Delta t), t_n) | , \nonumber
\end{eqnarray}
where we have the component index $i = 1, 2$ and the method index $j = \{ AB, ABA, iter \}$, and $T$ is given as an important time point, such as the end time-point $t = T$. Furthemore,r $\xi_{i, ref}$ is a reference solution, such as with very small $\Delta t_{ref}$, and $\xi_{i,j}$ is the numerical solution of the method $j$ with $\Delta t = \{ \Delta t_{coarse}, \Delta t_{coarse}/2,  \Delta t_{coarse}/4, \Delta t_{coarse}/8 \}$ and the finest time-step is $\Delta t_{coarse}/16$. In space, we compare to the coarsest grid, which means that we interpolate the
finer space solutions to the coarsest grid.
\end{itemize}

\item Vectorial for ${\bf \xi} = (\xi_1, \xi_2, \xi_3)^t$:
\begin{itemize}
\item Comparison in Time:
\begin{eqnarray}
 err_{j}(x, \Delta t) = \sum_{n=1}^N \Delta t \; ( \sum_{i=1}^3 | \xi_{i, ref}(x, t_n) - \xi_{i, j}(x, t_n) | ) ,
\end{eqnarray}
where the method index 
$j = \{ AB, ABA, iter \}$ and $x$ is given as an important spatial point, such as $x = 0.72$. Furthermore, $u_{ref}$ is a reference solution, such as with very small $\Delta t_{ref}$, and $u_j$ is the numerical solution of the method $j$ with $\Delta t = \{ \Delta t_{coarse}, \Delta t_{coarse}/2,  \Delta t_{coarse}/4, \Delta t_{coarse}/8 \}$ and the finest time-step is $\Delta t_{coarse}/16$. In space, we compare to the coarsest grid, which means that we interpolate the
finer space solutions to the coarsest grid.
\item Comparison in time and space:
\begin{eqnarray}
&& err_{j}(\Delta t) = \\
&& = \sum_{k=1}^{J_{coarse}} \Delta x_{coarse} \sum_{n=1}^N \Delta t ( \sum_{i=1}^3 | \xi_{i, ref}(x_k(\Delta t_{coarse}/16), t_n) - \xi_{i, j}(x_k(\Delta t), t_n) | ) , \nonumber
\end{eqnarray}
where the method index $j = \{ AB, ABA, iter \}$ and $T$ is given as an important time point, such as the end time-point $t = T$. Further $u_{ref}$ is a reference solution, such as with very small $\Delta t_{ref}$, and $u_j$ is the numerical solution of the method $j$ with $\Delta t = \{ \Delta t_{coarse}, \Delta t_{coarse}/2,  \Delta t_{coarse}/4, \Delta t_{coarse}/8 \}$ and the finest time-step is $\Delta t_{coarse}/16$. In space, we compare to the coarsest grid, which means that we interpolate the
finer space solutions to the coarsest grid.
\end{itemize}

\end{itemize}

Convergence-tableau for the different methods.

We have the following CFL-condition:
\begin{eqnarray}
\Delta t \le \frac{\Delta x^2}{2 D_{max}} \approx \Delta x^2 , 
\end{eqnarray}
where we have $2 D_{max} \approx 1$.

We write in the notation of the grid-points:
\begin{eqnarray}
J^2 \le N ,
\end{eqnarray}
where $J$ are the number of spatial grid-points and $N$ is the number
of time-points.

\begin{figure}[ht]
\begin{center}  
\includegraphics[width=5.0cm,angle=-0]{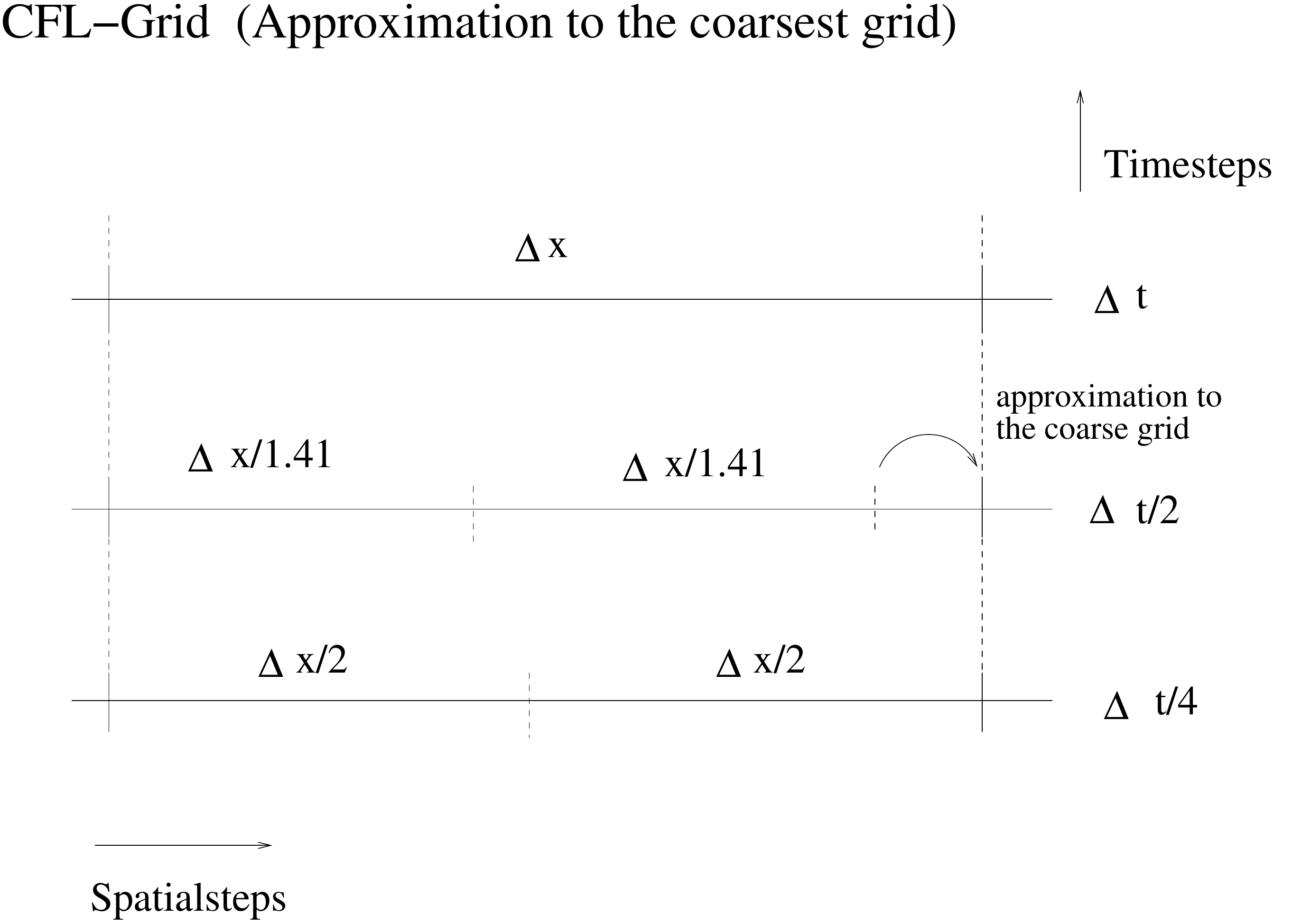} 
\end{center}
\caption{\label{fig_cfl} The optimal spatial-grid (CFL-condition) to the time-steps.}
\end{figure}

We have the following resolutions in Table \ref{table_1_ref} and Figure \ref{fig_cfl}.
\begin{table}[h]
\begin{center}
\begin{tabular}{|c|c|c|}
\hline 
Method & Spatial-Points $J$ & Time-points $N$ \\
\hline 
 iter3 & 190 & 80000\\
 (reference solution) & & \\
\hline
\hline
iter3, iter2,  & 140 & 40000\\
AB, ABA  	&&	 \\
\hline
iter3, iter2,  & 100 & 20000\\
AB, ABA  	&&	 \\
\hline
iter3, iter2,  & 70 & 10000\\
AB, ABA  	&&	 \\
\hline
iter3, iter2,  & 50 & 5000 \\
AB, ABA  	&&	 \\
  \hline
\end{tabular}
\caption{\label{table_1_ref} The spatial- and time-grid-points related to the
reference solution}
\end{center}
\end{table}

To compare the values only on the coarsest CFL-grid, we have to apply the
following approximation:

\begin{eqnarray}
&& x_k(\Delta t) , k = 0, \ldots, 50 , \nonumber \\
&&  \mbox{coarsest spatial grid with time-step} \; \Delta t , \\
&& x_k(\Delta t/2) , k(\Delta t/2) = \lfloor k \; \sqrt{2} \rfloor  = 1, 2, 4, \ldots, 70 ,  \nonumber  \\
&& \mbox{next finer grid with time-step} \; \Delta t/2 , \\
&& x_k(\Delta t/4) , k(\Delta t/4) = k \; 2 = 2, 4, 6, \ldots, 100 ,  \nonumber  \\
&&  \mbox{next finer grid with time-step} \; \Delta t/4 , \\
&& x_k(\Delta t/8) , k(\Delta t/8) = \lfloor k \; 2 \sqrt{2} \rfloor  = 2, 5, 8, \ldots, 140 ,  \nonumber  \\
&& \mbox{next finer grid with time-step} \; \Delta t/8 , \\
&& x_k(\Delta t/16) , k(\Delta t/16) = \lfloor k \; 4 \rfloor  = 4, 8, 12, \ldots, 200 ,  \nonumber  \\
&& \mbox{finest grid with time-step} \; \Delta t/16 ,
\end{eqnarray}
where $\lfloor x \rfloor = \max \{k \in \mathbb {Z} | k \le x  \}$.

\begin{remark}
We have the following computational times for the Picard's methods
in table \ref{table_1_hydro}.
\begin{table}[h]
\begin{center}
\begin{tabular}{|c|c|}
\hline 
Example & Computational Time   [sec]   \\
\hline 
\hline 
 1	&	6.8736e+03 \\
 2	&	7.3985e+03 \\
 3	&	8.8402e+03 \\
  \hline
\end{tabular}
\caption{\label{table_1_hydro} The computational time of the  three experiments with different 
reaction parameters with $N_{spatial} = 140$ number of spatial discretisation points, $N_{end} = 80000$ number of time-steps}
\end{center}
\end{table}

The convergence results are given in Figure \ref{hydrogen_conv}.
\begin{figure}[ht]
\begin{center}  
\includegraphics[width=5.0cm,angle=-0]{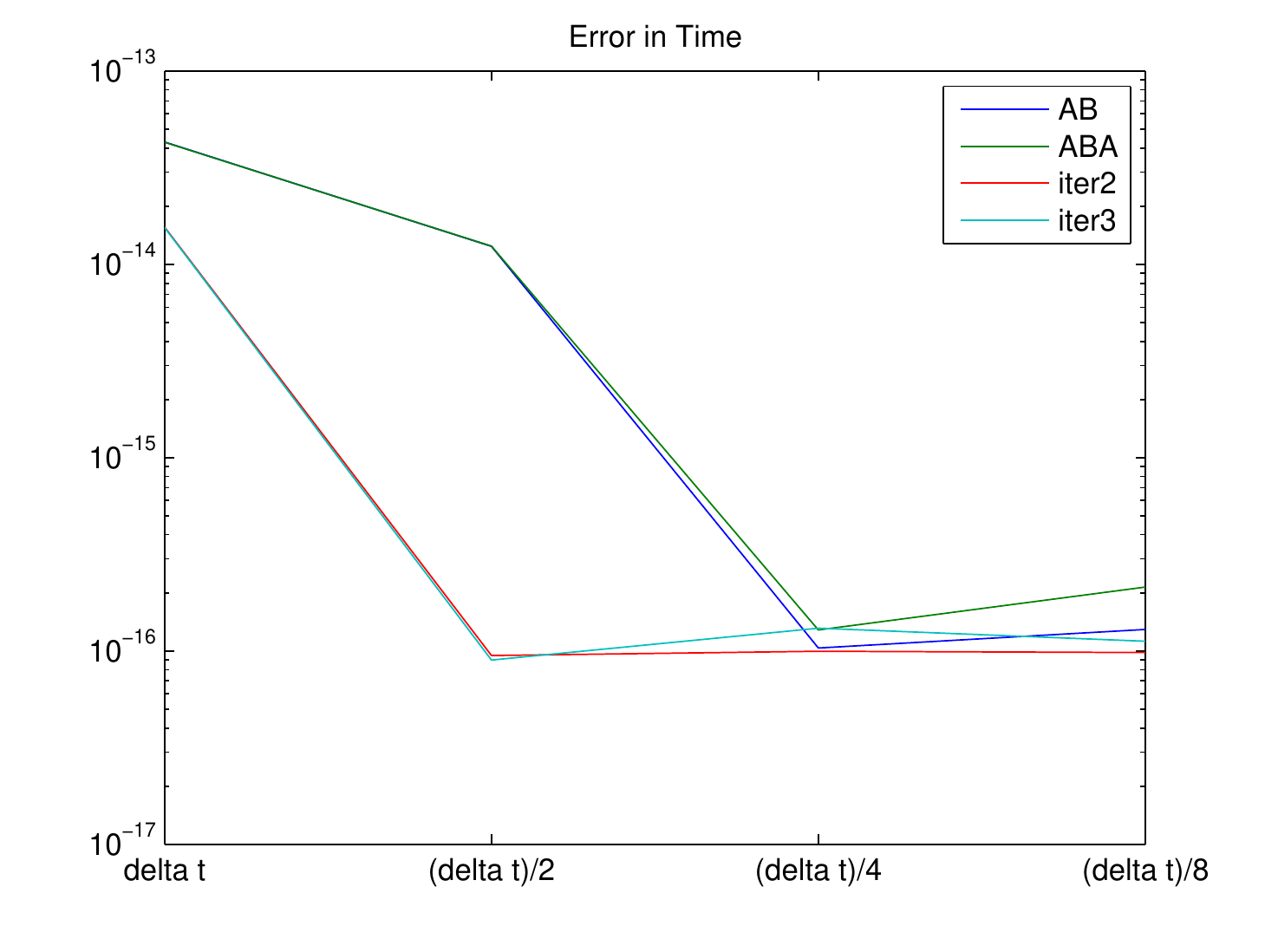} 
\end{center}
\caption{\label{hydrogen_conv} The convergence of the different methods are computed over the full domain and with different time-steps.}
\end{figure}

We have the following computational times for the Picard's methods
in table \ref{table_2_hydro}.
\begin{table}[h]
\begin{center}
\begin{tabular}{|c|c|}
\hline 
Example & Computational Time   [sec]   \\
\hline 
\hline 
AB &  	165.0895 [sec] \\
ABA &	262.3892 [sec]  \\
iter2 &	404.9536 [sec]  \\
iter3 &	578.3206 [sec] \\
  \hline
\end{tabular}
\caption{\label{table_2_hydro} The computational time of the different 
methods with $N_{spatial} = 140$ number of spatial discretisation points, $N_{end} = 80000$ number of time-steps}
\end{center}
\end{table}

Here, we see the additional work of the iterative methods.

\end{remark}
We obtain optimal solutions for the iterative methods, while we could
extend the time-step. For more detailed computations and smaller
time-steps, the noniterative splitting methods are more effective, while
we could obtain at least a second order approach, see also \cite{geiser2016_1}.

%AU: Is this remark (i.e. Remark 3) missing?
\begin{remark}

\end{remark}

\section{Conclusions and Discussions }
\label{concl}

We present the coupled model for a multi-component transport model
for reactive plasma. The nonlinear partial differential equations
are solved with iterative methods and a combination of splitting 
approaches. The numerical algorithms are presented and their numerical
convergences are shown. Although iterative splitting 
methods are more time-consuming, they are more
accurate than noniterative splitting approaches.
The benefits of noniterative methods when we apply explicit 
schemes include fast computation time 
and good resolution of space and time space. The implicit behavior of iterative methods allows larger
time-steps to be used and they could accelerate the solver process.
In the future we aim to study the numerical analysis of the different combined schemes
and we will simulate more delicate multicomponent models.

\bibliographystyle{plain}

\end{document}